%% file: A_Reissner-Mindlin_plate_formulation_using_symmetric_Hu-Zhang_elements_via_polytopal_transformations.tex
\documentclass[10pt]{article}
\usepackage{amsfonts,color}
\usepackage{amssymb}
\usepackage{footnote}
\usepackage{mathtools}
\usepackage{amsthm,wasysym}
\usepackage[english]{babel}
\usepackage{bbm}
\usepackage{colonequals}
\usepackage{setspace}
\usepackage{titlesec}
\usepackage{float}

\usepackage{graphicx}
\usepackage{amsfonts}
\usepackage{hyperref}
\usepackage{subcaption}
\usepackage{amsmath}
\usepackage{nicefrac}
\usepackage{gensymb}
\usepackage{enumerate}
\usepackage[nameinlink,capitalize]{cleveref}
\usepackage{chngcntr}
\usepackage{multirow}

\usepackage{cite}
\usepackage[nottoc,numbib]{tocbibind}

\usepackage{pgf,tikz,pgfplots}
\pgfplotsset{compat=1.14}
\usepackage{mathrsfs}
\usetikzlibrary{arrows}
\usetikzlibrary{arrows.meta}

\DeclareMathOperator{\at}{\bigg\vert}
\newcommand{\vb}[1]{\mathbf{#1}}
\newcommand{\bm}[1]{\boldsymbol{#1}}

\DeclareMathOperator{\sym}{\mathrm{sym}}
\DeclareMathOperator{\spa}{\mathrm{span}}

\DeclareMathOperator{\id}{\mathrm{id}}
\DeclareMathOperator{\cof}{\mathrm{cof}}

\newcommand{\one}{\bm{\mathbbm{1}}}
\newcommand{\con}[2]{\langle {#1} , \, {#2} \rangle}
\newcommand{\norm}[1]{\| {#1} \|}
\newcommand{\trnp}{\mathrm{tr}_n^\parallel}

\DeclareMathOperator{\airy}{\mathrm{airy}}

\newcommand{\dd}{\mathrm{d}}
\newcommand{\D}{\mathrm{D}}
\DeclareMathOperator{\di}{\mathrm{div}}
\DeclareMathOperator{\Di}{\mathrm{Div}}

\DeclareMathOperator{\curl}{\mathrm{curl}}

\newcommand{\SO}{\mathrm{SO}}

\newcommand{\Sym}{\mathrm{Sym}}
\newcommand{\Ned}{\mathcal{N}_{I}}
\newcommand{\AW}{\mathcal{AW}}
\newcommand{\HHJ}{\mathcal{HHJ}}
\newcommand{\HZ}{\mathcal{HZ}}
\newcommand{\RT}{\mathcal{RT}}
\newcommand{\BDM}{\mathcal{BDM}}

\newcommand{\Lag}{\mathcal{L}}

\newcommand{\Po}{\mathit{P}}
\newcommand{\Le}{{\mathit{L}^2}}

\newcommand{\Hone}{\mathit{H}^1}
\newcommand{\Honez}{\mathit{H}_0^1}
\newcommand{\HsD}[1]{\mathit{H}^\mathrm{sym}(\mathrm{Div}{#1})}
\newcommand{\Hd}[1]{\mathit{H}(\mathrm{div}{#1})}

\newcommand{\Hc}[1]{\mathit{H}(\mathrm{curl}{#1})}

\newcommand{\HdD}[1]{\mathit{H}(\di\mathrm{Div}{#1})}

\newcommand{\body}{V}
\newcommand{\surf}{A}
\newcommand{\curv}{s}

\renewcommand{\O}{\mathcal{O}}
\newcommand{\R}{\mathbb{R}}
\newcommand{\U}{\mathit{U}}
\newcommand{\X}{\mathit{X}}
\newcommand{\Y}{\mathit{Y}}
\newcommand{\Z}{\mathit{Z}}
\newcommand{\C}{\mathit{C}}
\renewcommand{\H}{\mathit{H}}
\newcommand{\tem}{\mathcal{T}}
\newcommand{\ver}{\mathcal{V}}
\newcommand{\edge}{\mathcal{E}}

\newcommand{\cell}{\mathcal{C}}

\newcommand{\Cm}{\mathbb{C}}
\newcommand{\J}{\mathbb{J}}
\newcommand{\A}{\mathbb{A}}

\newcommand{\ksmu}{k_s \mu}

\newcommand{\ud}{\vb{u}}

\allowdisplaybreaks[1]

\makeindex

\newtheoremstyle{break}
{\topsep}{\topsep}%
{\itshape}{}%
{\bfseries}{}%
{\newline}{}%
\theoremstyle{break}

\newtheorem{theorem}{Theorem}

\newtheorem{remark}{Remark}

\newtheorem{definition}{Definition}

\counterwithin{theorem}{section}
\counterwithin{lemma}{section}
\counterwithin{remark}{section}
\counterwithin{observation}{section}
\counterwithin{definition}{section}

\counterwithin{figure}{section}
\counterwithin{equation}{section}

\makeatletter
\let\@fnsymbol\@arabic
\makeatother

\crefname{Problem}{Problem.}{Problem.}

\title{A Reissner-Mindlin plate formulation using symmetric Hu-Zhang elements via polytopal transformations}
\author{\normalsize{Adam Sky}\thanks{Corresponding author: Adam Sky, Institute of Computational Engineering and Sciences, Department of Engineering, Faculty of Science, Technology and Medicine, University of Luxembourg, 6 Avenue de la Fonte, L-4362 Esch-sur-Alzette, Luxembourg, email: adam.sky@uni.lu}
    , \quad
	\normalsize{Michael Neunteufel}\thanks{Michael Neunteufel, Institute of Analysis and Scientific Computing, Technische Universit\"at Wien, Wiedner Hauptstr. 8-10 , 1040 Wien, Austria, email: michael.neunteufel@tuwien.ac.at}
    , \quad
	\normalsize{Jack S. Hale}\thanks{Jack S. Hale, Institute of Computational Engineering and Sciences, Department of Engineering, Faculty of Science, Technology and Medicine, University of Luxembourg, 6 Avenue de la Fonte, L-4362 Esch-sur-Alzette, Luxembourg, email: jack.hale@uni.lu}
    \quad 
    \normalsize{and} \quad
	\normalsize{Andreas Zilian}\thanks{Andreas Zilian, Institute of Computational Engineering and Sciences, Department of Engineering, Faculty of Science, Technology and Medicine, University of Luxembourg, 6 Avenue de la Fonte, L-4362 Esch-sur-Alzette, Luxembourg, email: andreas.zilian@uni.lu}
}

\begin{document}

\maketitle

\begin{abstract}
In this work we develop new finite element discretisations of the
shear-deformable Reissner--Mindlin plate problem based on the Hellinger-Reissner principle of symmetric stresses. Specifically, we use conforming Hu-Zhang elements to discretise the bending moments in the space of symmetric square integrable fields with a square integrable divergence $\bm{M} \in \HZ \subset \HsD{}$. The latter results in highly accurate approximations of the bending moments $\bm{M}$ and in the rotation field being in the discontinuous Lebesgue space $\bm{\phi} \in [\Le]^2$, such that the Kirchhoff-Love constraint can be satisfied for $t \to 0$. In order to preserve optimal convergence rates across all variables for the case $t \to 0$, we present an extension of the formulation using Raviart-Thomas elements for the shear stress $\vb{q} \in \RT \subset \Hd{}$.

We prove existence and uniqueness in the continuous setting and rely on exact complexes for inheritance of well-posedness in the discrete setting.

This work introduces an efficient construction of the Hu-Zhang base functions on the reference element via the polytopal template methodology and Legendre polynomials, making it applicable to hp-FEM. The base functions on the reference element are then mapped to the physical element using novel polytopal transformations, which are suitable also for curved geometries.

The robustness of the formulations and the construction of the Hu-Zhang element are tested for shear-locking, curved geometries and an L-shaped domain with a singularity in the bending moments $\bm{M}$. Further, we compare the performance of the novel formulations with the primal-, MITC- and recently introduced TDNNS methods.         
\\
\vspace*{0.25cm}
\\
\end{abstract}

\section{Introduction}
The subject of this paper is a new finite element discretisation of the
shear-deformable Reissner--Mindlin plate problem. The Reissner-Mindlin plate
problem describes the elastic behaviour of a body with thickness $t$ far smaller
than in its planar dimensions. It is well known that a naive discretisation of
the Reissner--Mindlin problem will `lock' as the thickness of the plate
approaches zero, which enforces the Kirchhoff-Love constraint. This loss of
robustness in the thickness means that the numerical method may converge
sub-optimally, or even produce a completely incorrect solution.

Over the past decades the problem of shear-locking in the finite element
discretisation of the Reissner--Mindlin plate problem has received a
significant amount of attention in the literature. We can point to the Mixed
Interpolation of Tensorial Components (MITC) approaches
\cite{bathe_mitc7_1989,hale_simple_2018}, Partial Selective Reduced Integration
(PSRI) approaches~\cite{chinosi_numerical_1995}, the Falk-Tu
element~\cite{falk_locking-free_2000}, the Discrete Shear Gap (DSG)
method~\cite{falk_locking-free_2000}, and the smoothed finite element approach \cite{Nguyen-Xuan2010,NGUYENXUAN20081184} as examples. We refer the reader to the
review~\cite{falk_finite_2008} for a full treatment. A common thread in both
the implementation and numerical analysis of these methods is the use of a
mixed formulation where a dual quantity, such as the shear stress, is
approximated in addition to the usual primal rotation and displacement
variables. For example, in the MITC approach the rotations are interpolated
from the discrete $[\Hone]^2$-conforming space into the larger discrete
$\Hc{}$-conforming space \cite{hale_simple_2018}, which is compatible with the
Kirchhoff--Love condition, compare \cref{sec:locking}. Typically, the
$\Hc{}$-conforming N\'ed\'elec spaces are employed
\cite{ndlec_new_1986,nedelec_mixed_1980}.

As it is the underlying approach we take in this paper, we turn our focus to
the smaller number of numerical methods for the Reissner-Mindlin problem that
use a mixed formulation involving the bending moments
\cite{daVeiga,Guo2013,pechstein_tdnns_2017}. The bending moments lie naturally
in the Hilbert space of symmetric tensors with row-wise square-integrable
divergence $\bm{M} \in \HsD{}$. Conforming finite element discretisations of
this space have been studied in the context of the two and three-dimensional
Hellinger--Reissner elasticity problem \cite{Lederer} where the Cauchy stress
must be discretised. Notable examples of conforming elements include the
Arnold--Winther ($\AW$) \cite{arnold_mixed_2002,arnold_finite_2008} and
Hu--Zhang ($\HZ$) finite elements
\cite{hu_family_2014,hu_family_2015,hu_finite_2016}. The computer
implementation of conforming discretisations of $\HsD{}$ also remains an
impediment to their widespread use. Their implementation is complicated because
the reference cell basis functions do not map `straightforwardly' to the
physical cell. Recent advancements in transformation theory
\cite{aznaran_transformations_2021,kirby_general_2018} and open source finite
element technology~\cite{kirby_code_2019} have bought the $\AW$ elements closer
to practical use for users of the open source Firedrake finite element
solver~\cite{rathgeber_firedrake_2016}. However, the question of how to map
this element onto curved geometries remains open. We also remark that the
lowest-order conforming three-dimensional analogues of the $\AW$ and $\HZ$
spaces~\cite{aznaran_transformations_2021} have large local dimensions, $162$
and $210$ respectively, and consequently are not used widely in practical
computations. As the Reissner-Mindlin problem is two-dimensional we avoid this
practical issue in our work.

To avoid having to deal with the difficulties of discretising $\HsD{}$ it is
also possible to employ methods where each row of the stress tensor is
discretised in $\Hd{}$ and the necessary symmetry of stress tensor is
subsequently enforced using, e.g., a Lagrange multiplier. Again, this technique
was first used in the context of the Hellinger--Reissner formulation of
elasticity~\cite{arnold_peers_1984,arnold_mixed_2007}, and was then used in the
Reissner--Mindlin problem in \cite{daVeiga}.

Most recently, the Tangential Displacement Normal-Normal Stress (TDNNS)
method, which was also conceived in the context of the standard elasticity
problem \cite{pechstein_anisotropic_2012,pechstein_analysis_2018}, has been
employed to alleviate locking in the Reissner-Mindlin problem
\cite{pechstein_tdnns_2017}. This method can be seen as a natural extension of
the ideas of Hellan, Herrmann and Johnson for the Kirchoff--Love plate problem
to the Reissner--Mindlin setting. The TDNNS formulation discretises the
rotations directly in the $\Hc{}$-space ensuring their exact compatibility with
the gradient of the transverse displacement $w \in \H^1$ which is given by $\nabla w
\in \Hc{}$, as per the exact de Rham sequence \cite{pauly_hilbert_2022,arnold_complexes_2021}. This is made possible by defining the bending moments $\bm{M} \in
\HdD{}$ the Hilbert space of tensor-valued functions with square-integrable
normal-normal components. The use of this space is justified by relying on the
duality of $\Hc{}$ with $\H^{-1}(\di)$ and the relation $\Di \HdD{} \subseteq
\H^{-1}(\di)$. In the discrete construction the authors of
\cite{pechstein_tdnns_2017} employ conforming N\'ed\'elec elements for the
rotations $\bm{\phi} \in \Ned \subset \Hc{}$ and the `slightly' non-conforming
Hellan--Herrmann--Johnson elements
\cite{neunteufel_hellanherrmannjohnson_2019,arnold_hellan--herrmann--johnson_2020}
for the bending moments $\bm{M} \in \HHJ \nsubseteq \HdD{}$. The convergence of
the formulation is proven in both discrete \cite{pechstein_tdnns_2017} and
natural norms \cite{pechstein_analysis_2018}.

In this paper we propose formulations where the shear-strains \cite{Voss2020} are elevated to
$\vb{q} \in [\Le]^2$ respectively $\vb{q} \in \Hd{}$, by defining the rotations in $\bm{\phi} \in [\Le]^2$, thus
circumventing shear-locking due to incompatibility of the spaces. We note that
this level of regularity on the rotation is lower than both the standard
displacement approach $\bm{\phi} \in [\Hone]^2$, and the TDNNS approach where
$\bm{\phi} \in \Hc{}$. The advantage of this approach lies in the employment of
the bending moments in $\bm{M} \in \HsD{}$. This allows us to handle all
existence and uniqueness proofs in the context of exact Hilbert
complexes~\cite{arnold_finite_2008,hu_family_2014,arnold_complexes_2021,PaulyDeRham}, which automatically assert their
validity in the conforming discrete setting~\cite{Bra2013}. For the conforming
discretisation of the bending moments we rely on Hu-Zhang elements
\cite{hu_family_2014} of arbitrary order, such that the bending moments read
$\bm{M} \in \HZ \subset \HsD{}$.

For the computer implementation of our method we extend the polytopal template
methodology for basis function construction, introduced in
\cite{sky_polytopal_2022,sky_higher_2023}, to the $\HZ \subset \HsD{}$ setting. Further, we
present novel polytopal transformations of the base functions from the
reference to the physical element based on fourth-order tensors. 
We stress that this novel transformation approach allows for curved finite elements on the physical domain. 
In addition, our construction employs orthogonal Legendre polynomials, thus making it
appropriate to the hp-finite element method~\cite{VOS20105161,KOPP2022115575}. The
implementation is carried out in the open source finite element software NGSolve\footnote{www.ngsolve.org}
\cite{Sch2014,Sch1997} and is available as supplementary material to this
paper\footnote{https://github.com/Askys/NGSolve\_HuZhang\_Element}. 

This paper is structured as follows. First we recap the derivation of the Reissner-Mindlin plate model and its corresponding variational problem in the primal setting. Next we discuss the phenomena of shear-locking and possible solutions. In the subsequent part we introduce two mixed variational approaches to the Reissner-Mindlin plate and prove their well-posedness. Section three is devoted to the introduction of the finite element formulations focusing on the Hu-Zhang element, along with its polytopal template and the novel transformation technique. In section four we test the formulations and the element construction for shear-locking, a curved boundary, and an L-shaped domain. Finally, we present our conclusions and outlook.


\subsection{Notation}
The following notation is used throughout this work.
Exception to these rules are made clear in the precise context.
\begin{itemize}
    \item Vectors are defined as bold lower-case letters $\vb{v}, \, \bm{\xi}$
    \item Matrices are bold capital letters $\bm{M}$
    \item Fourth-order tensors are designated by the blackboard-bold format $\mathbb{A}$
    \item We designate the Cartesian basis as $\{\vb{e}_1, \, \vb{e}_2, \, \vb{e}_3\}$
    \item Summation over indices follows the standard rule of repeating indices. Latin indices represent summation over dimension 3, whereas Greek indices define summation over dimension 2  
    \item The angle-brackets are used to define scalar products of arbitrary dimensions $\con{\vb{a}}{\vb{b}} = a_i b_i$, $\con{\bm{A}}{\bm{B}} = A_{ij}B_{ij}$
    \item The matrix product is used to indicate all partial-contractions between a higher-order and a lower-order tensor $\bm{A}\vb{v} = A_{ij} v_j \vb{e}_i$, $\mathbb{A}\bm{B} = A_{ijkl}B_{kl}\vb{e}_i \otimes \vb{e}_j$
    \item Subsequently, we define various differential operators based on the Nabla-operator $\nabla = \partial_i \vb{e}_i$, which is defined with respect to the dimension of the domain
    \item Volumes, surfaces and curves of the physical domain are identified via $\body$, $\surf$ and $\curv$, respectively. Their counterparts on the reference domain are $\Omega$, $\Gamma$ and $\gamma$
    \item Tangent and normal vectors on the physical domain are designated by $\vb{t}$ and $\vb{n}$, respectively. On the reference domain we use $\bm{\tau}$ and $\bm{\nu}$ 
\end{itemize}
We define the constant space of symmetric matrices as
\begin{align}
    \Sym(d) = \{ \bm{M} \in \R^{d \times d} \; | \; \bm{M} = \bm{M}^{T} \} \, .
\end{align}
Further, for our variational formulations we introduce the following Hilbert spaces and their respective norms 
\begin{subequations}
    \begin{align}
    \Hone(\body) &= \{ u \in \Le(\body) \; | \; \nabla u \in [\Le(\body)]^d \} \, , & \norm{u}^2_{\Hone} &= \norm{u}^2_\Le + \norm{\nabla u}^2_\Le \, , \\ 
    \H^2(\surf) &= \{ u \in \Hone(\surf) \; | \; \D \nabla u \in [\Le(\surf)]^{2\times 2} \} \, , & \norm{u}^2_{\H^2} &= \norm{u}^2_{\Hone} + \norm{\D \nabla u}^2_\Le \, , \\ 
    \Hc{,\surf} &= \{ \vb{v} \in [\Le(\surf)]^2 \; | \; \curl \vb{v} \in \Le(\surf) \} \, , & \norm{\vb{v}}^2_{\Hc{}} &= \norm{\vb{v}}^2_\Le + \norm{\curl \vb{v}}^2_\Le \, , \\
    \Hd{,\surf} &= \{ \vb{v} \in [\Le(\surf)]^2 \; | \; \di \vb{v} \in \Le(\surf) \} \, , & \norm{\vb{v}}^2_{\Hd{}} &= \norm{\vb{v}}^2_\Le + \norm{\di \vb{v}}^2_\Le \, ,
\end{align}
\end{subequations}
where $\D \vb{v} = \vb{v} \otimes \nabla$.
The spaces are based on the Lebesgue space
\begin{align}
    \Le(\body) &= \{ u : \body \to \mathbb{R} \; | \; \| u \|_\Le < \infty  \} \, , && \|u\|_\Le^2 = \int_\body u^2 \, \dd \body \, .
\end{align}
Hilbert spaces with vanishing traces are marked with a zero-subscript, for example $\Honez(\surf)$. 
Scalar products pertaining to the Hilbert spaces are indicated by a subscript on the angle-brackets
\begin{align}
    \con{u}{v}_{\Le(\body)} = \int_\body \con{u}{v} \, \dd \body \, .
\end{align}
When the domain is clear from context, we omit it from the subscript.
Finally, we define the space of symmetric square integrable matrices with a square integrable divergence as 
\begin{align}
    \HsD{,\surf} = \{ \bm{M} \in [\Hd{,\surf}]^2 \; | \; \bm{M} = \bm{M}^T  \} \, , 
\end{align}
where $\Di \bm{M} = \bm{M} \nabla = (\bm{M}_{,i}) \vb{e}_i$.

\section{The linear Reissner-Mindlin plate}

\subsection{Energy functional}

The energy functional of linear elasticity is given by the quadratic form
\begin{align}
    &I(\ud) = \int_\body \dfrac{1}{2} \langle \bm{\varepsilon}, \, \Cm \bm{\varepsilon} \rangle \, \dd \body - \langle \ud ,\, \vb{f} \rangle \, \dd \body \, , && \bm{\varepsilon} = \sym\D\ud \, , && \D \ud = \ud \otimes \nabla \, , 
\end{align}
where $\ud : \overline{\body} \subset \R^3 \to \R^3$ is the displacement field, $\Cm \in \R^{3 \times 3 \times 3 \times 3}$ is the tensor of material constants and $\vb{f} : \body \subset \R^3 \to \R^3$ represents the body forces.
For a linear isotropic Saint Venant-Kirchhoff material
\begin{align}
    \Cm = \lambda \one \otimes \one + 2 \mu \J  \, , 
\end{align}
where $\{\lambda, \mu \}$ are the Lam\'e constants, $\one : \R^3 \to \R^3$ is the second order identity tensor, and $\J: \R^{3\times 3} \to \R^{3\times 3}$ is the fourth order identity tensor, one can split the quadratic form of the internal energy between membrane and shear strains
\begin{align}
     \dfrac{1}{2}\int_\body  \langle \sym\D\ud, \, \Cm \sym\D \ud \rangle \, \dd \body = \dfrac{1}{2}\int_\body  \langle \bm{\epsilon}, \, (\lambda \one \otimes \one + 2 \mu \J) \bm{\epsilon}  \rangle + 2\langle \bm{\gamma}, \, \mu \J  \bm{\gamma}  \rangle \, \dd \body \, .
\end{align}
Under the plane stress assumption $\sigma_{33} = 0$, one can neglect out-of-plane component $\varepsilon_{33}$ in the energy functional as it does not produce energy, such that the strain tensors read
\begin{align}
            \bm{\epsilon} = \varepsilon_{\alpha \beta}\vb{e}_\alpha \otimes \vb{e}_\beta \, , && \bm{\gamma} = \varepsilon_{\alpha 3} (\vb{e}_\alpha \otimes \vb{e}_3 + \vb{e}_3 \otimes \vb{e}_\alpha) \, .
        \end{align}
Further, the in-plane relation between strains and stresses reduces to
\begin{align}
    &\bm{\sigma} = \mathbb{D} \bm{\epsilon} \, , && \mathbb{D} = \dfrac{E}{1- \nu^2}[\nu \one \otimes \one + (1-\nu) \mathbb{J}] \, ,
\end{align}
where the tensors are now with respect to the two-dimensional Euclidean space $\R^2$.
The Reissner-Mindlin plate formulation arises under the following kinematical assumption for the displacement field
\begin{align}
    &\ud(x,y,z) = w \vb{e}_3 - z \phi_\alpha \vb{e}_\alpha \, , &&  w = w (x,y) \, , && \phi_\alpha = \phi_\alpha (x,y)  \, ,
\end{align}
such that the deflection $w$ and the small rotations $\phi_\alpha$ are functions of the $x-y$ plane.
Consequently, the strain tensors take the form 
\begin{subequations}
\begin{align}
    \bm{\epsilon} &= -\dfrac{z}{2}(\phi_{\alpha,\beta} + \phi_{\beta,\alpha}) \vb{e}_\alpha \otimes \vb{e}_\beta = - z \sym \D \bm{\phi} \, , && \bm{\phi} = \phi_\alpha \vb{e}_\alpha \, , \\
    \bm{\gamma} &=  \dfrac{1}{2} (w_{,\alpha} - \phi_\alpha) (\vb{e}_\alpha \otimes \vb{e}_3 + \vb{e}_3 \otimes \vb{e}_\alpha)  \, ,
\end{align}
\end{subequations}
where the gradient operator $\D$ is now with respect to the in-plane variables $\{x,y\}$.
The quadratic form of the out-of-plane shear strains can be further simplified to
\begin{align}
    2\langle \bm{\gamma}, \, \mu \J  \bm{\gamma}  \rangle = 2\mu \langle \bm{\gamma}, \,  \bm{\gamma}  \rangle = 2 \mu \dfrac{1}{4} [ 2(w_{,x} - \phi_1)^2 +2(w_{,y} - \phi_2)^2   ] = \mu \| \nabla w - \bm{\phi} \|^2 \, ,
\end{align}
where $\nabla = \vb{e}_\alpha \partial_\alpha$.
The internal energy takes the form 
\begin{align}
    \dfrac{1}{2} \int_\body z^2 \langle \sym \D \bm{\phi}, \, \mathbb{D} \sym \D \bm{\phi} \rangle + \mu \| \nabla w - \bm{\phi} \|^2 \, \dd \body 
    &= \dfrac{1}{2} \int_\surf \int_{-t/2}^{t/2} z^2 \langle \sym \D \bm{\phi}, \, \mathbb{D} \sym \D \bm{\phi} \rangle + \mu \| \nabla w - \bm{\phi} \|^2 \, \dd z \, \dd \surf \notag \\
    &= \dfrac{1}{2} \int_\surf \dfrac{t^3}{12} \langle \sym \D \bm{\phi}, \, \mathbb{D} \sym \D \bm{\phi} \rangle + \mu  t \| \nabla w - \bm{\phi} \|^2 \, \dd \surf \, ,
\end{align}
by splitting the volume integral between the surface and out-of-plane variable $z$.
By defining the volume forces accordingly, $\vb{f} = f \vb{e}_3$, and integrating the external work over the out-of-plane variable one finds the minimisation problem of the linear Reissner-Mindlin plate
\begin{align}
    I(w, \bm{\phi}) &=  \dfrac{1}{2} \int_\surf \dfrac{t^3}{12} \langle \sym \D \bm{\phi}, \, \mathbb{D} \sym \D \bm{\phi} \rangle + \mu  t \| \nabla w - \bm{\phi} \|^2 \, \dd \surf - \int_\surf t w  f  \, \dd \surf  \quad \to \min \quad \text{w.r.t.} \quad \{w, \bm{\phi}\} \, .
\end{align}
\begin{remark}
    The kinematical assumption of the Reissner-Mindlin plate results in constant shear strains and stresses across the cross-section of the plate. However, the latter violates the assumption of vanishing shear stresses on the lower and upper parts of the plate and further, results in a higher stiffness. As such, we add the shear correction factor $k_s$ \cite{gruttmann_shear_2001,gruttmann_shear_2017}, which is usually $k_s = 5/6$, to rectify the shear stiffness of the formulation
    \begin{align}
        \int_\surf \mu  t \| \nabla w - \bm{\phi} \|^2 \, \dd \surf \quad \to \quad \int_\surf \ksmu  t \| \nabla w - \bm{\phi} \|^2 \, \dd \surf \, .
    \end{align}
\end{remark}

\subsection{Variational formulation}
The variational form of the Reissner-Mindlin plate is derived by taking variations of the energy functional with respect to the deflection $w$ and the rotations $\bm{\phi}$. The former yields
\begin{align}
    \delta_w I = \int_\surf \ksmu t \langle \nabla \delta w , \, \nabla w - \bm{\phi} \rangle - t \, \delta w f \, \dd \surf = 0 \, ,    
    \label{eq:var_w}
\end{align}
whereas the latter results in
\begin{align}
    \delta_{\bm{\phi}} I = \int_\surf \dfrac{t^3}{12} \langle \sym \D \delta \bm{\phi} , \, \mathbb{D} \sym \D \bm{\phi} \rangle - \ksmu t \langle \delta \bm{\phi} , \, \nabla w - \bm{\phi} \rangle \, \dd \surf  = 0 \, .
    \label{eq:var_phi}
\end{align}
We apply Green's formula to the variation with respect to $w$ (\cref{eq:var_w}) to find
\begin{align}
    \delta_w I = \int_{\curv_N^w} \ksmu t  \, \delta w \langle \vb{n} , \, \nabla w - \bm{\phi} \rangle \, \dd \curv   -\int_\surf \ksmu t \langle  \delta w , \, \di(\nabla w - \bm{\phi}) \rangle \, \dd \surf  - \int_\surf t \, \delta w f \, \dd \surf = 0 \, ,    
    \label{eq:varp_w}
\end{align}
where we split the boundary of the domain into Dirichlet and Neumann boundaries with respect to the deflection $\partial \surf = \curv_D^w \cup \curv_N^w$ such that $\curv_D^w \cap \curv_N^w = \emptyset$.
Analogously, applying partial integration to \cref{eq:var_phi} yields
\begin{align}
    \delta_{\bm{\phi}} I = \int_{\curv_N^\phi} \dfrac{t^3}{12} \langle \delta \bm{\phi} , \, (\mathbb{D} \sym \D \bm{\phi}) \vb{n} \rangle \, \dd \curv  -\int_\surf \dfrac{t^3}{12} \langle \delta \bm{\phi} , \, \Di (\mathbb{D} \sym \D \bm{\phi}) \rangle \, \dd \surf - \int_\surf \ksmu t \langle \delta \bm{\phi} , \, \nabla w - \bm{\phi} \rangle \, \dd \surf  = 0 \, ,
    \label{eq:varp_phi}
\end{align}
where we again split the boundary of the domain into Dirichlet and Neunmann boundaries such that $\partial \surf = \curv_D^\phi \cup \curv_N^\phi$ and $\curv_D^\phi \cap \curv_N^\phi = \emptyset$.
From \cref{eq:varp_w} and \cref{eq:varp_phi} we extract the boundary value problem 
\begin{subequations}
    \begin{align}
        - k_s  \mu t \di (\nabla w - \bm{\phi}) &= t f && \text{in} && \surf \, ,  \\
        - \dfrac{t^3}{12} \Di (\mathbb{D} \sym \D \bm{\phi}) - \ksmu t  (\nabla w - \bm{\phi}) &= 0 && \text{in} &&  \surf \, ,  \\
        w &= \widetilde{w} && \text{on} && \curv_D^w \, , \\
        \bm{\phi} &= \widetilde{\bm{\phi}} && \text{on} && \curv_D^\phi \, , \\
        \langle \vb{n}, \,  \nabla w - \bm{\phi} \rangle  &= 0 && \text{on} && \curv_N^w \, , \\
        \dfrac{t^3}{12} (\mathbb{D} \sym \D \bm{\phi}) \vb{n} &= \widetilde{\vb{m}} && \text{on} && \curv_N^\phi \, .
    \end{align}
\end{subequations}
The Neumann boundary condition for the gradient of the rotation is given by the bending moments $\widetilde{\vb{m}}$.
By summing the variations with respect to the deflection $w$ and the rotations $\bm{\phi}$ we find the variational problem  
\begin{align}
    \int_\surf \dfrac{t^3}{12} \langle \sym \D \delta \bm{\phi} , \, \mathbb{D} \sym \D \bm{\phi} \rangle + \ksmu t \langle \nabla \delta w - \delta \bm{\phi} , \, \nabla w - \bm{\phi} \rangle \, \dd \surf  =  \int_\surf t \, \delta w \, f \, \dd \surf \, .
    \label{eq:prob_primal}
\end{align}
The problem is well-posed for $\{w,\bm{\phi}\} \in \Honez(\surf) \times [\Honez(\surf)]^2$, but is susceptible to shear locking \cite{bathe_mitc7_1989}.

\subsection{Shear locking} \label{sec:locking}
In order to observe the problem of shear locking we reformulate the variational problem \cref{eq:prob_primal} into
\begin{align}
    \int_\surf \langle \sym \D \delta \bm{\phi} , \, \mathbb{D}_* \sym \D \bm{\phi} \rangle + \dfrac{\ksmu}{t^2} \langle \nabla \delta w - \delta \bm{\phi} , \, \nabla w - \bm{\phi} \rangle \, \dd \surf  =  \int_\surf  \delta w \, g \, \dd \surf \, ,
\end{align}
where we defined $\mathbb{D}_* = (1/12) \mathbb{D}$, divided the entire equation by $t^3$ and set the volume forces $f = t^2 g$.
Now, let the thickness of the plate approach zero $t \to 0$, then the term $\ksmu / t^2 \to +\infty$ may become infinitely large.
In order for the energy to remain finite, the equation must now satisfy the Kirchhoff-Love constraint in the limit $t \to 0$
\begin{align}
    \nabla w - \bm{\phi} = 0 \quad \iff \quad \nabla w = \bm{\phi} \, . 
\end{align}
However, the variables live in incompatible spaces
\begin{align}
    \nabla w \in \nabla \Hone(\surf) = \Hc{,\surf} \cap \ker(\curl) \nsubseteq [\Hone(\surf)]^2 \ni \bm{\phi} \, . 
\end{align}
In other words, the rotations $\bm{\phi}$ are defined on a space with a higher regularity than the gradient of the deflection $\nabla w$. For example, if the deflection is $\C^0(\surf)$-continuous, than its gradient is only tangential-continuous. At the same time, discretizations of the rotation $\bm{\phi}$ are at least $[\C^0(\surf)]^2$-continuous for $[\Hone(\surf)]^2$-conformity, such that at the limit $t \to 0$, the continuous rotations approximate the discontinuous gradient of the deflection. Depending on the fineness of the mesh and its polynomial power, that approximation is either sub-optimal \cite{sky_hybrid_2021} or impossible, such that locking occurs \cite{bathe_mitc7_1989,hale_simple_2018,pechstein_tdnns_2017}.  

The strong form of the reformulated variational problem \cref{eq:prob_primal} reads
\begin{subequations}
    \begin{align}
    - \dfrac{k_s  \mu }{t^2} \di (\nabla w - \bm{\phi}) &= g && \text{in} && \surf \, ,  \\
        -  \Di (\mathbb{D}_* \sym \D \bm{\phi}) -  \dfrac{\ksmu}{t^2}  (\nabla w - \bm{\phi}) &= 0 && \text{in} &&  \surf \, ,  
\end{align}
\end{subequations}
where we ignored boundary conditions. Now one can introduce a mixed formulation that circumvents the large constant $\lim_{t \to 0} \ksmu / t^2 \to + \infty$ that arises for very thin plates by introducing a new unknown for the shear stresses
\begin{align}
     \vb{q} = -\ksmu / t^2 (\nabla w - \bm{\phi}) \, ,
\end{align}
such that the constraint imposed on the variational problem reads
\begin{align}
    \int_\surf \langle \delta \vb{q} , \, \nabla w - \bm{\phi} \rangle + \dfrac{t^2}{\ksmu} \langle \delta \vb{q}, \, \vb{q} \rangle \, \dd \surf = 0 \, .
\end{align}
Clearly, at the limit $t \to 0$, the difference $\nabla w- \bm{\phi}$ is to vanish in the $\Le$-sense.
Alternatively, locking could be avoided by defining the rotations in the compatible space $\Hc{,\surf}$, which is accomplished in the context of the TDNNS method \cite{pechstein_analysis_2018,pechstein_anisotropic_2012,pechstein_tdnns_2017}, where a new unknown for the bending moments is introduced in $\HdD{,\surf}$
\begin{align}
    &\bm{M} = \mathbb{D}_* \sym \D \bm{\phi} \, , && \bm{M} \in \HdD{,\surf} = \{ \bm{M} \in [\Le(\surf)]^{2 \times 2} \; | \;  \bm{M} = \bm{M}^T \, , \; \di \Di \bm{M} \in \H^{-1}(\surf) \} \, .
\end{align}
The complexity of constructing $\HdD{,\surf}$-conforming subspaces is circumvented in \cite{pechstein_tdnns_2017} by instead relying on Hellan-Herrmann-Johnson finite elements and proving stability in the discrete setting.
\begin{remark}
    The definition of the $\HdD{,\surf}$-space in the context of the TDNNS method is such that $\Di \HdD{,\surf} \subseteq \H^{-1}(\di, \surf) = \Hc{,\surf}'$. In other words, the divergence of elements of $\HdD{,\surf}$ are in the dual space of $\Hc{,\surf}$. Alternatively, in the context of the $\di \Di$-complex \cite{CRMECA_2023__351_S1_A8_0,PaulyDiv,Hu2021,Chen2022} used for biharmonic equations, the $\HdD{,\surf}$ space is defined with a higher regularity as $\HdD{,\surf} = \{ \bm{M} \in [\Le(\surf)]^{2 \times 2} \; | \;  \di \Di \bm{M} \in \Le(\surf) \}$.
\end{remark}

\subsection{A mixed problem}
Interestingly, the problem of shear-locking can be partially mitigated by the Hellinger-Reissner principle, where the bending moments are defined as a new variable
\begin{align}
    \bm{M} = \mathbb{D}_* \sym \D \bm{\phi} \, .
\end{align}
This approach is analogous to the TDNNS-method in \cite{pechstein_tdnns_2017}. However, we employ different spaces with the exact sequence property \cite{arnold_finite_2008} 
\begin{subequations}
    \begin{align}
        \airy \H^2(\surf) &= \HsD{,\surf} \cap \ker(\Di) \, , && \airy(\cdot) = [\bm{R}\nabla(\cdot)] \otimes \nabla \bm{R}^T \, ,  && \bm{R} = \begin{bmatrix}
            0 & 1 \\ -1 & 0 
        \end{bmatrix} \, , \\
        \Di \HsD{,\surf} &= [\Le(\surf)]^2 \, ,
    \end{align}
\end{subequations}
which is depicted in \cref{fig:seq}. The latter is a sub-sequence of the elasticity sequence \cite{pauly_elasticity_2022,pauly_hilbert_2022,chen}.
\begin{figure}
		\centering
		\input{figs/seq.tex}
		\caption{An exact sub-sequence of the elasticity sequence. The range of the Airy-operator is exactly the kernel of the divergence operator, and the divergence operator applied to the $\HsD{,\surf}$-space is a surjection onto $[\Le(\surf)]^2$.}
		\label{fig:seq}
	\end{figure}
We adapt the strong form to 
\begin{subequations}
    \begin{align}
    - \dfrac{k_s  \mu }{t^2} \di (\nabla w - \bm{\phi}) &= g && \text{in} && \surf \, ,  \\
        - \Di \bm{M} -  \dfrac{\ksmu}{t^2}  (\nabla w - \bm{\phi}) &= 0 && \text{in} &&  \surf \, , \\
        \sym \D \bm{\phi} - \mathbb{A} \bm{M}  &= 0 && \text{in} &&  \surf \, ,
\end{align}
\end{subequations}
where we employ the compliance tensor 
\begin{align}
    \mathbb{A} = \mathbb{D}_*^{-1} = 12 \, \mathbb{D}^{-1} = \dfrac{12}{E} [(1+ \nu) \mathbb{J} - \nu \one \otimes \one ] \, .
\end{align}
Applying test functions to the equilibrium equations and the constraint equation yields
\begin{subequations}
    \begin{align}
        \int_\surf - \dfrac{k_s  \mu }{t^2} \langle \delta w ,\,  \di (\nabla w - \bm{\phi}) \rangle \, \dd \surf &= \int_\surf \delta w \,  g \, \dd \surf \, , \label{eq:mixed1} \\
        \int_\surf  -  \langle \delta \bm{\phi} , \,  \Di \bm{M} \rangle -  \dfrac{\ksmu}{t^2} \langle \delta \bm{\phi} , \,  \nabla w - \bm{\phi} \rangle \, \dd \surf &= 0 \, , \\
        \int_\surf   \langle \delta \bm{M} ,\, \sym \D \bm{\phi} \rangle - \langle \delta \bm{M} , \, \mathbb{A} \bm{M} \rangle \, \dd \surf  &= 0 \, . \label{eq:mixed3}
    \end{align}
\end{subequations}
Partial integration of \cref{eq:mixed1} and \cref{eq:mixed3} results in
\begin{subequations}
    \begin{align}
        \int_\surf  \dfrac{k_s  \mu }{t^2} \langle \nabla \delta w ,\,  \nabla w - \bm{\phi} \rangle \, \dd \surf - \int_{\curv_N^w}  \dfrac{k_s  \mu }{t^2} \delta w \langle \vb{n} ,\,  \nabla w - \bm{\phi} \rangle \, \dd \curv &= \int_\surf \delta w \,  g \, \dd \surf \, ,  \\
        \int_{\curv_N^M} \con{\delta \bm{M} \vb{n} }{ \bm{\phi} } \, \dd \curv - \int_\surf  \langle \Di \delta \bm{M} ,\, \bm{\phi} \rangle + \langle \delta \bm{M} , \, \mathbb{A} \bm{M} \rangle \, \dd \surf  &= 0 \, . 
    \end{align}
\end{subequations}
As such, we define the following boundary conditions 
\begin{subequations}
    \begin{align}
        w &= \widetilde{w} && \text{on} && \curv_D^w \, , \\
        \bm{M}\vb{n} &= \widetilde{\bm{M}} \vb{n} && \text{on} && \curv_D^M \, , \\
        \bm{\phi} &= \widetilde{\bm{\phi}} && \text{on} && \curv_N^M \, .
    \end{align}
\end{subequations}
\begin{remark}[Kinematical boundary conditions]
    The kinematical boundary conditions of the deflection are obvious.
    The rotations $\bm{\phi}$ however, are controlled implicitly through the boundary of the bending moments $\bm{M}$. Namely, if the normal projection of the bending moments is set to zero on the Dirichlet boundary $\widetilde{\bm{M}}\vb{n}|_{\curv_{D}^M} = 0$, this implies a kinematical hinge for the rotations $\bm{\phi}$ on $\curv_{D}^M$. Conversely, on the Neumann boundary of the bending moments $\curv_{N}^M$, the rotations are prescribed weakly via $\int_{\curv_N^M} \con{\delta \bm{M} \vb{n} }{ \bm{\phi} } \, \dd \curv$, such that omitting the term implies $\bm{\phi}|_{\curv_N^M} = 0$.  
\end{remark}
The complete variational problem can now be written as
\begin{subequations}
    \label{eq:prob_mixed}
    \begin{align}
    \int_\surf \con{\delta \bm{M}}{\mathbb{A} \bm{M}} +  \con{\Di \delta \bm{M}}{\bm{\phi}} \, \dd \surf &= \int_{\curv_N^M} \con{\delta \bm{M} \, \vb{n}}{\widetilde{\bm{\phi}}} \, \dd \curv && \forall \, \delta \bm{M} \in \HsD{,\surf} \, , \\ 
    \int_\surf \con{\delta \bm{\phi}}{\Di \bm{M}} - \dfrac{\ksmu}{t^2} \con{\nabla \delta w - \delta \bm{\phi}}{\nabla w - \bm{\phi}} \, \dd \surf &= - \int_\surf \delta w \, g \, \dd \surf && \forall \, \{\delta w, \delta \bm{\phi}\} \in \Honez(\surf) \times [\Le(\surf)]^2 \, .
\end{align}
\end{subequations}
With the definition of $\bm{\phi} \in [\Le(\surf)]^2$ there holds $\nabla w \in \Hc{,\surf} \subset [\Le(\surf)]^2$ such that for the limit $t \to 0$ there exists $\bm{\phi} = \nabla w$ in the $\Le$-sense and shear-locking can be avoided by $\nabla w - \bm{\phi} = 0$ as per the Kirchhoff-Love condition.
\begin{theorem}[Well-posedness for $t > 0$]
    The variational problem in \cref{eq:prob_mixed} is well-posed in the space $\X(\surf) = \Y(\surf) \times \HsD{,\surf}$ with $\Y(\surf) = \Honez(\surf) \times [\Le(\surf)]^2$ such that $\{\delta w, \delta \bm{\phi}, \delta \bm{M}\} \in \X(\surf)$ and $\{w, \bm{\phi}, \bm{M}\} \in \X(\surf)$, assuming a contractible domain. The spaces are equipped with the norms
    \begin{subequations}
        \begin{align}
            \| \{w, \bm{\phi}\} \|_\Y^2 &= \|w \|_{\Hone}^2 + \|\bm{\phi} \|_{\Le}^2 \, , \\
            \| \{w, \bm{\phi}, \bm{M} \} \|_\X^2 &= \| \{w, \bm{\phi}\} \|_\Y^2 + \|\bm{M} \|_{\HsD{}}^2 \, ,
        \end{align}
    \end{subequations}
    and there holds the stability estimate
    \begin{align}
        \norm{\{w, \bm{\phi}\}}_\Y + \norm{ \{ \bm{M}, \nabla w - \bm{\phi} \} }_{\HsD{} \times \Le } \leq c \norm{g}_{\Y'} \, ,
    \end{align}
    where $c = c(\mathbb{A},\ksmu, t)$, and $\Y'$ is the dual space of $\Y$.
\end{theorem}

\begin{proof}
    For the proof we introduce the shear stresses $\vb{q} = (-\ksmu / t^2) (\nabla w - \bm{\phi}) \in [\Le(\surf)]^2$ as additional unknown and $\Z(\surf) = \HsD{,\surf} \times [\Le(\surf)]^2$. Then the variational problem reads
    \begin{subequations}
        \begin{align}
            \int_\surf \con{\delta \bm{M}}{\A \bm{M}} + \dfrac{t^2}{\ksmu} \con{\delta \vb{q}}{ \vb{q}} + \con{\Di \delta \bm{M}}{\bm{\phi}} + \con{\delta \vb{q}}{\nabla w - \bm{\phi}} \, \dd \surf &=  0 && \forall \, \{\delta \bm{M}, \delta \vb{q}\} \in \Z (\surf) \, ,\\
            \int_\surf \con{\delta \bm{\phi}}{ \Di  \bm{M}} + \con{\nabla \delta w - \delta\bm{\phi}  }{ \vb{q}} \, \dd \surf &=  - \int_\surf \delta w \, g \, \dd \surf  && \forall \, \{\delta w, \, \delta \bm{\phi}\} \in \Y (\surf) \, , 
\end{align}
\end{subequations}
    where we neglected Neumann boundary terms for simplicity. 
    The space $\Z$ is equipped with the norm $\norm{\{ \bm{M} , \vb{q} \}}_\Z^2 = \norm{\bm{M}}_{\HsD{}}^2 + \norm{\vb{q}}_\Le^2$.
    The variational problem gives rise to the following bilinear and linear forms 
    \begin{subequations}
        \begin{align}
        a(\{\delta \bm{M}, \delta \vb{q} \}, \{ \bm{M} , \vb{q} \}) &= \int_\surf \con{\delta \bm{M}}{\A \bm{M}} + \dfrac{t^2}{\ksmu} \con{\delta \vb{q}}{ \vb{q}} \, \dd \surf \, , \\ 
        b(\{ \delta w, \delta \bm{\phi} \}, \{ \bm{M} , \vb{q} \}) &= \int_\surf \con{\delta \bm{\phi}}{ \Di  \bm{M}} + \con{\nabla \delta w - \delta\bm{\phi}  }{ \vb{q}} \, \dd \surf \, , \\
        l( \delta w) &= - \int_\surf \delta w \, g \, \dd \surf \, , 
    \end{align}
    \end{subequations}
    such that it can be written as the saddle point problem
    \begin{subequations}
        \begin{align}
            a(\{\delta \bm{M}, \delta \vb{q} \}, \{ \bm{M} , \vb{q} \}) + b(\{ \delta\bm{M} , \delta\vb{q} \},\{ w,  \bm{\phi} \}) &= 0 && \forall \, \{\delta \bm{M}, \delta \vb{q}\} \in \Z (\surf) \, , \\ 
            b(\{ \delta w, \delta \bm{\phi} \}, \{ \bm{M} , \vb{q} \}) &= l( \delta w)  && \forall \, \{\delta w, \, \delta \bm{\phi}\} \in \Y (\surf) \, .
        \end{align}
        \label{eq:problem}
    \end{subequations}
    Existence and uniqueness follows by the Brezzi theorem \cite{Bra2013}. The continuity of the linear form $l(\cdot)$ is obvious. The continuity of $a(\cdot,\cdot)$ can be shown via
    \begin{align}
        a(\{\delta \bm{M}, \delta \vb{q} \}, \{ \bm{M} , \vb{q} \}) &= \con{\delta \bm{M}}{\A \bm{M}}_\Le + \dfrac{t^2}{\ksmu} \con{\delta \vb{q}}{ \vb{q}}_\Le  \notag \\
        &\overset{CS}{\leq}  c_A\| \delta \bm{M} \|_\Le \norm{\bm{M}}_\Le + \dfrac{t^2}{\ksmu} \norm{\delta \vb{q}}_\Le \norm{\vb{q}}_\Le  \notag \\
        & \overset{\phantom{a}}{\leq} 2 c_1 \norm{\{\delta \bm{M}, \delta \vb{q}\}}_\Z \norm{\{ \bm{M},  \vb{q}\}}_\Z \, ,
        \label{eq:aconti}
    \end{align}
    where we used the positive definiteness of $\A$ on symmetric matrices
    \begin{align}
        \exists \, \{k_A, c_A\} > 0: \quad  k_A \norm{\bm{S}}^2 \leq \con{\bm{S}}{\A \bm{S}} \leq c_A \norm{\bm{S}}^2 \quad \forall \, \bm{S} \in \Sym(2) \, ,
    \end{align}
    and the Cauchy-Schwarz inequality. As such, the continuity constant reads $\alpha_1 = 2 c_1$ with $c_1 = \max\{ c_A, t^2 / (\ksmu)\}$. The continuity of $b(\cdot, \cdot)$ is given by
    \begin{align}
        b(\{ \delta w, \delta \bm{\phi} \}, \{ \bm{M} , \vb{q} \}) &= \con{\delta \bm{\phi}}{ \Di  \bm{M}}_\Le + \con{\nabla \delta w - \delta\bm{\phi}  }{ \vb{q}}_\Le \notag \\
        &\overset{CS}{\leq} \norm{\delta \bm{\phi}}_\Le\norm{ \Di  \bm{M}}_\Le + \norm{\nabla \delta w - \delta\bm{\phi}  }_\Le \norm{ \vb{q}}_\Le \notag \\
        &\overset{T}{\leq} \norm{\delta \bm{\phi}}_\Le\norm{ \Di  \bm{M}}_\Le + (\norm{\nabla \delta w}_\Le + \norm{\delta\bm{\phi}  }_\Le) \norm{ \vb{q}}_\Le \notag \\
        &\overset{\phantom{a}}{\leq} 3\norm{ \{ \delta w, \delta \bm{\phi} \}}_\Y \norm{ \{ \bm{M}, \vb{q} \}}_\Z \, ,
    \end{align}
    using again the Cauchy-Schwarz and triangle inequalities such that the continuity constant reads $\alpha_2 = 3$.
    Next we need to show coercivity of $a(\cdot, \cdot)$ on the kernel of $b(\{\delta w, \delta \bm{\phi}\}, \cdot)$. 
    The kernel is characterised by
        \begin{align}
        \ker(B) &= \{ \{\bm{M}, \vb{q}\} \in \Z(\surf) \; | \; b(\{ \delta w, \delta \bm{\phi} \}, \{ \bm{M} , \vb{q} \}) = 0 \quad \forall \, \{ \delta w, \delta \bm{\phi} \} \in \Y(\surf) \} \, , 
        \label{eq:kernel}
    \end{align}
    and implies
        \begin{align}
            \con{\delta \bm{\phi}}{ \Di  \bm{M} - \vb{q}}_\Le &= 0 && \forall \, \delta \bm{\phi} \in [\Le(\surf)]^2 \, , 
        \end{align}
    holds on the kernel. The result follows by setting $\delta w$ to zero.
    The kernel-coercivity of $a(\cdot ,\cdot)$ can now be shown via
    \begin{align}
        a(\{ \bm{M}, \vb{q} \}, \{ \bm{M} , \vb{q} \}) &= \con{ \bm{M}}{\A \bm{M}}_\Le + \dfrac{t^2}{\ksmu} \con{\vb{q}}{ \vb{q}}_\Le \notag \\
        &\overset{\phantom{a}}{\geq} k_A \norm{\bm{M}}_\Le^2 + \dfrac{t^2}{\ksmu} \norm{ \vb{q}}_\Le^2 \notag \\ 
        &\overset{\phantom{a}}{=} k_A \norm{\bm{M}}_\Le^2 + \dfrac{t^2}{2\ksmu} (\norm{ \Di \bm{M} }_\Le^2 + \norm{ \vb{q} }_\Le^2) \notag \\
        &\overset{\phantom{a}}{\geq} \beta_1 \norm{ \{ \bm{M}, \vb{q} \} }_\Z^2  \, ,  
    \end{align}
    since on the kernel of $B$ there holds $\Di \bm{M} = \vb{q}$ in the $\Le$-sense. The coercivity constant reads $\beta_1 = \min\{ k_A , t^2 / (2\ksmu) \}$. Lastly, $b(\cdot, \cdot)$ must satisfy the Ladyzhenskaya–Babuška–Brezzi (LBB) condition: $\exists \beta_2>0$ such that for all $\{ w,\bm{\phi}\}\in Y(\surf)$
    \[
    \sup_{\{\bm{M}, \vb{q}\} \in \Z} \dfrac{b(\{  w,  \bm{\phi} \}, \{ \bm{M} , \vb{q} \})}{\norm{\{\bm{M},\vb{q}\}}_Z}\geq \beta_2\,\norm{\{ w,\bm{\phi}\}}_Y.
    \]
    Let $\{ w,\bm{\phi}\}\in \Y(\surf)$ be given. We define $\vb{q}:= \nabla w-\bm{\phi}\in[\Le(\surf)]^2$ and by the exact sequence property \cref{fig:seq} $\bm{M}\in\HsD{,\surf}$ such that $\Di \bm{M}=\bm{\phi}$ and \cite{arnold_finite_2008,Lederer}
    \begin{align}
        &\norm{\bm{M}}_{\HsD{}} \leq c_H \norm{ \bm{\phi}}_\Le.
        \label{eq:HR}
    \end{align}
    Then, there holds
    \begin{align}
        \sup_{\{\bm{M}, \vb{q}\} \in \Z} \dfrac{b(\{  w,  \bm{\phi} \}, \{ \bm{M} , \vb{q} \})}{\norm{\{\bm{M},\vb{q}\}}_Z}&\geq \dfrac{\norm{\bm{\phi}}_{\Le}^2+\norm{\nabla w-\bm{\phi}}_{\Le}^2}{\sqrt{\norm{\bm{M}}_{\Le}^2+\norm{\Di\bm{M}}_{\Le}^2+\norm{\vb{q}}_{\Le}^2}}\notag\\
        &=\dfrac{2\norm{\bm{\phi}}_{\Le}^2+\norm{\nabla w}_{\Le}^2-2\con{\nabla w}{\bm{\phi}}_{\Le}}{\sqrt{\norm{\bm{M}}_{\Le}^2+\norm{\bm{\phi}}_{\Le}^2+\norm{\nabla w -\bm{\phi}}_{\Le}^2}}\notag\\
        &\overset{Y,\,T}{\geq}\dfrac{(2-\varepsilon^{-1})\norm{\bm{\phi}}_{\Le}^2+(1-\varepsilon)\norm{\nabla w}_{\Le}^2}{\sqrt{c_H^{2}\norm{\bm{\phi}}_{\Le}^2+3\norm{\bm{\phi}}_{\Le}^2+2\norm{\nabla w}_{\Le}^2}}\notag\\
        &\overset{\varepsilon=3/4}{\geq}\dfrac{\dfrac{2}{3}\norm{\bm{\phi}}_{\Le}^2+\dfrac{1}{4}\norm{\nabla w}_{\Le}^2}{\sqrt{c_H^{2}+3}\,\norm{\{\bm{\phi}, w\}}_Y}\notag\\
        &\geq\dfrac{\min\left\{\dfrac{2}{3},\dfrac{1}{1+c_F^2}\right\}\norm{\{\bm{\phi}, w\}}_Y^2}{\sqrt{c_H^{2}+3}\,\norm{\{\bm{\phi}, w\}}_Y}=\beta_2\norm{\{\bm{\phi}, w\}}_Y,
    \end{align}
    where we applied the Young\footnote{Young's inequality: $\exists \, \varepsilon > 0: \quad -\langle x, y \rangle \geq -\dfrac{\norm{x}^2}{2\varepsilon} - \dfrac{\varepsilon \norm{y}^2 }{2}$} and Poincar\'e-Friedrich\footnote{Poincar\'e-Friedrich's inequality: $\exists \, c_F > 0: \quad  \|x\|_{\Le} \leq c_F \|\nabla x \|_\Le \quad \forall x \in \Honez$} inequalities.
    Consequently, the constant $\beta_2$ reads 
    \begin{align}
         &\beta_2 = \dfrac{\min\left\{\dfrac{2}{3},\dfrac{1}{1+c_F^2}\right\}}{\sqrt{c_H^{2}+3}} \, .
    \end{align}
    Thus, there holds by the Brezzi theorem the stability estimate 
    \begin{align}
        c \norm{g}_{\Y'} &\geq \norm{\{w, \bm{\phi}\}}_\Y + \norm{ \{ \bm{M}, \vb{q}\} }_{\Z}    \notag \\
        &\overset{\phantom{a}}{=} \norm{\{w, \bm{\phi}\}}_\Y + \sqrt{ \norm{\bm{M}}_{\HsD{}}^2 + \norm{\vb{q}}_{\Le}^2 } \notag \\
        &\geq \norm{\{w, \bm{\phi}\}}_\Y + \sqrt{ \norm{\bm{M}}_{\HsD{}}^2 +  \left (\dfrac{\ksmu}{t^2} \right )^2 \norm{\nabla w - \bm{\phi}}_{\Le}^2 } \notag \\
        &\geq  \min \left \{1, \dfrac{\ksmu}{t^2} \right \} (\norm{\{w, \bm{\phi}\}}_\Y + \norm{ \{ \bm{M}, \vb{q}\} }_{\Z}) \, ,
    \end{align}
    where we inserted the definition of $\vb{q}$. This concludes the proof.
\end{proof}
Unfortunately, while the formulation is well-posed, its stability constant is not independent of the thickness $t$ such that for $t \to 0$ instability may occur, even though locking in the sense of shear-locking is circumvented by the compatibility of the spaces. In order to derive a robust formulation also for $t \to 0$ we reformulate \cref{eq:problem} via integration by parts of the shear stress, yielding the following variational problem
\begin{subequations}
    \begin{align}
        \int_\surf \con{\delta \bm{M}}{\A \bm{M}} + \dfrac{t^2}{\ksmu} \con{\delta \vb{q}}{ \vb{q}} + \con{\Di \delta \bm{M}}{\bm{\phi}} - (\di \delta \vb{q}) \, w  -  \con{\delta \vb{q}}{\bm{\phi}} \, \dd \surf &= 0 && \forall \, \{\delta \bm{M}, \delta \vb{q} \} \in \Z(\surf) \, , \\
        \int_\surf \con{\delta \bm{\phi}}{\Di \bm{M}} -  \delta  w \, (\di  \vb{q}) - \con{\delta \bm{\phi}}{\vb{q}} \, \dd \surf &= - \int_\surf \delta w \, g \, \dd \surf && \forall \, \{\delta w, \delta \bm{\phi} \} \in \Y(\surf) \, ,
    \end{align}
    \label{eq:fourfield}
\end{subequations}
with $\Z(\surf) = \HsD{,\surf} \times \Hd{,\surf}$ and $\Y(\surf) = \Le(\surf) \times [\Le(\surf)]^2$.
In addition, there arises the Neumann term for the right-hand side $\int_{\curv_N^q} \con{\delta \vb{q}}{\vb{n}} \, \widetilde{w} \, \dd \curv$, which controls the prescribed deflection $\widetilde{w}$ on the boundary of the domain, since $w \in \Le(\surf)$ does not incorporate Dirichlet boundary conditions.
\begin{theorem}[Robustness in $t$]
    Under the assertion of a contractible domain, the variational problem in \cref{eq:fourfield} is well-posed in the space 
    \begin{align}
        \{w, \bm{\phi}, \bm{M}, \vb{q}\}  \in \X(\surf) &= \Y(\surf) \times \Z(\surf) = \Le(\surf) \times [\Le(\surf)]^2 \times  \HsD{,\surf} \times \Hd{,\surf}   \, , 
    \end{align}
    and there holds the stability estimate
    \begin{align}
        \norm{\{w, \bm{\phi}\}}_\Y + \norm{ \{ \bm{M}, \vb{q} \} }_{\Z_t} \leq c \norm{g}_{\Y'} \, ,
    \end{align}
    where the constant $c$ is independent of the thickness $t$ such that $c = c(\mathbb{A},\ksmu)$, and $\norm{\cdot}_{\Z_t}$ is the $t$-dependent norm 
    \begin{align}
        \norm{ \{ \bm{M}, \vb{q} \} }^2_{\Z_t} = \norm{\bm{M}}_{\Le}^2+\|\Di \bm{M}-\vb{q}\|^2_{\Le} + t^2 \norm{\vb{q}}_{\Le}^2 + \norm{\di \vb{q}}^2_{\Le} \, .
    \end{align}
    The space $\Y$ is equipped with the natural norm
    \begin{align}
        \norm{\{w, \bm{\phi}\}}_\Y = \norm{w}_{\Le}^2 + \norm{\bm{\phi}}_{\Le}^2 \, .
    \end{align}
\end{theorem}
\begin{proof}
    The proof follows the ideas of \cite{daVeiga}. We define the bi- and linear forms 
    \begin{subequations}
        \begin{align}
            a(\{\delta \bm{M}, \delta \vb{q} \}, \{ \bm{M} , \vb{q} \}) &= \int_\surf \con{\delta \bm{M}}{\A \bm{M}} + \dfrac{t^2}{\ksmu} \con{\delta \vb{q}}{ \vb{q}} \, \dd \surf \, , \\ 
            b(\{ \delta w, \delta \bm{\phi} \}, \{ \bm{M} , \vb{q} \}) &= \int_\surf \con{\delta \bm{\phi}}{\Di \bm{M}} -  \delta  w \, (\di  \vb{q}) - \con{\delta \bm{\phi}}{\vb{q}} \, \dd \surf \, , \\
        l( \delta w) &= - \int_\surf \delta w \, g \, \dd \surf \, ,
        \end{align}
    \end{subequations}
    such that the problem has the same form as \cref{eq:problem}, but on different spaces.
    The continuity of $a(\cdot,\cdot)$ is given analogously to \cref{eq:aconti}, where the continuity constant is modified to $\alpha_1 = 2\max\{ c_A, 1 / (\ksmu)\}$.
    The continuity of $b(\cdot, \cdot)$ follows via
    \begin{align}
        b(\{ \delta w, \delta \bm{\phi} \}, \{ \bm{M} , \vb{q} \}) &=  \con{\delta \bm{\phi}}{\Di \bm{M}-\vb{q}}_{\Le} -  \con{\delta  w} {\di  \vb{q}}_{\Le}\notag \\ 
        &\overset{CS}{\leq} \norm{\delta \bm{\phi}}_\Le\norm{ \Di  \bm{M}-\vb{q}}_\Le + \norm{ \delta w  }_\Le \norm{ \di \vb{q}}_\Le \notag \\
        &\overset{\phantom{a}}{\leq} 2\norm{ \{ \delta w, \delta \bm{\phi} \}}_\Y \norm{ \{ \bm{M}, \vb{q} \}}_{\Z_t}  \, , 
    \end{align}
    such that $\alpha_2 = 2$.
    The kernel of $b(\{\delta w, \delta \bm{\phi}\}, \cdot)$ is of similar form to \cref{eq:kernel} and implies
    \begin{subequations}
        \begin{align}
            \con{\delta \bm{\phi}}{ \Di  \bm{M} - \vb{q}}_\Le &= 0 && \forall \, \delta \bm{\phi} \in [\Le(\surf)]^2 \, ,  \\
            -\con{\delta w}{ \di \vb{q}}_\Le &= 0 && \forall \, \delta w \in \Le(\surf) \, ,
        \end{align} 
    \end{subequations} 
    such that on the kernel there holds 
    \begin{align}\norm{\{\bm{M},\vb{q}\}}_{\Z_t}^2=\norm{\bm{M}}_{\Le}^2+t^2\norm{\vb{q}}_{\Le}^2\qquad \forall \{\bm{M},\vb{q}\} \in \ker(B).
    \end{align}
    Consequently, $a(\cdot, \cdot)$ is uniformly kernel-coercive in $t$ 
    \begin{align}
        a(\{ \bm{M}, \vb{q} \}, \{ \bm{M} , \vb{q} \}) &= \con{ \bm{M}}{\A \bm{M}}_\Le + \dfrac{t^2}{\ksmu} \con{\vb{q}}{ \vb{q}}_\Le \overset{\phantom{a}}{\geq} k_A \norm{\bm{M}}_\Le^2 + \dfrac{t^2}{\ksmu} \norm{ \vb{q}}_\Le^2
        \overset{\phantom{a}}{\geq} \beta_1 \norm{ \{ \bm{M}, \vb{q} \} }_{\Z_t}^2  \, ,  
    \end{align}
    such that the coercivity constant reads $\beta_1 = \min\{ k_A, 1/(\ksmu) \}$.

     For the LBB-condition let $\{ w,\bm{\phi}\}\in \Y(\surf)$ be given. We define $\vb{q}:=\nabla z\in \Hd{,\surf}$, where $z\in \Honez(\surf)$ solves the  Poisson equation $-\Delta z=w$ such that
     \begin{align}
        \di \vb{q} = -w,\qquad \norm{\vb{q}}_{\Le}+\norm{\di \vb{q}}_{\Le}\le c_1\norm{w}_{\Le}.
     \end{align}
    Further, by the exact sequence property we choose $\bm{M}\in \HsD{,\surf}$ such that $\Di \bm{M}=c^2_1\,\bm{\phi}$ and $\norm{\bm{M}}_{\Le}+\norm{\Di \bm{M}}_{\Le}\le c_H\|\bm{\phi}\|_{\Le}$.
    As such, we obtain 
    \begin{align}
        \sup_{\{\bm{M}, \vb{q}\} \in \Z} \dfrac{b(\{  w,  \bm{\phi} \}, \{ \bm{M} , \vb{q} \})}{\norm{\{\bm{M},\vb{q}\}}_{\Z_t}} &\geq \dfrac{\con{ \bm{\phi}}{\Di \bm{M}}_{\Le} -  \con{  w} {\di  \vb{q}}_{\Le} - \con{ \bm{\phi}}{\vb{q}}_{\Le}}{\sqrt{\norm{\bm{M}}_{\Le}^2 +\norm{\Di \bm{M}- \vb{q}}^2_{\Le}+ t^2 \norm{\vb{q}}_{\Le}^2 + \norm{\di \vb{q}}^2_{\Le} }} \notag \\
        &\overset{Y}{\geq} \dfrac{ c_1^2\norm{\bm{\phi}}_{\Le}^2 + \norm{ w}_{\Le}^2 - \dfrac{\varepsilon\norm{\bm{\phi}}_{\Le}^2}{2} - \dfrac{\norm{\vb{q}}_{\Le}^2}{2\varepsilon} }{\sqrt{\norm{\bm{M}}_{\Le}^2 +\norm{\Di \bm{M}- \vb{q}}^2_{\Le} + t^2 \norm{\vb{q}}_{\Le}^2 + \norm{\di \vb{q}}^2_{\Le} }} \notag \\
        &\overset{\phantom{PE}}{\geq} \dfrac{ \left (c_1^2-\dfrac{\varepsilon}{2}\right)\norm{\bm{\phi}}_{\Le}^2 + \left (1- \dfrac{ c_1^2 }{2 \varepsilon} \right)\norm{ w}_{\Le}^2 }{ \sqrt{2c_H^2\norm{\bm{\phi}}_{\Le}^2 + 2c_1^2(2 + t^2) \norm{ w}_{\Le}^2 }} \notag \\
        &\overset{\varepsilon=c_1^2}{\geq} \dfrac{ \dfrac{c_1^2}{2}\norm{\bm{\phi}}_{\Le}^2 +  \dfrac{ 1 }{2 }\norm{ w}_{\Le}^2 }{ \sqrt{\max\{2c_H^2,6c_1^2\}} \sqrt{ \norm{ \bm{\phi}}_{\Le}^2 + \norm{ w}_{\Le}^2 }}\overset{\phantom{a}}{\geq} \beta_2 \norm{\{ w,  \bm{\phi}\}}_{\Y} \, ,
    \end{align}
    where we applied Young's inequality and the stability estimates of the constructed $\vb{q}$ and $\bm{M}$.
    The constant reads $\beta_2 = \dfrac{\min\{c_1^2,1\}}{2\sqrt{\max\{2c_H^2,6c_1^2\}}}$. 
    Finally, the stability estimate follows by Brezzi's theorem \cite{Bra2013}.
\end{proof}
The proof follows analogously for the case $t > 1$ using natural, $t$-independent norms. Further, we note that the resulting variational problem is in fact comparable to the formulation introduced in \cite{daVeiga}. However, due to our reliance on the symmetric $\HsD{}$-space, we do not need to impose the symmetry of the bending moments $\bm{M}$ weakly, and thus, we do not introduce an auxiliary variable for symmetry. 

\section{Finite element formulation}
The first conforming subspace for $\HsD{,\surf}$ was given by the Arnold-Winther element \cite{arnold_finite_2008}. The element relied on an enriched polynomial space on symmetric tensors for the purpose of stability such that its dimension in the lowest order reads $\dim \AW^2(\surf) = 24$. Further, the element allows for an elasticity sub-complex with a commuting property for sufficiently smooth functions
\begin{subequations}
    \begin{align}
        \airy \Pi_a w &= \Pi_s \airy w \, , & w &\in \H^2(\surf) \, , \\
        \Di \Pi_s \bm{M} &= \Pi_o \Di \bm{M} \, , &  \bm{M} &\in \HsD{,\surf} \, ,
    \end{align}
\end{subequations}
which is depicted in \cref{fig:com_aw}. 
Alternatively, the Hu-Zhang element defines an $\HsD{,\surf}$-conforming subspace using the full polynomial space for symmetric tensors in the polynomial order $p \geq 3$ \cite{hu_family_2014} such that its dimension reads $\dim \HZ^p(\surf) = \dim [\Po^p(\surf) \otimes \Sym(2)] = [\dim\Po^p(\surf)][\dim \Sym(2)] = 3[\dim\Po^p(\surf)] = 3(p+2)(p+1)/2$. 
As such, it can be used to improve the elasticity sub-complex, as shown in \cref{fig:com_hz}.
\begin{figure}
		\centering
        \begin{subfigure}{0.48\linewidth}
    		\centering
    		\input{figs/com_aw}
            \caption{}
            \label{fig:com_aw}
    	\end{subfigure}
        \begin{subfigure}{0.48\linewidth}
    		\centering
    		\input{figs/com_hz}
            \caption{}
            \label{fig:com_hz}
    	\end{subfigure}
		\caption{An exact sub-complex of the elasticity sub-sequence, such that the interpolation operators commute. Here, $\mathcal{A}^5(\surf)$ represents the quintic Argyris space of $\C^1(\surf)$-continuous polynomials, and $\mathit{D}^p(\surf)$ are element-wise discontinuous linear polynomials. In (a) we present the complex using the quadratic Arnold-Winther element, whereas (b) depicts the alternative exact complex using the cubic Hu-Zhang space, such that the discontinuous polynomial space is improved to quadratic power.}
  \label{fig:com_elast}
\end{figure}
We mention that a lowest order construction with $\dim [\Po(\surf)^2 \otimes \Sym(2)] + 3 = 21$ degrees of freedom is also mentioned in \cite{arnold_finite_2008} and \cite{hu_finite_2016}. However, that element formulation is stable for the sequence-pairing with the space of rigid body motions, which is a subspace of $\Po^1(\surf)$.

The finite element formulation of the three-field problem reads
    \begin{align}
    \int_\surf \con{\delta \bm{M}}{\mathbb{A} \bm{M}} +  \con{\Di \delta \bm{M}}{\bm{\phi}} \, \dd \surf &= 0 && \forall \, \delta \bm{M} \in \HZ^p(\surf) \, ,  \label{eq:TFSRM} \\ 
    \int_\surf \con{\delta \bm{\phi}}{\Di \bm{M}} - \dfrac{\ksmu}{t^2} \con{\nabla \delta w - \delta \bm{\phi}}{\nabla w - \bm{\phi}} \, \dd \surf &= - \int_\surf \delta w \, g \, \dd \surf && \forall \, \{\delta w, \delta \bm{\phi}\} \in \U^p(\surf) \times [\mathit{D}^{p-1}(\surf)]^2 \, , \notag 
\end{align}
where $\HZ^p(\surf)$ is the Hu-Zhang element of order $p$, $\U^p(\surf)$ represent the polynomial $\C^0(\surf)$-continuous space of order $p$, and $\mathit{D}^{p-1}(\surf)$ is the space of discontinuous, piece-wise polynomials of order $p-1$.

The discrete quad-field finite element formulation reads
\begin{align}
        &\int_\surf \con{\delta \bm{M}}{\A \bm{M}} + \dfrac{t^2}{\ksmu} \con{\delta \vb{q}}{ \vb{q}} + \con{\Di \delta \bm{M}}{\bm{\phi}} - (\di \delta \vb{q}) \, w  -  \con{\delta \vb{q}}{\bm{\phi}} \, \dd \surf = 0 \quad \forall \, \{\delta \bm{M}, \delta \vb{q} \} \in \HZ^p(\surf) \times \RT^{p-1}(\surf) \, , \notag \\
        &\int_\surf \con{\delta \bm{\phi}}{\Di \bm{M}} -  \delta  w \, (\di  \vb{q}) - \con{\delta \bm{\phi}}{\vb{q}} \, \dd \surf = - \int_\surf \delta w \, g \, \dd \surf \qquad  \forall \, \{\delta w, \delta \bm{\phi} \} \in \mathit{D}^{p-1}(\surf) \times [\mathit{D}^{p-1}(\surf)]^2 \, , 
        \label{eq:QFSRM}
    \end{align}
where $\RT^{p-1}(\surf)$ is the $\Hd{,\surf}$-conforming Raviart-Thomas finite element space \cite{Raviart} of order $p-1$. The latter is equipped with a commuting interpolant \cite{DEMKOWICZ2005267,Demkowicz2000} in the exact de-Rham complex \cite{PaulyDeRham,arnold_complexes_2021} such that
\begin{subequations}
    \begin{align}
    \bm{R}\nabla \Pi_g w &= \Pi_d \bm{R} \nabla w \, , && w \in \Hone(\surf) \, , \\ 
    \di \Pi_d \vb{q} &= \Pi_o \di \vb{q} \, , && \vb{q} \in \Hd{,\surf} \, ,
\end{align}
\end{subequations}
holds for sufficiently smooth functions.
The complex is depicted in \cref{fig:com_rt}.
\begin{figure}
		\centering
        \input{figs/com_rt}
		\caption{A two-dimensional exact de Rham complex with commuting interpolants with respect to the Raviart-Thomas space.}
  \label{fig:com_rt}
\end{figure}
For the definition of base functions of the Raviart-Thomas element we refer to \cite{Anjam2015,Zaglmayr2006,Joachim2005,sky_polytopal_2022}.
\begin{theorem}[Discrete well-posedness]
    Given the discrete sub-spaces of the product-spaces
    \begin{subequations}
        \begin{align}
            \X_h^{p}(\surf) &= \U_h^{p}(\surf) \times [\mathit{D}_h^{p-1}(\surf)]^2 \subset \Hone(\surf) \times [\Le(\surf)]^2 = \X(\surf) \, , \\
            \Y_h^{p-1}(\surf) &= \mathit{D}_h^{p-1}(\surf) \times [\mathit{D}_h^{p-1}(\surf)]^2 \subset \Le(\surf) \times [\Le(\surf)]^2 = \Y(\surf) \, ,  \\
            \Z_h^p(\surf) &= \HZ_h^p(\surf) \times \RT_h^{p-1}(\surf)  \subset \HsD{,\surf} \times \Hd{,\surf} = \Z(\surf) \, ,
        \end{align}
    \end{subequations}
    then the discrete formulations \cref{eq:TFSRM} and \cref{eq:QFSRM} are well-posed in 
    \begin{align}
        \{w^h, \bm{\phi}^h, \bm{M}^h\} &\in \X_h^{p}(\surf) \times \HZ_h^{p}(\surf) \, , \\ 
        \{w^h, \bm{\phi}^h, \bm{M}^h, \vb{q}^h\} &\in \Y^{p-1}_h(\surf) \times \Z_h^p(\surf)   \, , 
    \end{align}
    respectively. 
\end{theorem}
\begin{proof}
    For both problems the discrete well-posedness is directly inherited from the continuous one due to the existence of commuting interpolants as per Fortin's criterion \cite[Thm. 4.8]{Bra2013}. 
\end{proof}

\begin{remark}
    We note that in \cref{eq:QFSRM}, one could also employ the Brezzi-Douglas-Marini \cite{BDM} element $\BDM^{p-1}(\surf)$ of order $p-1$ instead of the Raviart-Thomas element. However, that would require to discretise the deflection $w$ in the lower order space $\mathit{D}^{p-2}(\surf)$.
\end{remark}

Due to the well-posedness of the continuous and discrete problems we directly obtain the quasi best-approximation result \cite{Bra2013}.
\begin{theorem}[Quasi best approximation]
    Let $\{w,\bm{\phi},\bm{M}\}\in \X(\surf) \times \HsD{,\surf}$ be the exact solution of \cref{eq:prob_mixed}, and $\{w^h, \bm{\phi}^h, \bm{M}^h\} \in \X_h^p(\surf) \times \HZ_h^p(\surf)$ the discrete solution, then there holds
    \begin{align}
         \norm{\{w,\bm{\phi}\}-\{w^h, \bm{\phi}^h\}}_{\X} +  \norm{\bm{M}- \bm{M}^h}_{\HsD{}}\leq  C_1 \inf_{\{\delta w,\delta \bm{\phi},\delta \bm{M}\} \in \X_h^{p} \times \HZ_h^{p}} \begin{aligned}
             (&\norm{\{w,\bm{\phi}\}-\{\delta w, \delta\bm{\phi}\}}_{\X} \\ &+ \norm{\bm{M}- \delta\bm{M}}_{\HsD{}}) \, ,
         \end{aligned}
    \end{align}
    where $C_1=C_1(t)$. Analogously, let  $\{w,\bm{\phi},\bm{M},\vb{q}\}\in \Y(\surf) \times \Z(\surf)$ be the exact solution of \cref{eq:fourfield}, and $\{w^h, \bm{\phi}^h, \bm{M}^h, \vb{q}^h\} \in \Y_h^{p-1}(\surf)\times \Z_h^p(\surf)$ the discrete solution, then the approximation satisfies
    \begin{align}
        \norm{\{w,\bm{\phi}\}-\{w^h, \bm{\phi}^h\}}_{\Y} +  \norm{\{\bm{M},\vb{q}\}- \{\bm{M}^h,\vb{q}^h\}}_{\Z_t}\leq  C_2 \inf_{\{\delta w,\delta\bm{\phi},\delta\bm{M},\delta \vb{q}\} \in \Y_{h}^{p-1} \times \Z_h^{p}} \begin{aligned}
             (&\norm{\{w,\bm{\phi}\}-\{\delta w, \delta\bm{\phi}\}}_{\Y} \\ &+ \norm{\{\bm{M},\vb{q} \}- \{\delta\bm{M} , \delta \vb{q} \} }_{\Z_t}) \, ,
         \end{aligned}
    \end{align}
        where $C_2$ is independent of the thickness $t$.
    \label{th:cea}
\end{theorem}
We are now in a position to compute a priori error estimates. 
\begin{theorem}[Convergence of the three-field formulation]
    Assume the exact solution has the regularity $w\in H^{p+1}(\surf)$, $\bm{\phi}\in [H^p(\surf)]^2$, and $\bm{M}\in [H^{p}(\surf)]^{2\times 2}$, $\Di \bm{M}\in [H^p(\surf)]^2$. Then the discrete formulation \cref{eq:TFSRM} exhibits the following convergence rate for a uniform triangulation 
    \begin{align}
         \norm{\{w,\bm{\phi}\}-\{w^h, \bm{\phi}^h\}}_{\X} +  \norm{\bm{M}- \bm{M}^h}_{\HsD{}} &\leq  C h^p \,( |w|_{\H^{p+1}}+ | \bm{\phi}|_{\H^p} + |\bm{M}|_{\H^{p}} + |\Di \bm{M} |_{\H^p} ) ,
    \end{align} 
    with $C = C(t)$.
\end{theorem}
\begin{proof}
    Using the approximation error in \cref{th:cea} one finds via the interpolants from \cite{Demkowicz2000,hu_family_2014} with their approximation and commutative properties
    \begin{align}
        (\norm{\{w,\bm{\phi}\}-\{w^h, \bm{\phi}^h\}}_{\X} +  \norm{\bm{M}- \bm{M}^h}_{\HsD{}})^2 &\leq  C\inf_{\{\delta w,\delta \bm{\phi},\delta \bm{M}\} \in \X_h^{p} \times \HZ_h^{p}} \begin{aligned}
             (&\norm{\{w,\bm{\phi}\}-\{\delta w, \delta\bm{\phi}\}}_{\X} \\ &+ \norm{\bm{M}- \delta\bm{M}}_{\HsD{}})^2
         \end{aligned} \notag\\
        &\overset{\phantom{a}}{\leq}  C ( \norm{\{w,\bm{\phi}\}-\{\Pi_g^p w, \Pi^{p-1}_o \bm{\phi}\}}_{\X} +  \norm{\bm{M}- \Pi_s^p \bm{M}}_{\HsD{}} )^2
        \notag \\
        &\overset{Y}{\leq}  C ( \norm{\{w,\bm{\phi}\}-\{\Pi_g^p w, \Pi^{p-1}_o \bm{\phi}\}}_{\X}^2 +  \norm{\bm{M}- \Pi_s^p \bm{M}}_{\HsD{}}^2 ) \notag \\
        &\overset{\phantom{a}}{=}  C ( \norm{w-\Pi_g^p w}_{\Hone}^2 +\norm{\bm{\phi}-\Pi_o^{p-1} \bm{\phi}}_{\Le}^2 
             \notag \\
             &\qquad +\norm{\bm{M}- \Pi_s^p \bm{M}}_{\Le}^2 + \norm{\Di \bm{M}- \Di \Pi_s^p \bm{M}}_{\Le}^2
         )
        \notag \\
        &\overset{\phantom{a}}{=}  C ( \norm{w-\Pi_g^p w}_{\Hone}^2 +\norm{\bm{\phi}-\Pi_o^{p-1} \bm{\phi}}_{\Le}^2 
             \notag \\
             &\qquad +\norm{\bm{M}- \Pi_s^p \bm{M}}_{\Le}^2 + \norm{ (\id - \Pi_o^{p-1}) \Di \bm{M}}_{\Le}^2
         )
        \notag \\
        &\overset{\phantom{a}}{\leq}  C ( h^{2p}|w|_{\H^{p+1}}^2 + h^{2p}| \bm{\phi}|_{\H^p}^2 + h^{2p}|\bm{M}|^2_{\H^{p}} + h^{2p}|\Di \bm{M} |^2_{\H^p} )
        \notag \\
        &\overset{\phantom{a}}{\leq}  Ch^{2p}  ( |w|_{\H^{p+1}}^2 + | \bm{\phi}|_{\H^p}^2 + |\bm{M}|^2_{\H^{p}} + |\Di \bm{M} |^2_{\H^p} )  \, , 
    \end{align}
    where we used Young's inequality.
\end{proof}
\begin{theorem}[Convergence of the quad-field formulation]
    Assume the exact solution has the regularity $w\in H^{p+1}(\surf)$, $\bm{\phi}\in [H^p(\surf)]^2$, $\bm{M}\in [H^{p}(\surf)]^{2\times 2}$, $\Di \bm{M}\in [H^p(\surf)]^2$, $\vb{q}\in [H^p(\surf)]^2$, and $\di \vb{q}\in H^p(\surf)$. Then the discrete formulation \cref{eq:QFSRM} satisfies the following convergence rate
    \begin{align}
         \norm{\{w,\bm{\phi}\}-\{w^h, \bm{\phi}^h\}}_{\Y} +  \norm{\{\bm{M},\vb{q}\}- \{\bm{M}^h,\vb{q}^h\}}_{\Z_t} &\leq  C h^p \,\begin{aligned}(& |w |_{\H^p} + | \bm{\phi}|_{\H^p} + |\bm{M}|_{\H^{p}} + |\Di \bm{M} |_{\H^p}   \\&+(1 + t) |\vb{q} |_{\H^p} + |\di \vb{q} |_{\H^p} ) ,
         \end{aligned}
    \end{align} 
    under a uniform triangulation with $C$ independent of the thickness $C \neq C(t)$.
\end{theorem}
\begin{proof}
    We employ \cref{th:cea} to derive
    \begin{align}
        (\norm{\{w,\bm{\phi}\}-\{w^h, \bm{\phi}^h\}}_{\Y} +  &\,\norm{\{\bm{M},\vb{q}\}- \{\bm{M}^h,\vb{q}^h\}}_{\Z_t})^2 \leq  C \inf_{\{\delta w,\delta\bm{\phi},\delta\bm{M},\delta \vb{q}\} \in \Y_{h}^{p-1} \times \Z_h^{p}} \begin{aligned}
             (&\norm{\{w,\bm{\phi}\}-\{\delta w, \delta\bm{\phi}\}}_{\Y} \\ &+ \norm{\{\bm{M},\vb{q} \}- \{\delta\bm{M} , \delta \vb{q} \} }_{\Z_t})^2
             \end{aligned}
        \notag \\
        &\overset{\phantom{a}}{\leq}  C
             (\norm{\{w,\bm{\phi}\}-\{\Pi_o^{p-1}w, \Pi_o^{p-1}\bm{\phi}\}}_{\Y} + \norm{\{\bm{M},\vb{q} \}- \{\Pi_s^p\bm{M} , \Pi_d^{p-1} \vb{q} \} }_{\Z_t})^2
        \notag \\
        &\overset{Y}{\leq}  C ( \norm{\{w,\bm{\phi}\}-\{\Pi_o^{p-1}w, \Pi_o^{p-1}\bm{\phi}\}}_{\Y}^2 +  \norm{\{\bm{M},\vb{q} \}- \{\Pi_s^p\bm{M} , \Pi_d^{p-1} \vb{q} \} }_{\Z_t}^2 )
        \notag \\
        &\overset{\phantom{a}}{=}C ( \norm{w - \Pi_o^{p-1}w}_{\Le}^2 + \norm{\bm{\phi} - \Pi_o^{p-1} \bm{\phi}}_{\Le}^2 + \norm{\bm{M} - \Pi_s^p\bm{M}}_{\Le}^2  
        \notag \\
        & \qquad +\|(\Di \bm{M} - \Di \Pi_s^p \bm{M})-(\vb{q} - \Pi_d^{p-1}\vb{q})\|^2_{\Le}+ t^2 \norm{\vb{q}-\Pi_d^{p-1}\vb{q}}_{\Le}^2  
        \notag \\ & \qquad + \norm{\di \vb{q} - \di \Pi_d^{p-1}\vb{q}}^2_{\Le} )\notag\\
        &\overset{CS}{\leq}  C ( h^{2p}|w |_{\H^p}^2 + h^{2p}| \bm{\phi}|_{\H^p}^2 + \norm{\bm{M} - \Pi_s^p\bm{M}}_{\Le}^2 +\|(\id - \Pi_o^{p-1}) \Di \bm{M}\|^2_{\Le} 
        \notag \\ & \qquad  + \norm{\vb{q} - \Pi_d^{p-1}\vb{q}}^2_{\Le} + t^2 \norm{\vb{q}-\Pi_d^{p-1}\vb{q}}_{\Le}^2 + \norm{(\id - \Pi_o^{p-1})\di \vb{q} }^2_{\Le} )
        \notag \\ 
        &\overset{\phantom{a}}{\leq}  C ( h^{2p}|w |_{\H^p}^2 + h^{2p}| \bm{\phi}|_{\H^p}^2 + h^{2p}|\bm{M}|^2_{\H^{p}} 
       + h^{2p}|\Di \bm{M} |^2_{\H^{p}}  \notag \\ & \qquad  + h^{2p}|\vb{q} |^2_{\H^p} + t^2 h^{2p}|\vb{q} |^2_{\H^p} + h^{2p}|\di \vb{q} |^2_{\H^p} )
        \notag \\
        &\overset{\phantom{a}}{\leq}  C h^{2p} ( |w |_{\H^p}^2 + | \bm{\phi}|_{\H^p}^2 + |\bm{M}|^2_{\H^{p}} + |\Di \bm{M} |^2_{\H^p}    +(1 + t^2) |\vb{q} |^2_{\H^p} + |\di \vb{q} |^2_{\H^p} ),
    \end{align}
    where we again used Young's inequality and the interpolants from \cite{Demkowicz2000} and \cite{hu_family_2014}, and applied the Cauchy-Schwarz inequality.
\end{proof}
\begin{remark}
    We note that for certain smoothness of the data, the convergence estimates can be further improved using the Aubin-Nitsche technique, as our numerical example demonstrates, see \cref{sec:shear_example}.   
\end{remark}

\begin{remark}
    In the three-field formulation one can eliminate the rotations $\bm{\phi} \in [\mathit{D}^{p-1}(\surf)]^2$ element-wise with static condensation, since the basis is discontinuous. Although the quad-field formulation introduces another variable for the shear stress $\vb{q} \in \RT^{p-1}(\surf)$, it elevates the deflections to the discontinuous space $w \in \mathit{D}^{p-1}(\surf)$, such both the rotations $\bm{\phi}$ and the deflections $w$ can now be statically condensated element-wise. Applying the latter implies that for both formulations the global system is solved only in two-variables, $\{w,\bm{M}\}$ or $\{\bm{M}, \vb{q}\}$, respectively.    
\end{remark}

\subsection{Polytopal template for $\HZ^p$}
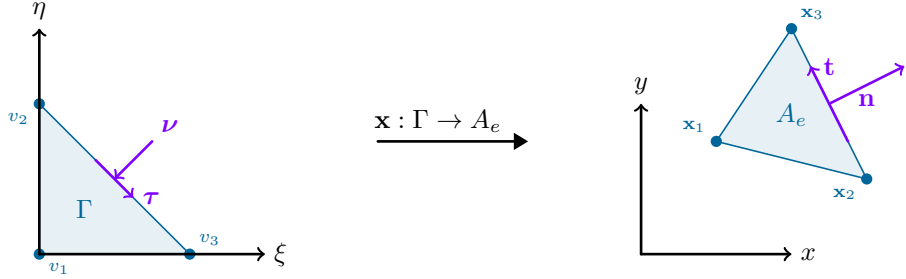
\begin{figure}
		\centering
		\input{figs/trimap}
		\caption{Barycentric mapping of the reference triangle to an element in the physical domain.}
		\label{fig:trimap}
	\end{figure}

In this work, we introduce an efficient approach of constructing $\HZ^p(\Gamma)$ base functions on the reference element (see \cref{fig:trimap}) 
\begin{align}
    \Gamma = \{ (\xi, \eta) \in [0,1] \; | \; \xi + \eta \leq 1 \} \, , 
\end{align}
via the polytopal template methodology \cite{sky_polytopal_2022,sky_higher_2023},
and mapping them to the physical element by a polytope-specific transformation, which allows for curved element geometries.
We mention that an alternative mapping approach for the Arnold-Winther element can be found in \cite{kirby_general_2018,aznaran_transformations_2021}, although it does not treat the question of curved geometries.   
We restrict ourselves to $p \geq 3$, which allows for a straight-forward construction of the basis.

In order to introduce the construction, we first decompose an arbitrary subspace for $\Hone$-discretisations $\U^p$ into its parts associated with the polytopes of the triangle.
\begin{definition}[Triangle $\U^p(\surf)$-polytopal spaces]
Each polytope is associated with a space of base functions as follows:
\begin{itemize}
    \item Each vertex is associated with the space of its respective base function $\ver^p_i$. As such, there are three spaces in total $i \in \{1,2,3\}$ and each one is of dimension one, $\dim \ver^p_i = 1 \quad \forall \, i \in \{1,2,3\}$.
    \item For each edge there exists a space of edge functions $\edge^p_j$ with the multi-index $j \in \mathcal{J} = \{(1,2),(1,3),(2,3)\}$.
    The dimension of each edge space is given by $\dim \edge^p_j = p-1$.
    \item Lastly, the cell is equipped with the space of cell base functions $\cell_{123}^p$ with $\dim \cell_{123}^p = (p-2)(p-1)/2$.
\end{itemize}
The association with a respective polytope reflects the support of the trace operator for $\Hone$-spaces.
\end{definition}
The polynomial space for the $\HZ^p(\surf)$ element reads 
\begin{align}
    &\mathit{S}^p(\surf) = \Po^p(\surf) \otimes \Sym(2) \, , && p \geq 3 \, .
\end{align}
The Hu-Zhang element is defined in \cite{hu_family_2014} with a polytopal framework on each physical element using Lagrangian base functions $\Lag^p(\surf) = \spa \{n_i\}$ such that on vertices the base functions are given by 
\begin{align}
    &\bm{\rho}^p_{ij} = n^p_i\bm{S}_j  \, , && n^p_i \in \Lag^p(\surf) \cap \ver_i^p(\surf) \, , && \bm{S}_j \in \{\bm{S} \in \Sym(2) \; | \; \| \bm{S} \| = 1\} \, . 
\end{align}
The same method is used to construct the internal cell functions with $n_i^p \in \Lag^p(\surf) \cap \cell^p(\surf) = \Lag_0^p(\surf)$. On each edge, the Hu-Zhang element defines two types of base functions with connectivity 
\begin{align}
    &\bm{\rho}^p_{ij} = n^p_i\bm{S}_j \, , && n^p_i \in \Lag^p(\surf) \cap \edge_i^p(\surf) \, , && \bm{S}_j \in \{ \sym(\vb{t} \otimes \vb{n}), \, \vb{n} \otimes \vb{n} \} \, ,
\end{align}
which preserve the symmetric normal continuity of the space, and one edge-cell function type
\begin{align}
    &\bm{\rho}^p_{i} = n^p_i \vb{t} \otimes \vb{t} \, , && n^p_i \in \Lag^p(\surf) \cap \edge_i^p(\surf) \, ,
\end{align}
for the tangential-tangential component, which may jump between elements.
Clearly, the construction is linearly independent, since each base function is multiplied with a set of linearly independent tensors, such that the linear independence of the Lagrangian space is inherited on the tensorial level.

\begin{figure}
		\centering
		\input{figs/template.tex}
		\caption{Template tensors for the reference Hu-Zhang triangle element on their corresponding polytope.}
		\label{fig:template}
	\end{figure}
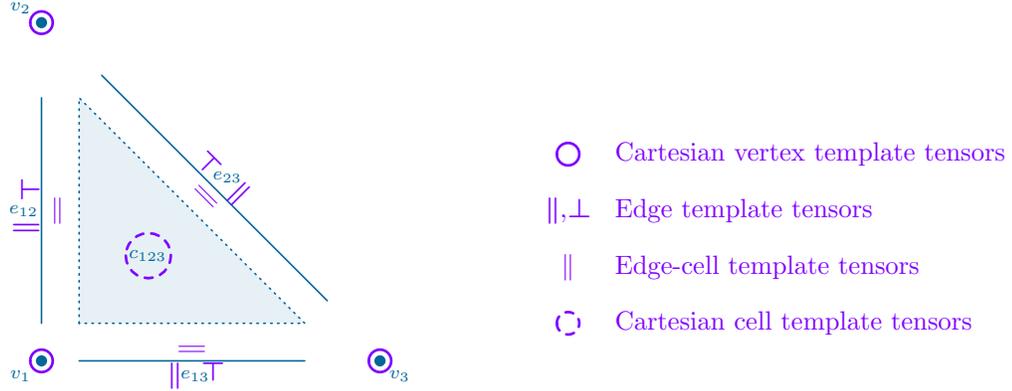

Following a similar approach, we define a polytopal template \cite{sky_polytopal_2022,sky_higher_2023} for the Hu-Zhang element on the reference triangle
\begin{align}
    \tem_1 &= \tem_2 = \tem_3 = \tem_{123} = \{\vb{e}_1 \otimes \vb{e}_1 , \sym(\vb{e}_1 \otimes \vb{e}_2), \vb{e}_2 \otimes \vb{e}_2 \} \, , \notag \\ \tem_{12} &= \{ \vb{e}_2 \otimes \vb{e}_2 , \sym (\vb{e}_1 \otimes \vb{e}_2), \vb{e}_1 \otimes \vb{e}_1 \} \, , \notag \\ \tem_{13} &= \{ \vb{e}_1 \otimes \vb{e}_1 , -\sym (\vb{e}_1 \otimes \vb{e}_2), \vb{e}_2 \otimes \vb{e}_2 \} \, , \notag \\ \tem_{23} &= \{ (\vb{e}_1 - \vb{e}_2) \otimes (\vb{e}_1 - \vb{e}_2) , \sym [(\vb{e}_2 - \vb{e}_1) \otimes (\vb{e}_1 + \vb{e}_2)], (\vb{e}_1 + \vb{e}_2) \otimes (\vb{e}_1 + \vb{e}_2) \} \, .
\end{align}
On the vertices and in the cell the template tensors are simply a normalised basis for the space of two-dimensional symmetric tensors $\Sym(2) = \spa[\vb{e}_1 \otimes\vb{e}_1, \sym(\vb{e}_1 \otimes \vb{e}_2), \vb{e}_2 \otimes \vb{e}_2]$. 
On edges, the template is given by the dyadic products of the non-normalised normal $\bm{\nu}$ and tangent $\bm{\tau}$ vectors belonging to each edge, respectively. 
The total template is given by 
\begin{align}
    \tem = \{\tem_1,\tem_2,\tem_3,\tem_{12},\tem_{13},\tem_{23},\tem_{123}\} \, .
\end{align}
The association of the template tensors with the polytopes of the reference element is depicted in \cref{fig:template}.
The Hu-Zhang space on the reference element is now given by
\begin{align}
    &\HZ^p(\Gamma) = \left \{ \bigoplus_{i=1}^3 \ver_i^p(\Gamma) \otimes \tem_i \right \} \oplus \left \{ \bigoplus_{j\in \mathcal{J}} \edge_j^p(\Gamma) \otimes \tem_j \right \} \oplus \cell_{123}^p(\Gamma) \otimes \tem_{123} \, , && \mathcal{J} = \{(1,2),(1,3),(2,3)\} \, , && p \geq 3 \, ,
\end{align}
where $\ver^p(\Gamma)_i$ are spaces of vertex-, $\edge^p(\Gamma)_j$ are spaces of edge-, and $\cell^p(\Gamma)_{123}$ is the space of cell base functions on the reference element of power $p \geq 3$. The summation over $j$ is understood in the sense of multi-indices.

We can now define the base functions of the Hu-Zhang element on the reference domain using some $\Hone$-conforming subspace equipped with the scalar base functions $\spa\{n_i\}$ such that $n_i = n_i(\xi,\eta)$.
\begin{definition}[Triangle $\HZ^p(\Gamma)$ base functions]
	The base functions of the Hu-Zhang element with $p \geq 3$ are defined on their respective polytopes as follows.
 \begin{itemize}
     \item On each vertex $v_i$ the base functions read
     \begin{align}
         \bm{\varrho}(\xi,\eta) &= n \bm{T} \, , && n \in \ver^p_i(\Gamma) \, , && \bm{T} \in \tem_i \, . 
     \end{align}
     \item On each edge $e_{ij}$ with vertices $v_i$ and $v_j$ the base functions are given by
         \begin{align}
             \bm{\varrho}(\xi,\eta) &= n \bm{T} \, , && n \in \edge^p_{ij}(\Gamma) \, , && \bm{T} \in \{\bm{T} \in \tem_{ij} \; | \;  \bm{T}\bm{\nu} \neq 0 \} \, ,
     \end{align}
     where $\bm{\nu}$ is the normal vector on the respective edge. As such, each scalar base function $n$ defines one symmetric tangent-normal base function $\bm{\varrho}_{\tau \nu}$ and one normal-normal base function $\bm{\varrho}_{\nu \nu}$ on each edge. 
     \item The cell base functions read
     \begin{subequations}
         \begin{align}
             \bm{\varrho}(\xi,\eta) &= n \bm{T} \, , && n \in \edge^p_{j}(\Gamma) \, , && \bm{T} \in \{\bm{T} \in \tem_{j} \; | \;  \bm{T}\bm{\nu} = 0 \} \, , && j \in  \{(1,2),(1,3),(2,3)\}  \, , \\
             \bm{\varrho}(\xi,\eta) &= n \bm{T} \, , && n \in \cell^p_{123}(\Gamma) \, , && \bm{T} \in \tem_{123} \, ,
     \end{align}
     \end{subequations}
     where the first three definitions are edge-cell base functions $\bm{\varrho}_{\tau \tau}$.
 \end{itemize}
 \end{definition}
The polytopal construction allows for an arbitrary choice of an $\Hone$-conforming subspace in the construction of the Hu-Zhang element. However, in this work, we rely on the Barycentric-, Legendre- and scaled integrated Legendre polynomials \cite{Zaglmayr2006,Joachim2005}
\begin{align}
    \lambda_1(x,y) &= 1-x-y \, , & \lambda_2(x,y) &= y \, , & \lambda_3(x,y) &= x \, , \\  
    l^0(x) &= 1 \, , & l^1(x) &= x \, , & l^p(x) &= (2p - 1)x \, l^{p-1}(x) - (p-1)l^{p-2}(x) \, , \\
    L^1_s(x,t) &= x \, , & L^2_s(x,t) &= \dfrac{1}{2} (x^2 - t^2) \, , & L^p_s(x,t) &= (2p - 3) x \, L^{p-1}_s(x,t) - (p-3) t^2 L^{p-2}_s(x,t) \, , 
\end{align}
where the scaled integrated Legendre polynomials are related to integrated Legendre polynomials via
\begin{align}
    &L_s^p(x, t) = (t)^{p} L^p \left ( \dfrac{x}{t} \right ) \, , && L^p(x) = \int_{-1}^x l^{p-1}(x) \, \dd x \, , && p \geq 2 \, .
\end{align}
Consequently, our construction is applicable to the hp-finite element method \cite{Zaglmayr2006,VOS20105161,KOPP2022115575}.
\begin{definition}[Hu-Zhang-Legendre base functions]
    In the following we define the Hu-Zhang-Legendre base functions on their respective polytopes.
    \begin{itemize}
        \item On each vertex $v_i$ the base functions are constructed using the corresponding barycentric coordinate
        \begin{align}
            &\bm{\varrho}_{ij} = \lambda_i \bm{T}_j \, , && \bm{T}_j \in \tem_i \, ,
        \end{align}
        such that there are three base functions on each vertex.
        \item On every edge $e_{ij}$ with $(i,j) \in \mathcal{J}$ and a corresponding normal vector $\bm{\nu}$ we apply the scaled integrated Legendre polynomials in the construction
        \begin{align}
             & \bm{\varrho}_{ijk}^p = L_s^p(\lambda_i, \lambda_j) \bm{T}_k  \, , && \bm{T}_k  \in \{\bm{T} \in \tem_{ij} \; | \;  \bm{T}\bm{\nu} \neq 0 \} \, , && 
             (i,j) \in \{(1,2),(1,3),(2,3)\}  \, ,
        \end{align}
        with $p \geq 2$. There are two base functions for each polynomial power on each edge. 
        \item Finally, cell base functions are given by
        \begin{subequations}
            \begin{align}
        \bm{\varrho}_{ijk}^p &= L_s^p(\lambda_i, \lambda_j) \bm{T}_k  \, , && \bm{T}_k  \in \{\bm{T} \in \tem_{ij} \; | \;  \bm{T}\bm{\nu} = 0 \} \, , && 
             (i,j) \in \{(1,2),(1,3),(2,3)\}  \, , \\
        \bm{\varrho}_{ik}^p &= \eta \, L_s^p(\lambda_1, \lambda_3) \, l^k(2 \xi - 1) \bm{T}_i \, , && \bm{T}_i \in \tem_{123} \, ,
        \end{align}
        \end{subequations}
        where the first term represents edge-cell base functions, such that $\bm{\nu}$ is the respective normal vector of said edge. The second row defines pure-cell base functions.
    \end{itemize}
\end{definition}

\subsection{Polytopal transformations}
Observe that on each edge, the base functions are simply given by the multiplication of a scalar base functions with dyadic products of the tangent and normal vectors of the edge. As such, we can define one-to-one maps to these tensors using double contractions with fourth order tensors.
On each edge of a physical element with tangent and normal vectors $\vb{t}$ and $\vb{n}$ we define the tensors 
\begin{subequations}
    \begin{align}
        &\mathbb{T}^{tt}: \bm{\tau} \otimes \bm{\tau} \to \vb{t} \otimes \vb{t} \, , & \mathbb{T}^{tt} &= \dfrac{1}{\norm{\bm{\tau}}^4} \vb{t} \otimes \vb{t} \otimes \bm{\tau} \otimes \bm{\tau} \, , \\
        &\mathbb{T}^{tn}: \sym(\bm{\tau} \otimes \bm{\nu}) \to \sym(\vb{t} \otimes \vb{n}) \, , & \mathbb{T}^{tn} &= \dfrac{1}{\norm{\bm{\tau}}^4} \sym(\vb{t} \otimes \vb{n}) \otimes \sym(\bm{\tau} \otimes \bm{\nu}) \, , \\
        &\mathbb{T}^{nn}: \bm{\nu} \otimes \bm{\nu} \to \vb{n} \otimes \vb{n} \, , & \mathbb{T}^{nn} &= \dfrac{1}{\norm{\bm{\tau}}^4} \vb{n} \otimes \vb{n} \otimes \bm{\nu} \otimes \bm{\nu} \, ,
    \end{align}
\end{subequations}
where $\bm{\tau}$ and $\bm{\nu}$ are the tangent and normal vectors on the corresponding reference edge of the element.
Due to the orthogonality of the tensorial basis 
\begin{align}
    \con{\bm{\tau} \otimes \bm{\tau}}{\sym(\bm{\tau} \otimes \bm{\nu})} = \con{\bm{\tau} \otimes \bm{\tau}}{\bm{\nu} \otimes \bm{\nu}} = \con{\sym(\bm{\tau} \otimes \bm{\nu})}{\bm{\nu} \otimes \bm{\nu}} = 0 \, ,
\end{align}
we can combine the three fourth-order transformation tensors into one transformation tensor for each edge
\begin{align}
        &\mathbb{T}= \mathbb{T}^{tt} + \mathbb{T}^{tn} + \mathbb{T}^{nn} \, .
\end{align}
The vertex base functions do not require any transformation since full symmetric-continuity is imposed at vertices and the Cartesian basis is global. As for the cell base functions, these do not affect the continuity of the construction since their underlying scalar base functions vanish on all edges of each element, such that no transformation is needed in order to maintain conformity.
We summarise the transformation in the following definition.
\begin{definition}[Transformations from the reference to the physical element]
    Only the base functions on edges require a transformation. All other base functions are mapped by the identity operator $ \bm{\rho} = \id(\bm{\varrho}) = \mathbb{J}\bm{\varrho} = \bm{\varrho}$. 
    \begin{itemize}
        \item On each edge $e_{ij}$ with $(i,j) = \{(1,2),(1,3),(2,3)\}$ equipped with the tangent and normal vectors $\bm{\tau}$ and $\bm{\nu}$, such that $\vb{t} = \bm{J} \bm{\tau}$ and $\vb{n} = (\cof \bm{J}) \bm{\nu}$
        the transformation tensor is given by
        \begin{align}
            \mathbb{T}_{ij}: \bm{\varrho}_{ijk}^p \to \bm{\rho}_{ijk}^p \, , && \mathbb{T}_{ij} = \dfrac{1}{\norm{\bm{\tau}}^4}(\vb{t} \otimes \vb{t} \otimes \bm{\tau} \otimes \bm{\tau} + \sym(\vb{t} \otimes \vb{n}) \otimes \sym(\bm{\tau} \otimes \bm{\nu}) + \vb{n} \otimes \vb{n} \otimes \bm{\nu} \otimes \bm{\nu}) \, .
        \end{align}
        Therefore, edge base functions on the physical edge of the element are generated by the double contraction of the corresponding transformation tensor with the base functions of the reference edge $\bm{\rho}_{ijk}^p = \mathbb{T}_{ij} \bm{\varrho}_{ijk}^p$. 
    \end{itemize}
    \label{de:trans}
\end{definition}
The divergence of the base functions follows via the chain-rule
\begin{align}
    \Di_x \bm{\rho} = (\mathbb{T} \bm{\varrho})_{,i}(\bm{J}^{-T} \vb{e}_i) \, .
\end{align}
In the simpler case of affine transformations, such that the Jacobian is a constant matrix $\bm{J} = const$, the divergence can be expressed as  
\begin{align}
    \Di_x \bm{\rho} = (\mathbb{T} \bm{\varrho} ) \nabla_x = [\mathbb{T} (n \bm{T})]( \bm{J}^{-T} \nabla_\xi) = (\mathbb{T} \,  \bm{T} \bm{J}^{-T})\nabla_\xi n \, .
\end{align}
 
\subsection{Boundary conditions}
In the following we assume the boundary data is smooth enough, such that point evaluation is possible. 
Consequently, we can impose the Dirichlet data at vertices via the functionals
\begin{align}
         &\con{\vb{e}_1 \otimes \vb{e}_1}{\bm{M}}\at_{v_i} \, , && \con{\sym(\vb{e}_1 \otimes \vb{e}_2)}{\bm{M}}\at_{v_i} \, , && \con{\vb{e}_2 \otimes \vb{e}_2}{\bm{M}}\at_{v_i} \, ,  
    \end{align}
for each vertex $v_i$.  
\begin{remark}
    The $\Hd{}$-space is characterised by the trace operator $\trnp(\cdot) = \con{\vb{n}}{\cdot}$ which considers the normal projection of fields on the boundary of the domain. Therefore, in the case of matrix-valued $\HsD{}$-fields, data does not have to be prescribed with respect to the tangent-tangent component. In order to leave the tangent-tangent component of the Hu-Zhang basis free even at vertices, it is required to identify its corresponding base function. This motivates the construction of the Hu-Zhang vertex base functions as 
    \begin{align}
        &\bm{\rho}_{ijk} = \sym(n_i \vb{d}_j \otimes \vb{d}_k) \, , && n_i \in \ver_i(\surf) \, , && \vb{d}_j = \bm{Q} \vb{e}_j \, , && \spa\{ \sym( \vb{d}_j \otimes \vb{d}_k) \} = \Sym(2) \, ,
    \end{align}
    for each vertex $v_i$ where $\{\vb{d}_j, \vb{d}_k\}$ represent an orthonormal basis such that $\vb{d}_1 \parallel \vb{t}$ and $\bm{Q} = \begin{bmatrix} \vb{d}_1 & \vb{d}_2 \end{bmatrix} \in \SO(2)$ is a rotation matrix, see \cref{fig:curved_ver}.
    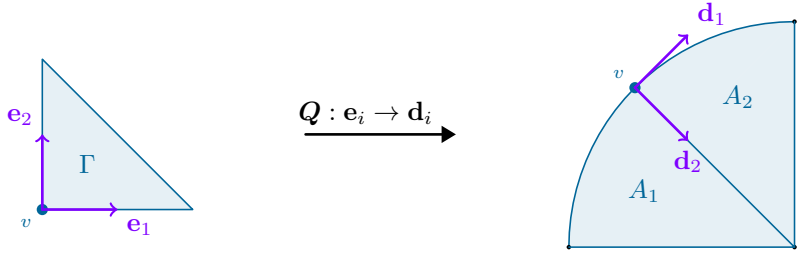
\begin{figure}
		\centering
		\input{figs/curved.tex}
		\caption{Mapping of the Cartesian basis belonging to vertex $v$ to a rotated orthonormal basis on the boundary of the domain, such that the tangent-tangent base function can be clearly defined as $\bm{\rho}_{111}= n_1\vb{d}_1 \otimes \vb{d}_1$ and its corresponding degree of freedom left unprescribed.}
		\label{fig:curved_ver}
	\end{figure}
    Thus, the evaluation of Dirichlet data at the vertex is carried out via the functionals
    \begin{align}
        &\con{\sym(\vb{d}_1 \otimes \vb{d}_2)}{\bm{M}}\at_{v_i} \, , && \con{\vb{d}_2 \otimes \vb{d}_2}{\bm{M}}\at_{v_i} \, ,
    \end{align}
    for some tensor field $\bm{M}$. Clearly, its tangent-tangent component $\vb{d}_1 \otimes \vb{d}_1$ is left free.
    If the Dirichlet boundary is $\C^1(\curv_D)$-continuous, then the tangent vector $\vb{d}_1$ is unique at all points. If the boundary is approximated by a $\C^0(\curv_D)$-continuous discretisation, then one can employ the average of two normal vectors belonging to neighbouring elements on the boundary to construct the $\vb{d}_i$-basis 
    \begin{align}
        &\vb{n}_* = \norm{\vb{n}_1}^{-1}\vb{n}_1 + \norm{\vb{n}_2}^{-1} \vb{n}_2 \, , && \vb{d}_2 = \norm{\vb{n}_*}^{-1} \vb{n}_* \, , &&
        \vb{d}_1 = \bm{R}^T \vb{d}_2 \, , 
    \end{align}
    where $\bm{R}$ is the ninety-degree rotation matrix. 
\end{remark}

If the data is fully known at every edge $e_{ij}$ corresponding with a curve $\curv_{ij} \subset \curv_D$, then one can simply apply a localised $\Hone$-projection of all non-tangent-tangent components  
\begin{align}
    \con{\bm{\rho}_k}{ \bm{M}}_{\Hone} = \con{\bm{\rho}_k}{\widetilde{\bm{M}}}_{\Hone}  \quad \forall \, \bm{\rho}_k  \in \HZ(\curv_{ij}) \setminus \left\{ \bm{\rho} \in \HZ(\curv_{ij}) \; | \; \con{\vb{t} \otimes \vb{t}}{\bm{\rho}}\at_{\curv_{ij}} \neq 0 \right \}  \, .
\end{align}
Note that vertex-values are prescribed a priori and as such, are considered Dirichlet-data for the above variational problem. 
Consequently, the semi-norm $|\cdot|_{\Hone}$ is norm-equivalent to $\|\cdot\|_\Le$ on the edge due to $\{\bm{M} \vb{n},\widetilde{\bm{M}}\vb{n}\} \in [\Honez(\curv_{ij})]^2$, such that no directional-derivatives of the base functions are required. Thus, one can define the equivalent linear problem as
\begin{align}
    k_{kl} &= \int_{\curv_{ij}}\con{\bm{\rho}_k}{\bm{\rho}_l} \, \dd \curv \, , && f_k = \int_{\curv_{ij}} \con{\bm{\rho}_k}{\widetilde{\bm{M}}}  \, \dd \curv \, , 
\end{align}
where $k_{kl}$ are the components of the stiffness matrix and $f_k$ is the corresponding right-hand-side.
Note that unless the vertex-data is split between tangential-tangential and non-tangential-tangential, the tangential-tangential component at the vertices is also incorporated into the right-hand-side.

\section{Numerical examples}
In the following we compute examples and compare between the primal (PRM), the MITC, the TDNNS, and our newly introduced mixed formulations with symmetric finite elements for the bending moments, designated here triple-field-symmetric-Reissner-Mindlin (TFSRM) for \cref{eq:TFSRM} and the quad-field-symmetric-Reissner-Mindlin (QFSRM) for \cref{eq:QFSRM}. The computations were performed in the open source finite element software NGSolve\footnote{www.ngsolve.org}
\cite{Sch2014,Sch1997}, and the implementation is available as additional material to this work\footnote{https://github.com/Askys/NGSolve\_HuZhang\_Element}. 

We briefly introduce the discretisations of the aforementioned variational formulations
\begin{subequations}
    \begin{align}
    &\text{PRM:} && \begin{aligned} 
    &\int_\surf \dfrac{t^3}{12} \langle \sym \D \delta \bm{\phi} , \, \mathbb{D} \sym \D \bm{\phi} \rangle + \ksmu \, t \langle \nabla \delta w - \delta \bm{\phi} , \, \nabla w - \bm{\phi} \rangle \, \dd \surf   \\  &\hspace{5cm} =\int_\surf t \, \delta w \, f \, \dd \surf \, , \quad \forall\, \{\delta w, \delta \bm{\phi}\} \in \U^p(\surf) \times [\U^p(\surf)]^2 \, ,
    \end{aligned}  \\ 
    &\text{MITC:} && \begin{aligned} 
    &\int_\surf \dfrac{t^3}{12} \langle \sym \D \delta \bm{\phi} , \, \mathbb{D} \sym \D \bm{\phi} \rangle + \ksmu \, t \langle \Pi_c^{p-1} (\nabla \delta w - \delta \bm{\phi}) , \, \Pi_c^{p-1} (\nabla w - \bm{\phi}) \rangle \, \dd \surf   \\ & \hspace{5cm} =\int_\surf t \, \delta w \, f \, \dd \surf \, , \quad \forall\, \{\delta w, \delta \bm{\phi}\} \in \U^p(\surf) \times [\mathit{S}^p(\surf)]^2 \, ,
    \end{aligned}  \\ 
    &\text{TDNNS:} && \begin{aligned} 
     \int_\surf \con{\delta \bm{M}}{\mathbb{A} \bm{M}} \, \dd \surf +  \con{\Di \delta \bm{M}}{\bm{\phi}}_{\mathcal{T}}  &= 0 \quad \forall \, \delta \bm{M} \in \HHJ^{p-1}(\surf) \, , \\
      \con{\delta \bm{\phi}}{\Di \bm{M}}_{\mathcal{T}}  - \int_\surf \dfrac{\ksmu}{t^2} \con{\nabla \delta w - \delta \bm{\phi}}{\nabla w - \bm{\phi}} \, \dd \surf &= - \int_\surf \delta w \, g \, \dd \surf \quad \forall \, \{\delta w, \delta \bm{\phi}\} \in \U^p(\surf) \times \Ned^{p-1}(\surf) \, , 
    \end{aligned}  
    \end{align}
    \label{eq:forms}
\end{subequations}
where $\mathit{S}^p(\surf)$ is given by $\U^p(\surf)$ enriched with bubble functions on each element $\mathit{S}^p(\surf) = \U^p(\surf) \oplus \sum_{e} \cell^{p+1}(\surf_e)$. The operator $\Pi_c^{p-1}$ defines the interpolant into the N\'ed\'elec element $\Ned^{p-1}(\surf)$ of the first type \cite{nedelec_mixed_1980}. We note that the present definition of the MITC element is equivalent to the one presented in \cite{bathe_mitc7_1989,Bathe1990Dis}. 
\begin{remark}
    The interpolation operator into the N\'ed\'elec space $\Pi_c^{p-1}$ in the MITC formulation does not need to be applied to $\nabla w$ but rather only to $\bm{\phi}$, if commuting interpolants \cite{DEMKOWICZ2005267} are employed for the finite element spaces, such that $\nabla \Pi^p_g w = \Pi_c^{p-1} \nabla w$, where $\Pi_g^p$ is the interpolant into the discrete $\C^0(\surf)$-continuous space spanned by polynomials of order $p$.
\end{remark}
\begin{remark}
    The scalar product $\con{\Di \delta \bm{M}}{\bm{\phi}}_{\mathcal{T}}$ is to be understood in the distributional sense \cite{pechstein_tdnns_2017} and includes boundary terms on each element $\con{\Di \delta \bm{M}}{\bm{\phi}}_{\mathcal{T}} = \sum_{T \in \mathcal{T}} \langle \Di \delta\bm{M} , \, \bm{\phi} \rangle_{\Le(T)} - \langle \delta \bm{M} \vb{n} , ( \vb{t}\otimes \vb{t} )\bm{\phi}  \rangle_{\Le(\partial T)}$. 
\end{remark}

In the following examples, relative errors are measured in the $\Le$-norm
\begin{align}
    &\norm{\widetilde{w} - w^h}_{\Le} / \norm{\widetilde{w}}_{\Le} \, , && \norm{\widetilde{\bm{\phi}} - \bm{\phi}^h}_{\Le} / \norm{\widetilde{\bm{\phi}}}_{\Le} \, , && \norm{\widetilde{\bm{M}} - \bm{M}^h}_{\Le} / \norm{\widetilde{\bm{M}}}_{\Le} \, , && \norm{\widetilde{\vb{q}} - \vb{q}^h}_{\Le} / \norm{\widetilde{\vb{q}}}_{\Le} \, ,
\end{align}
where $\widetilde{\cdot}$ represent analytical solutions and $\cdot^h$ are the obtained discrete solutions.

\subsection{Shear-locking on a rectangular plate} \label{sec:shear_example}
In the following we compare the behaviour of the five plate formulations. 
We define the domain $\overline{\surf} = [0,1]^2$ and set the boundary conditions such that the deflections $w$ and rotations $\bm{\phi}$ vanish on the boundary. The embedding of the clamped boundary condition for the rotations in the different formulations is summarised in the following table.
    \begin{table}[H]
    \centering
\begin{tabular}{l| c | c}
                     & $\bm{\phi}$ & $\bm{M}$ \\  \hline
                     \vspace{-0.2cm} &  &   \\   
              TFSRM/QFSRM &       $-$      &   $\int_{\curv_N^M} \con{\delta \bm{M} \, \vb{n}}{\bm{\phi}} \, \dd \curv = 0 \, , \qquad \curv_N^M = \partial \surf$         \\ 
              \vspace{-0.2cm} &  & \\ \hline 
              \vspace{-0.2cm} &  &   \\
              TDNNS                &      $\con{\vb{t}}{\bm{\phi}}\at_{\curv_D^\phi} = 0 \, , \qquad \curv_D^\phi = \partial \surf $  &  $\int_{\curv_N^M} \con{ \delta \bm{M} \, \vb{n}}{ (\vb{n} \otimes \vb{n}) \bm{\phi}} \, \dd \curv = 0 \, , \qquad \curv_N^M = \partial \surf$        \\  \vspace{-0.2cm} &  & \\ \hline
              \vspace{-0.2cm} &  &   \\
          PRM/MITC                &   $\bm{\phi}\at_{\curv_D^\phi} = 0 \, , \qquad \curv_D^\phi = \partial \surf    $       &      $-$   
\end{tabular}
\end{table} 
In order to compare the formulations we employ the analytical solution from \cite{Guo2013}. The deflection and rotation fields can be introduced in a concise manner using the following auxiliary functions
\begin{align}
         &f_0(\alpha) =  (\alpha-1)\alpha \, , &&
         f_1(\alpha) = 5\alpha^2 - 5\alpha + 1 \, , &&
         f_2(\alpha) = 2\alpha - 1 \, ,
\end{align}
such that the analytical solution reads
\begin{subequations}
        \begin{align}
    \widetilde{w}(x,y) &= \dfrac{100}{3} [ f_0(x)^3 f_0(y)^3] - \dfrac{2t^2}{5(1-\nu)} [f_0(y)^3  f_0(x)  f_1(x) + f_0(x)^3  f_0(y) f_1(y)] \, ,  \\
    \widetilde{\bm{\phi}}(x,y) &= 100
    \begin{bmatrix}
        f_0(y)^3 f_0(x)^2 f_2(x) \\ f_0(x)^3 f_0(y)^2 f_2(y)
    \end{bmatrix}
    \, ,   \\
    \widetilde{\bm{M}}(x,y) &= \mathbb{D}_* \sym \D \widetilde{\bm{\phi}} \, , \\
    \widetilde{\vb{q}}(x,y) &= -\dfrac{\ksmu}{t^2} (\nabla \widetilde{w} - \widetilde{\bm{\phi}}) \, .
\end{align}
\end{subequations}
The corresponding forces are given by
\begin{align}
    g(x,y) = \dfrac{200 E}{1 - \nu ^ 2}  [ f_0(y) f_1(x)  [2  f_0(y)^2 + f_0(x)  f_1(y)] + f_0(x)  f_1(y)  [2 f_0(x)^2 + f_0(y)  f_1(x)] ] \, .
\end{align}
We set the material parameters to 
\begin{align}
    &E = 1 \, , && \nu = 0.3 \, , && k_s = \dfrac{5}{6} \, , 
\end{align}
and vary the thickness $t \in \{10^{-1},10^{-5}\} $ across six regular meshes with $2^k$-elements such that $k \in \{1,3,\dots,11\}$.
The resulting convergence rates of the deflection $w$, the rotations $\bm{\phi}$, the bending moments $\bm{M}$ and the shear stress $\vb{q}$ under h-refinement for the case $t = 10^{-1}$ are depicted in \cref{fig:ex12}.
\begin{figure}
    	\centering
    	\begin{subfigure}{0.48\linewidth}
    		\centering
    		\input{figs/ex1_e-1_w}
    		\caption{}
    	\end{subfigure}
    	\begin{subfigure}{0.48\linewidth}
    		\centering
    		\input{figs/ex1_e-1_p}
    		\caption{}
    	\end{subfigure}
        \begin{subfigure}{0.48\linewidth}
    		\centering
    		\input{figs/ex1_e-1_M}
    		\caption{}
    	\end{subfigure}
     \begin{subfigure}{0.48\linewidth}
    		\centering
    		\input{figs/ex1_e-1_q}
    		\caption{}
    	\end{subfigure}
    	\caption{Relative error of the cubic formulations for $t = 10^{-1}$ in the deflection $w$ (a), the rotations $\bm{\phi}$ (b) ,the bending moments $\bm{M}$ (c) and the shear stress $\vb{q}$ (d).}
    	\label{fig:ex12}
    \end{figure}
\begin{figure}
    	\centering
    	\begin{subfigure}{0.48\linewidth}
    		\centering
    		\input{figs/ex1_e-5_w}
    		\caption{}
    	\end{subfigure}
    	\begin{subfigure}{0.48\linewidth}
    		\centering
    		\input{figs/ex1_e-5_p}
    		\caption{}
    	\end{subfigure}
        \begin{subfigure}{0.48\linewidth}
    		\centering
    		\input{figs/ex1_e-5_M}
    		\caption{}
    	\end{subfigure}
     \begin{subfigure}{0.48\linewidth}
    		\centering
    		\input{figs/ex1_e-5_q}
    		\caption{}
    	\end{subfigure}
    	\caption{Relative error of the cubic formulations for $t = 10^{-5}$ in the deflection $w$ (a), the rotations $\bm{\phi}$ (b), the bending moments $\bm{M}$ (c), and the shear-stress (d).}
    	\label{fig:ex122}
    \end{figure}
For $t = 10^{-1}$ all formulations perform well. The primal and MITC forms achieve faster convergence in the rotations $\bm{\phi}$ due to their increased polynomial orders. In contrast, the TFSRM and QFSRM formulations achieve only cubic convergence in the rotations $\bm{\phi}$, but quartic convergence in the bending moments $\bm{M}$. 
Considering the case $t = 10^{-5}$ depicted in \cref{fig:ex122}, we can clearly observe shear-locking in the primal formulation for meshes with $2^1$ and $2^3$ elements although we are using cubic polynomials. In fact, its convergence rates deteriorate for all variables and the shear stress $\vb{q}$ does not converge at all. The TFSRM formulation maintains optimal convergence rates for the deflection $w$ and the rotations $\bm{\phi}$. However, it is reduced to quadratic convergence in the bending moments $\bm{M}$ and linear convergence in the shear stress $\vb{q}$. Both the QFSRM and TDNNS formulations circumvent this problem and converge optimally in all variables. 
Interestingly, the MITC formulation is optimal in all variables aside from the shear stress $\vb{q}$, for which sub-optimal convergence whatsoever is achieved. In fact, testing the MITC formulation also for $t = 10^{-3}$ results in a sub-optimal convergence rate of $\O(h^2)$. We note that similar results are observed in \cite{Bathe1990Dis}. Further, for both the MITC and TDNNS formulations we observe a loss of convergence on the finest mesh with $2^{11}$-elements. This is due to the correlation of the extremely small element size $h \ll 1$ with the extremely small thickness $t \ll 1$, such that a loss in floating-point precision causes the solution diverge. Observe that the latter is not related to the robustness of the formulations, but rather the capacity of the solver.  

\begin{figure}
    	\begin{subfigure}{0.24\linewidth}
    		\centering
    		\includegraphics[width=1\linewidth]{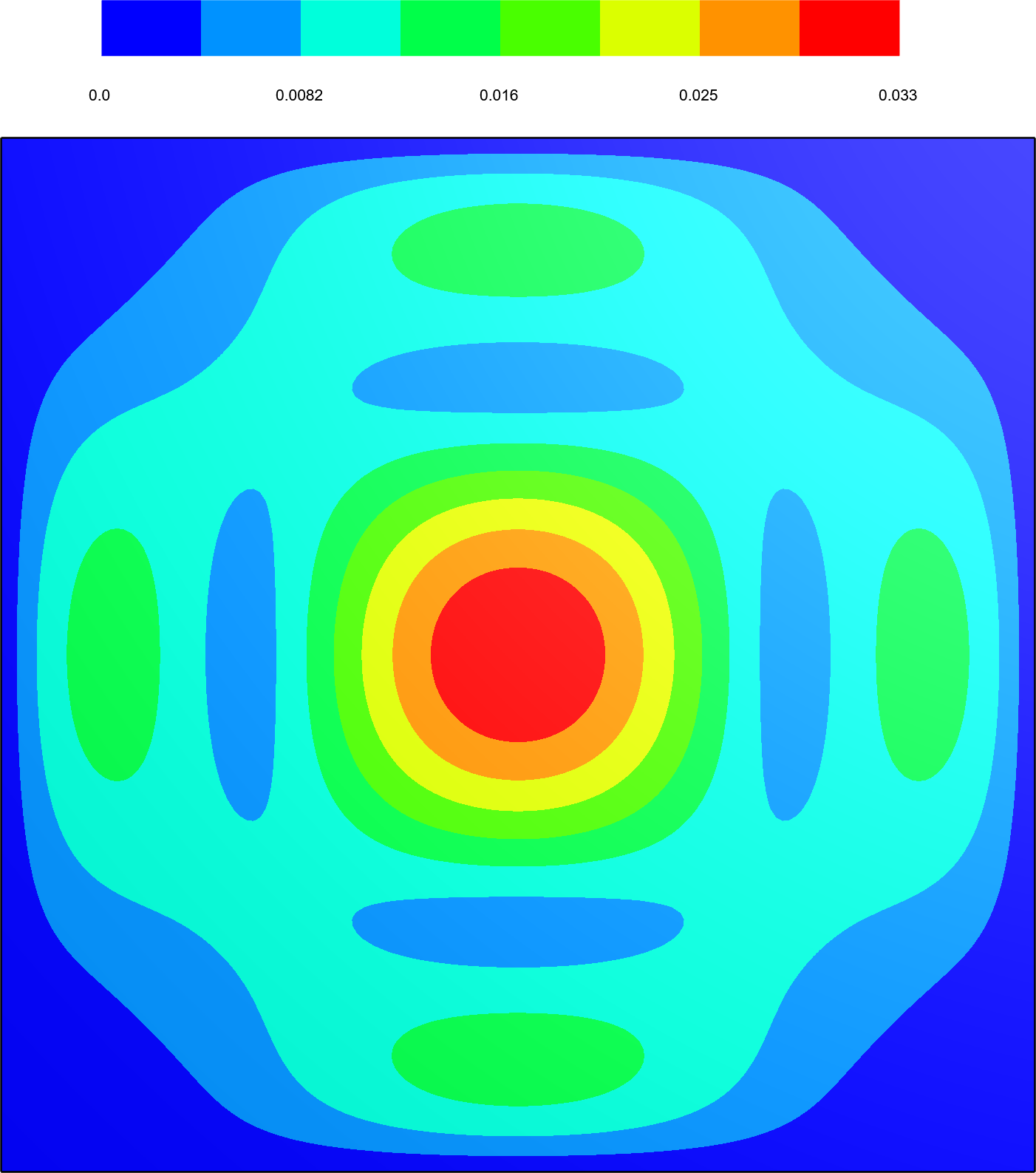}
    		\caption{}
    	\end{subfigure}
    	\begin{subfigure}{0.24\linewidth}
    		\centering
    		\includegraphics[width=1\linewidth]{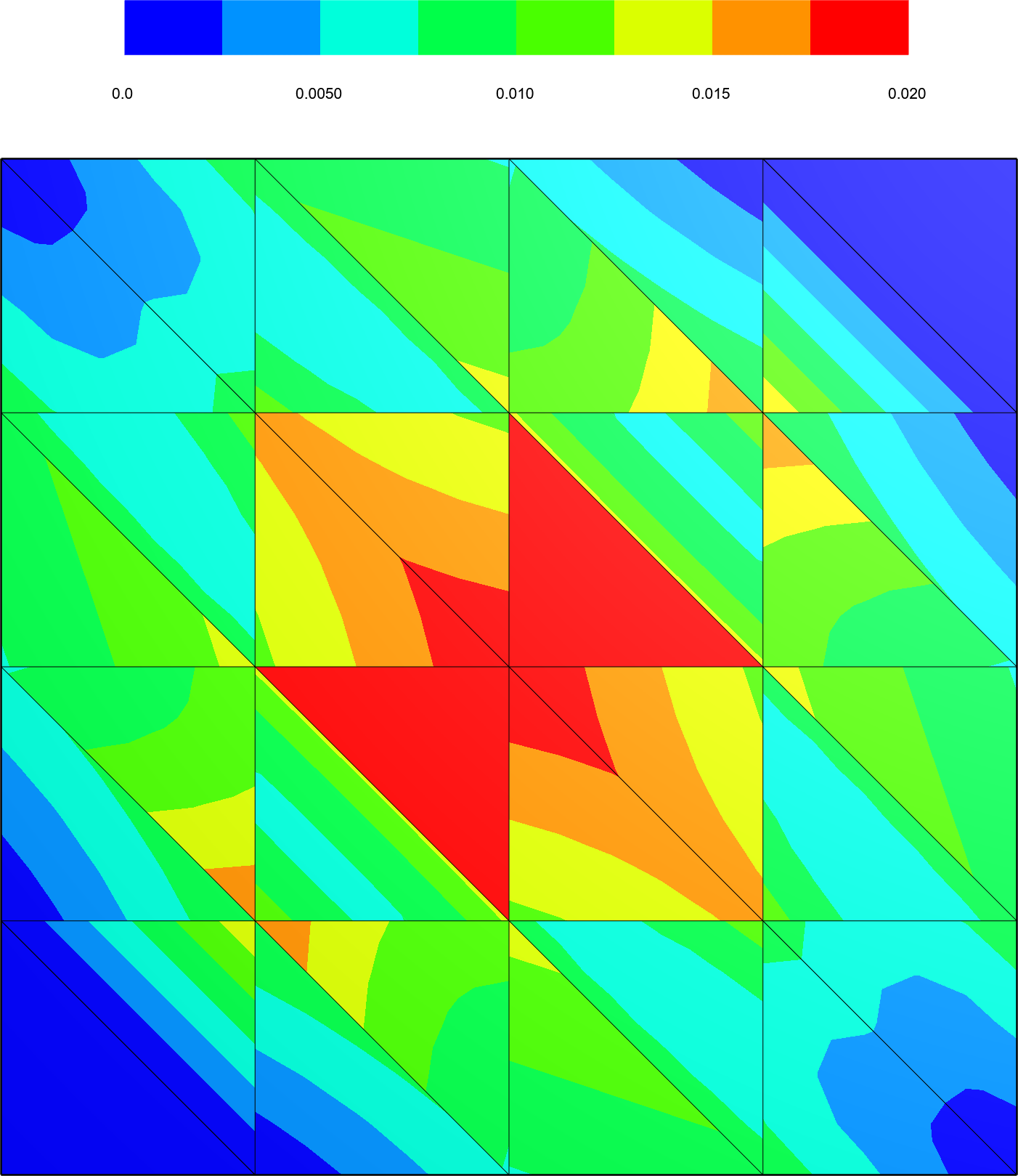}
            \caption{}
    	\end{subfigure}
        \begin{subfigure}{0.24\linewidth}
    		\centering
    		\includegraphics[width=1\linewidth]{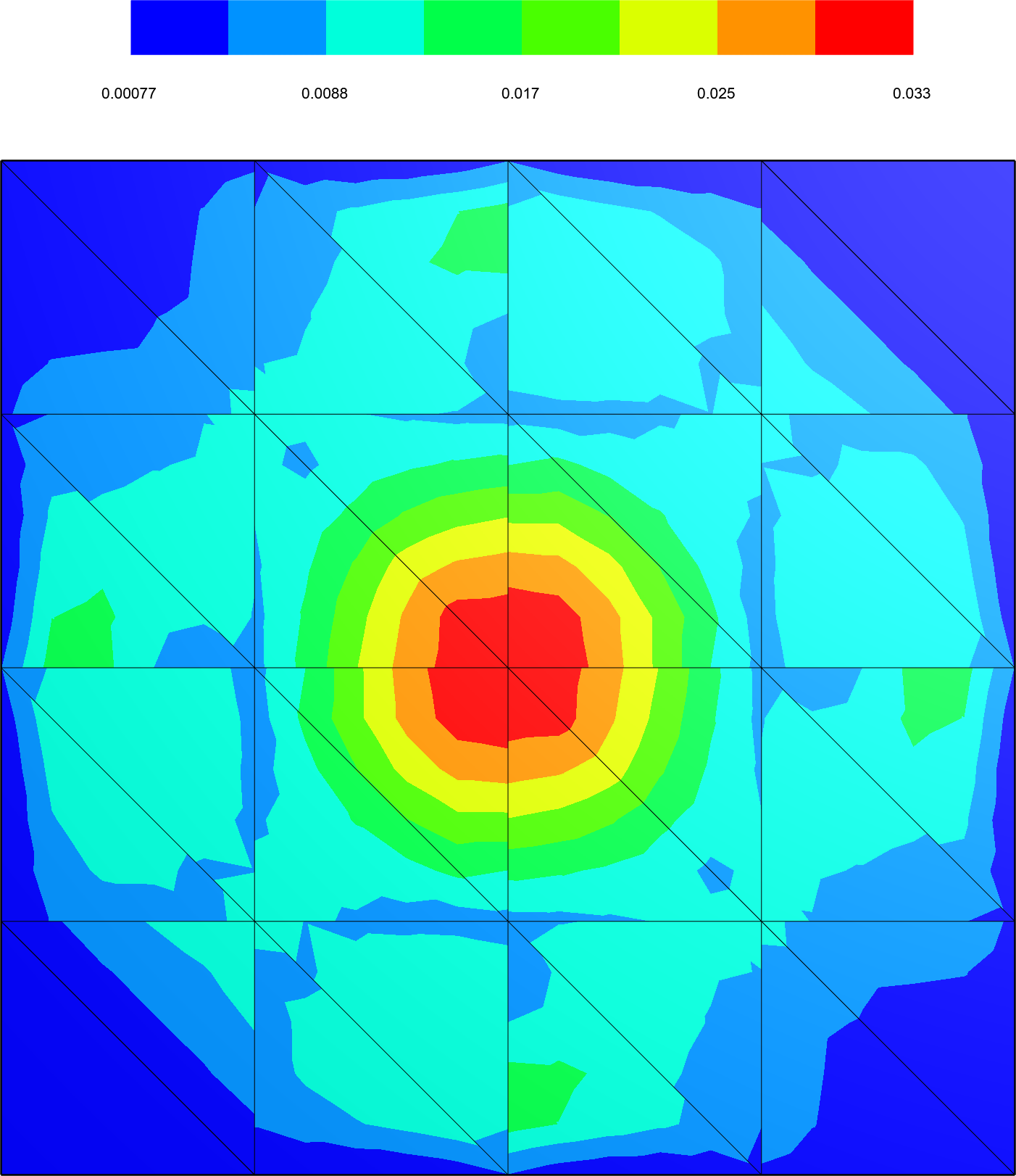}
            \caption{}
    	\end{subfigure}
        \begin{subfigure}{0.24\linewidth}
    		\centering
    		\includegraphics[width=1\linewidth]{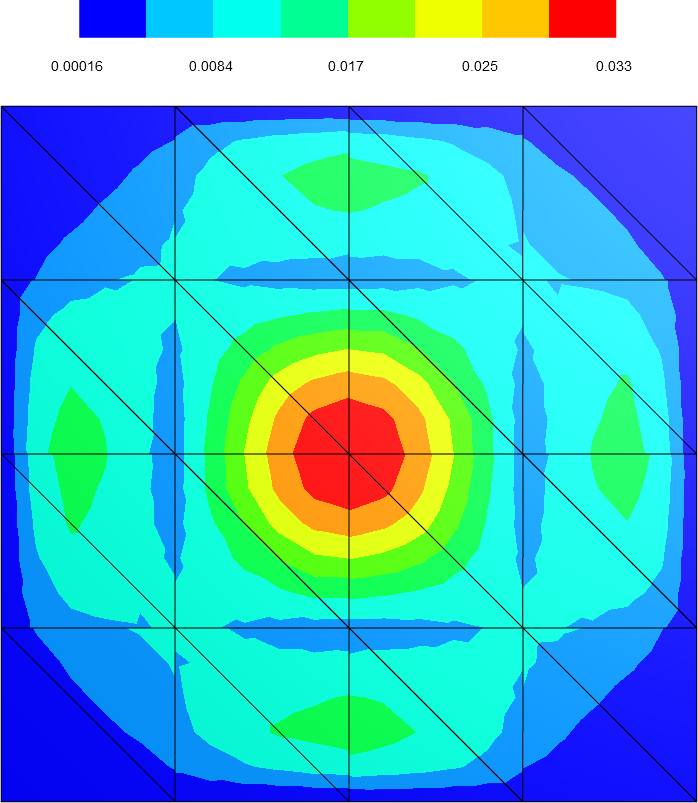}
            \caption{}
    	\end{subfigure}
    	\caption{Solution of the bending moments $\bm{M}$ analytically (a), via PRM, TFSRM and QFSRM, respectively. The approximations are depicted on a 32-element-mesh with $t = 10^{-5}$.}
    	\label{fig:ex1figs}
    \end{figure}

The result for the Frobenius norm of the bending moment field $\|\bm{M}\|$ with $2^5$-elements for the case $t = 10^{-5}$ is depicted in \cref{fig:ex1figs}.
Clearly, even when the primal formulation with cubic polynomials converges, it can produce rather poor approximations of the bending moments for a very small thickness $t \ll 1$. This is essential, since the bending moments are often used to determine the bearing capacity of plates in design processes. From \cref{fig:ex1figs} it becomes clear that the primal formulation underestimates the maximal bending moments by a factor of $\approx 1.5$, whereas both the TFSRM and QFSRM formulations find the correct maximum.     

\subsection{A curved circular plate}
\begin{figure}
    	\centering
    	\begin{subfigure}{0.3\linewidth}
    		\centering
    		\includegraphics[width=1\linewidth]{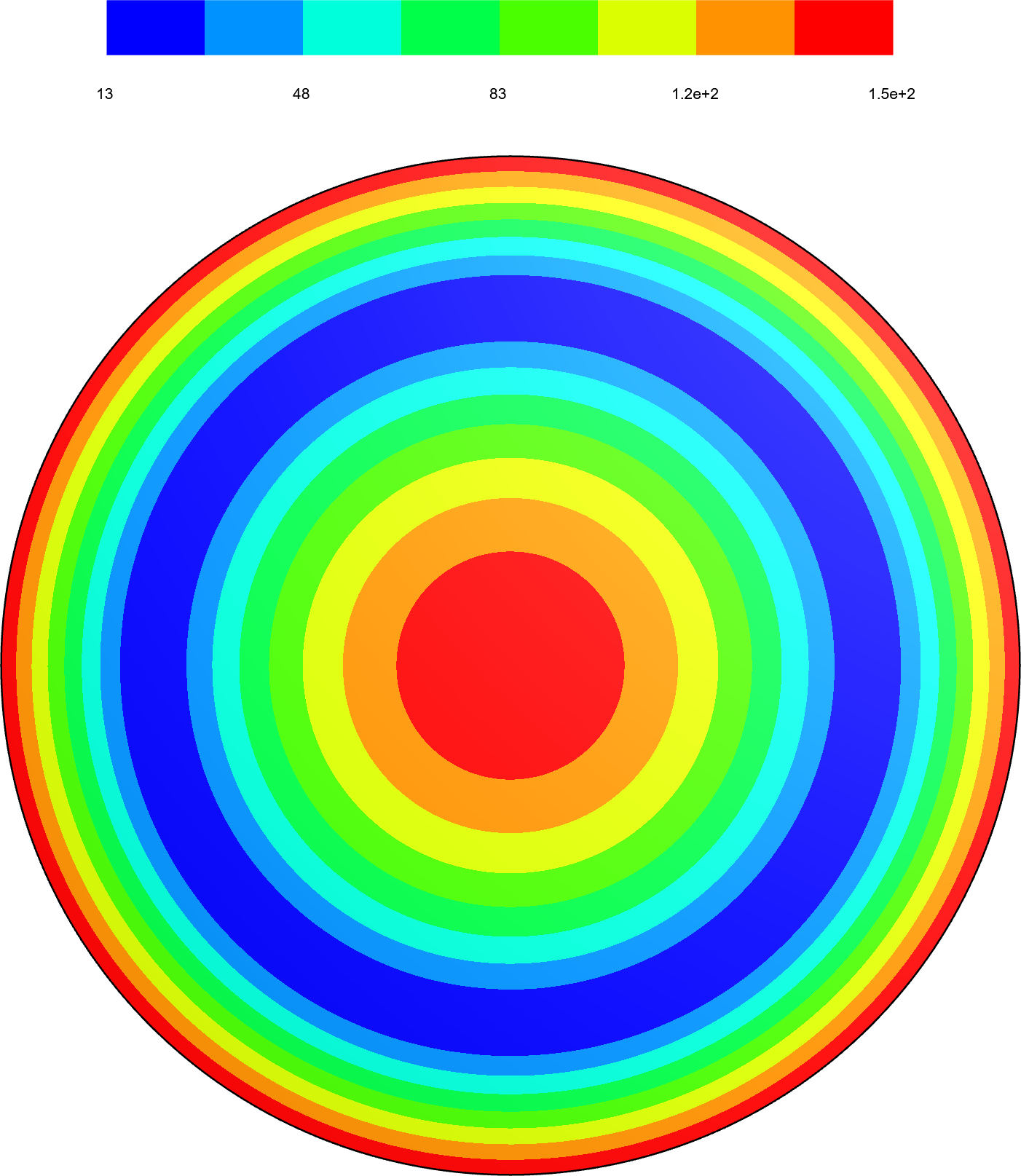}
    		\caption{}
    	\end{subfigure}
    	\begin{subfigure}{0.3\linewidth}
    		\centering
    		\includegraphics[width=1\linewidth]{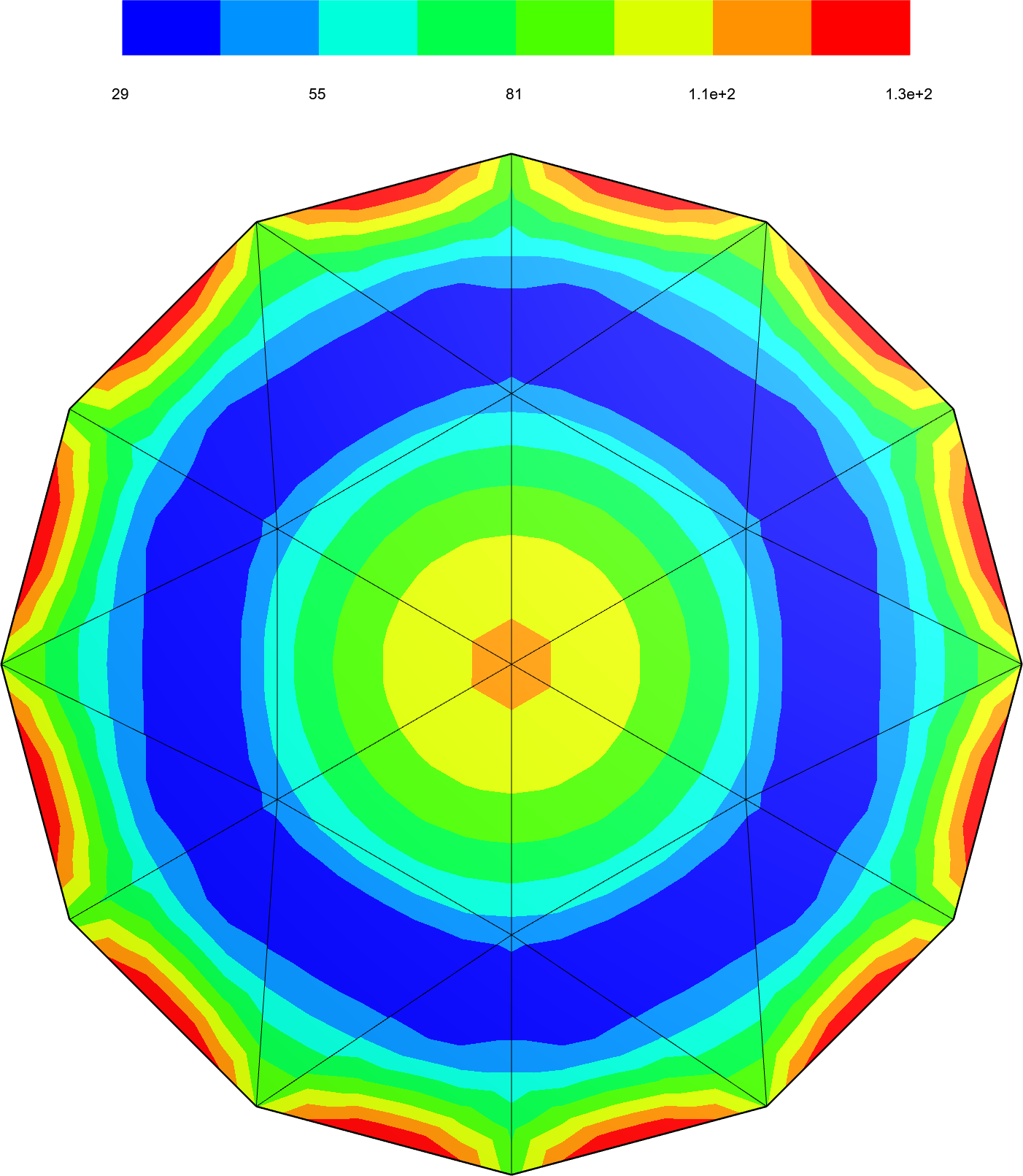}
            \caption{}
    	\end{subfigure}
        \begin{subfigure}{0.3\linewidth}
    		\centering
    		\includegraphics[width=1\linewidth]{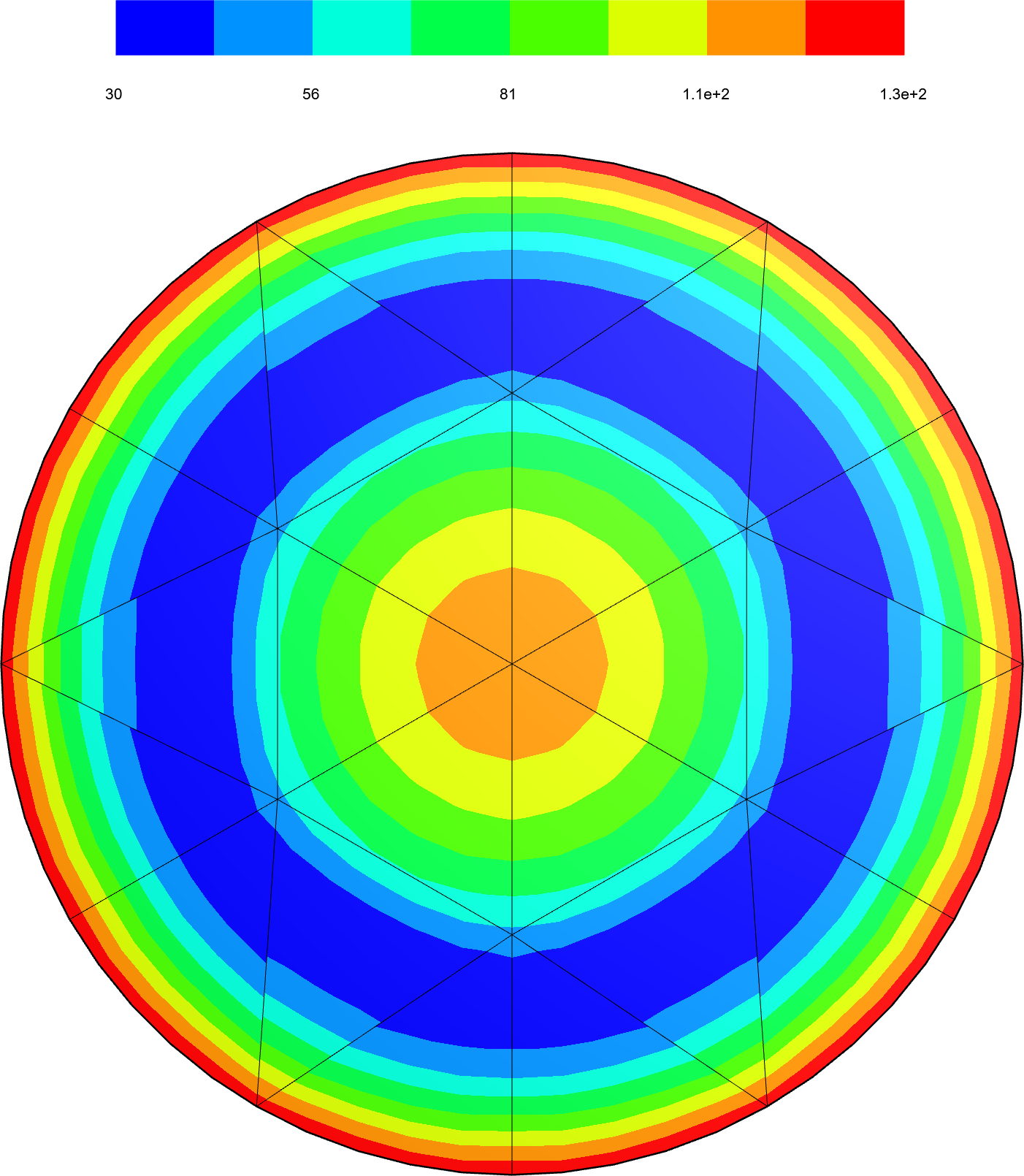}
            \caption{}
    	\end{subfigure}
    	\caption{Analytical solution of the bending moments $\bm{M}$ in the (a), TFSRM-formulation with a piece-wise linear (b) and piece-wise cubic (c) boundary using $24$ elements.}
    	\label{fig:ex2}
    \end{figure}

The transformation introduced in \cref{de:trans} allows to map the Hu-Zhang element from the reference triangle to a curved triangle on the physical domain. 
In order to demonstrate the effectiveness of said feature, in this example we consider the unit circle domain $\overline{\surf} = \{ (x,y) \in [-1,1]^2 \; | \; x^2 + y^2 \leq 1 \}$ with vanishing deflections and rotations on the boundary $\widetilde{w}|_{\partial \surf} = \widetilde{\bm{\phi}}|_{\partial \surf} = 0$ set analogously as in the previous example.
The constant forces are given by $f(x,y) = -1$, for which the analytical solution reads \cite{chinwuba_ike_mathematical_2018}
\begin{subequations}
    \begin{align}
    \widetilde{w} &= \dfrac{12(\nu^2-1)}{64 Et^3} (1-x^2 - y^2)^2 - \dfrac{1}{4 k_s \mu t} (1-x^2 -y^2) \, , \\
    \widetilde{\bm{\phi}} &= \dfrac{12(\nu^2-1)}{16 Et^3} (1-x^2 - y^2) \begin{bmatrix}
        x \\ y
    \end{bmatrix}   \, .
\end{align}
\end{subequations}
We set the material parameters to
\begin{align}
    &E = 240 \, , && \nu = 0.3 \, , && k_s = \dfrac{5}{6} \, , && t = 10^{-1} \, , 
\end{align}
and compute the formulations on a domain with a piece-wise linear and a piece-wise cubic boundary using $24$ cubic elements. 
The errors are listed in the following table and the deflection field for the TFSRM-formulation is depicted in \cref{fig:ex2}.
\begin{table}[H]
\centering
\begin{tabular}{l|cc|cc}
\multirow{2}{*}{} & \multicolumn{2}{c|}{} & \multicolumn{2}{c}{} \\ [-0.2cm]
\multirow{2}{*}{} & \multicolumn{2}{c|}{Piece-wise linear $\partial \surf^h$} & \multicolumn{2}{c}{Piece-wise cubic $\partial \surf^h$} \\ [0.2cm] \cline{2-5} 
\multirow{2}{*}{} & \multicolumn{2}{c|}{} & \multicolumn{2}{c}{}  \\ [-0.2cm]
                  & \multicolumn{1}{c|}{ 
                  $\| \widetilde{w} - w^h \|_{\Le} / \| \widetilde{w} \|_{\Le} $ }       & $\| \widetilde{\bm{\phi}} - \bm{\phi}^h \|_{\Le} / \| \widetilde{\bm{\phi}} \|_{\Le} $       & \multicolumn{1}{c|}{$\| \widetilde{w} - w^h \|_{\Le} / \| \widetilde{w} \|_{\Le} $}         & $\| \widetilde{\bm{\phi}} - \bm{\phi}^h \|_{\Le} / \| \widetilde{\bm{\phi}} \|_{\Le} $        \\ [0.2cm] \hline
                  \multirow{2}{*}{} & \multicolumn{2}{c|}{} & \multicolumn{2}{c}{} \\ [-0.2cm] 
TFSRM              & \multicolumn{1}{c|}{$11.6 \, \%$}        &   $11.1 \, \%$      & \multicolumn{1}{c|}{$0.2 \, \%$}          &   $ 1.1 \, \%$       \\[0.2cm] \hline
\multirow{2}{*}{} & \multicolumn{2}{c|}{} & \multicolumn{2}{c}{} \\ [-0.2cm] 
QFSRM             & \multicolumn{1}{c|}{$11.6 \, \%$}        &      $11.1 \, \%$   & \multicolumn{1}{c|}{$1.2 \, \%$}          &     $1.1 \, \%$     \\[0.2cm] \hline
\multirow{2}{*}{} & \multicolumn{2}{c|}{} & \multicolumn{2}{c}{} \\ [-0.2cm] 
TDNNS             & \multicolumn{1}{c|}{$11.3 \, \%$}        &      $10.8 \, \%$   & \multicolumn{1}{c|}{$0.2 \, \%$}          &     $1.7 \, \%$     \\[0.2cm] \hline
\multirow{2}{*}{} & \multicolumn{2}{c|}{} & \multicolumn{2}{c}{} \\ [-0.2cm] 
MITC              & \multicolumn{1}{c|}{$11.9 \, \%$}        &    $11.3 \, \%$     & \multicolumn{1}{c|}{$0.2 \, \%$}          &  $0.02 \, \%$        \\[0.2cm] \hline
\multirow{2}{*}{} & \multicolumn{2}{c|}{} & \multicolumn{2}{c}{} \\ [-0.2cm] 
PRM               & \multicolumn{1}{c|}{$12.1 \, \%$}        &      $11.6 \, \%$   & \multicolumn{1}{c|}{$0.3 \, \%$}          &   $0.8 \, \%$       \\[0.2cm] 
\end{tabular}
\end{table}
We note that for all formulations aside from QFSRM the difference in the exactness of the geometry yields a factor of approximately $\approx 100$ in the relative error of the deflection $w$. The polynomial order for the deflection $w$ in QFSRM is lower, such that a factor of $\approx 10$ is retrieved.
Excluding the MITC formulation, all the other formulation exhibit an improvement factor of $\approx 10$ in the relative error of the rotations $\bm{\phi}$. The MITC formulation is unique in its enrichment of the rotations $\bm{\phi} \in [\mathit{S}^p(\surf)]^2$, such that an improvement factor of $\approx 100$ is achieved. The vast improvement in the results due to the use of curved elements
illustrates the importance and necessity of corresponding mappings.

\subsection{Singularity on an L-shaped domain}
\begin{figure}
    	\centering
    	\begin{subfigure}{0.48\linewidth}
    		\centering
    		\input{figs/ex3_M}
            \vspace{-0.5cm}
    		\caption{}
      \label{fig:ex3a}
    	\end{subfigure}
         \begin{subfigure}{0.48\linewidth}
    		\centering
    		\includegraphics[width=0.7\linewidth]{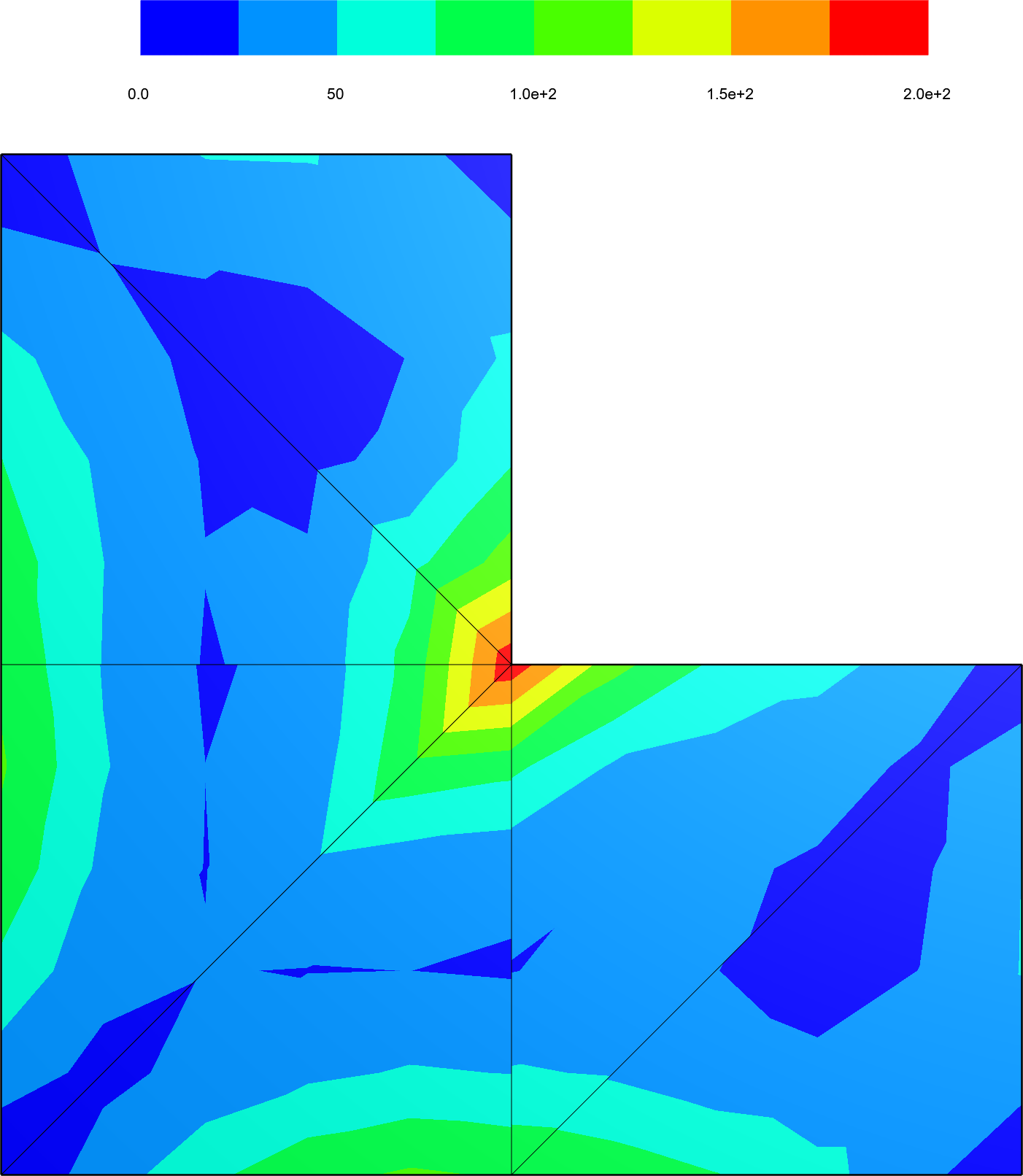}
    		\caption{}
    	\end{subfigure}
        \begin{subfigure}{0.48\linewidth}
    		\centering
    		\includegraphics[width=0.7\linewidth]{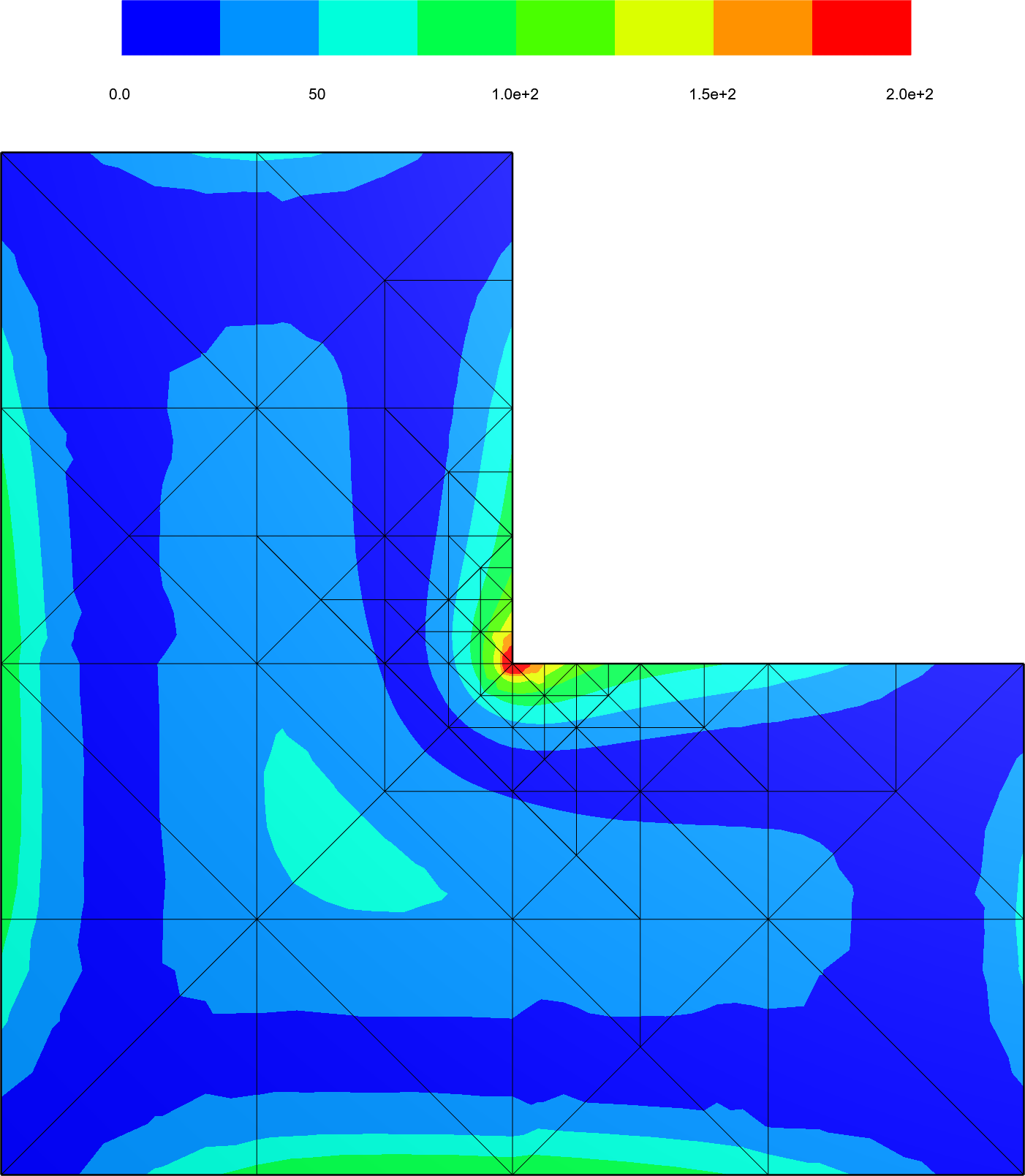}
    		\caption{}
    	\end{subfigure}
        \begin{subfigure}{0.48\linewidth}
    		\centering
    		\includegraphics[width=0.7\linewidth]{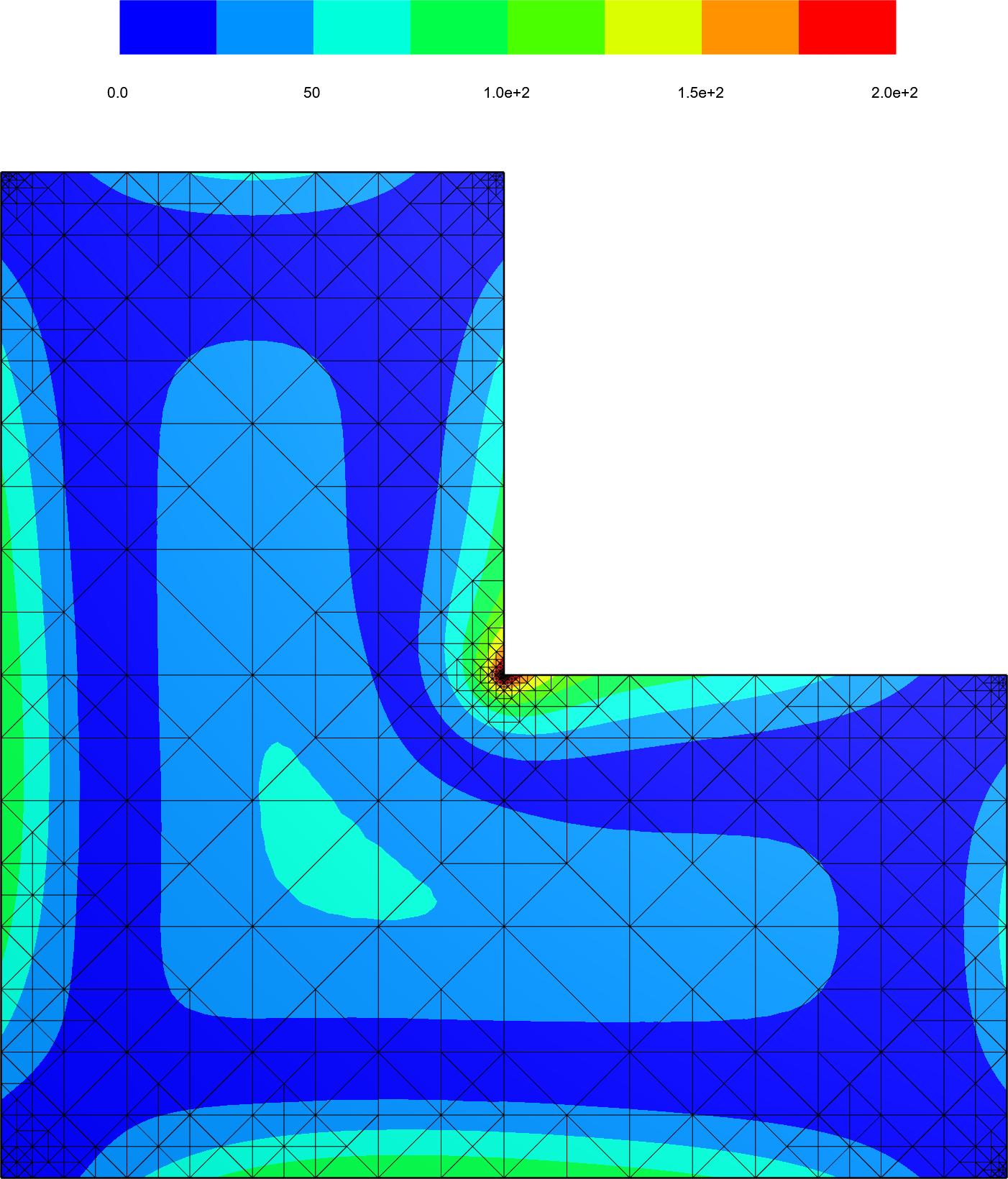}
    		\caption{}
    	\end{subfigure}
    	\caption{Relative error of the TFSRM formulation for $p \in \{3,5,7\}$ in the bending moments $\bm{M}$ (a). Depiction of the norm of the bending moments for $p = 5$ with $6$ (b), $114$ (c) and $1246$ (d) elements.}
    	\label{fig:ex3}
    \end{figure}

In this last example we demonstrate exponential convergence using h- and p-refinements in the presence of a singularity. 
We consider a fully clamped L-shaped domain $\overline{\surf} = [-1,1]^2 \setminus (0,1]^2$.
The boundary conditions of the TFSRM formulation are therefore complete Dirichlet for the deflection $\curv_D^w = \partial \surf$ and complete Neumann for the bending moments $\curv_N^M = \partial \surf$. At the re-entrant corner $(0,0)$, the bending moments produce a singularity, such that pure p-refinement no longer produces exponential convergence. This problem is alleviated by using adaptive h-refinement at the area of the singularity, therefore localising it. Since no analytical solution is available for this specific problem, we rely on a recovery-based error estimator \cite{GRATSCH2005235}. Under the assumption that the bending moments produce a smooth field we can define the recovery error estimator
\begin{align}
    \norm{\Pi_g^p \bm{M}^h - \bm{M}^h}_\Le / \norm{\Pi_g^p \bm{M}^h}_\Le \, , 
\end{align}
where $\Pi_g^g$ interpolates $\bm{M}^h$ into the $\C^0(\surf)$-continuous polynomial space $[\U^p(\surf)]^{2\times 2}$.
We set the material parameters to 
\begin{align}
    &E = 240 \, , && \nu = 0.3 \, , && t = 10^{-1} \, , 
\end{align}
and the constant force to $g(x,y) = -1000$. The convergence rates over h-refinement with $p \in \{3,5,7\}$, along with three solutions of the norm of the bending moment with $6$, $114$ and $1246$ elements, are depicted in \cref{fig:ex3}.
In \cref{fig:ex3a} we observe that the optimal convergence rates $\O(h^p)$ are retrieved under adaptive h-refinement. As expected, unlike for problems with smooth analytical solutions, we cannot expect convergence of the order $\O(h^{p+1})$.
In contrast, the quintic formulation achieves sub-optimal convergence under uniform h-refinement.
From the depictions of the bending moment in \cref{fig:ex3} it is apparent that the adaptive h-refinement scheme concentrates on the re-entrant corner, therefore reducing the distortion of the solution by the singularity.   

\section{Conclusions and outlook}
This paper proposes new mixed formulations for the Reissner-Mindlin plate based on the Hellinger-Reissner principle of symmetric stresses. The proposed formulations define the bending moments in the $\bm{M} \in \HsD{}$-space and rely on conforming Hu-Zhang finite elements for their discretisation $\HZ \subset \HsD{}$. Consequently, the rotations are intrinsic to the vector-valued discontinuous Lebesgue space $\bm{\phi} \in [\Le]^2$, such that the Kirchhoff-Love constraint can be satisfied for $t \to 0$ and locking in the sense of shear-locking is alleviated. The performance of the formulation was demonstrated in the first example, where for a non-too small thickness $t$, the TFSRM-formulation yields optimal convergence rates with higher accuracy in the bending moments $\bm{M}$. However, for $t \to 0$, the TFSRM formulation retrieves sub-optimal convergence rates in both the bending moments $\bm{M}$ and the shear stress $\vb{q}$. 
The QFSRM-formulation allows to alleviate this problem and exhibits optimal convergence rates across all variables also for $t \to 0$. We observe that in TFSRM, optimal convergence is maintained for $t \to 0$ in the deflection $w$ and the rotations $\bm{\phi}$ also without the additional field for the shear stress $\vb{q}$. Interestingly, our investigation also demonstrates the sensitivity of the MITC formulation in the thickness $t$ for approximations of the shear stress $\vb{q}$, such that for $t \to 0$ sub-optimal convergence in $\vb{q}$ is found (see also \cite{Bathe1990Dis}).  

Our second example accentuated the necessity of transformations from reference elements to curved elements on the physical domain by comparing the relative error induced via a piece-wise linear- and piece-wise cubic boundary. Clearly, the novel transformation proposed in this work allows to map the Hu-Zhang element to curved triangles.    

Finally, the construction of the Hu-Zhang element using Legendre polynomials allows to employ the element in the context of hp-FEM. The usefulness of said approach is demonstrated in the last example on an L-shaped domain with a singularity in the bending moments.

The mixed formulations proposed in this work were implemented using triangular elements. However, well-posedness of the discrete formulations is automatically inherited in the presence of commuting interpolants. As such, the formulation can be easily implemented for quadrilateral elements as well. For the Raviart-Thomas and Brezzi-Douglas-Marini elements we refer to \cite{Zaglmayr2006}. For symmetric $\HsD{}$-conforming elements we cite \cite{Hu2014Q,Hu2015Q}. 
This work did not consider the lowest order Arnold-Winther element, which could be used to further reduce the necessary amount of degrees of freedom.  

\section*{Acknowledgements}
Michael Neunteufel acknowledges support by the Austrian Science Fund (FWF) project F65.

\bibliographystyle{spmpsci}   

\footnotesize{
\bibliography{ref}  
} 


\end{document}

%% file: figs/seq.tex
\begin{tikzpicture}[line cap=round,line join=round,>=triangle 45,x=1.0cm,y=1.0cm]
				\clip(2,0) rectangle (10,1.5);
				\draw (3.3,0.9) node[anchor=north east] {$\H^2(\surf)$};
				\draw [-Triangle,line width=.5pt] (3.3,0.6) -- (4.8,0.6);
				\draw (4.,1.1) node[anchor=north] {$\airy$};
				\draw (4.8,0.9) node[anchor=north west] {$\HsD{,\surf}$};
				\draw [-Triangle,line width=.5pt] (7.1,0.6) -- (8.6,0.6);
				\draw (7.8,1.1) node[anchor=north] {$\Di$};
				\draw (8.6,0.9) node[anchor=north west] {$[\Le(\surf)]^2$};
			\end{tikzpicture}

%% file: figs/com_aw.tex
\begin{tikzpicture}[line cap=round,line join=round,>=triangle 45,x=1.0cm,y=1.0cm]
				\clip(2,-1.7) rectangle (10.5,1.5);
				\draw (3.3,0.9) node[anchor=north east] {$\H^2(\surf)$};
				\draw [-Triangle,line width=.5pt] (3.3,0.6) -- (4.8,0.6);
				\draw (4.,1.1) node[anchor=north] {$\airy$};
				\draw (4.8,0.9) node[anchor=north west] {$\HsD{,\surf}$};
				\draw [-Triangle,line width=.5pt] (7.1,0.6) -- (8.6,0.6);
				\draw (7.8,1.1) node[anchor=north] {$\Di$};
				\draw (8.6,0.9) node[anchor=north west] {$[\Le(\surf)]^2$};

                \draw [-Triangle,line width=.5pt] (2.7,0.3) -- (2.7,-1.2);
                \draw (2.7,-0.4) node[anchor=east] {$\Pi_a^5$};
                \draw [-Triangle,line width=.5pt] (6,0.3) -- (6,-1.2);
                \draw (6,-0.4) node[anchor=east] {$\Pi_s^2$};
                \draw [-Triangle,line width=.5pt] (9.3,0.3) -- (9.3,-1.2);
                \draw (9.3,-0.4) node[anchor=east] {$\Pi_o^1$};

                \draw (3.3,0.9-2) node[anchor=north east] {$\mathcal{A}^5(\surf)$};
				\draw [-Triangle,line width=.5pt] (3.3,0.6-2) -- (4.8,0.6-2);
				\draw (4.,1.1-2) node[anchor=north] {$\airy$};
				\draw (5.3,0.9-2) node[anchor=north west] {$\AW^2(\surf)$};
				\draw [-Triangle,line width=.5pt] (7.1,0.6-2) -- (8.6,0.6-2);
				\draw (7.8,1.1-2) node[anchor=north] {$\Di$};
				\draw (8.6,0.9-2) node[anchor=north west] {$[\mathit{D}^1(\surf)]^2$};
			\end{tikzpicture}

%% file: figs/com_hz.tex
\begin{tikzpicture}[line cap=round,line join=round,>=triangle 45,x=1.0cm,y=1.0cm]
				\clip(2,-1.7) rectangle (10.5,1.5);
				\draw (3.3,0.9) node[anchor=north east] {$\H^2(\surf)$};
				\draw [-Triangle,line width=.5pt] (3.3,0.6) -- (4.8,0.6);
				\draw (4.,1.1) node[anchor=north] {$\airy$};
				\draw (4.8,0.9) node[anchor=north west] {$\HsD{,\surf}$};
				\draw [-Triangle,line width=.5pt] (7.1,0.6) -- (8.6,0.6);
				\draw (7.8,1.1) node[anchor=north] {$\Di$};
				\draw (8.6,0.9) node[anchor=north west] {$[\Le(\surf)]^2$};

                \draw [-Triangle,line width=.5pt] (2.7,0.3) -- (2.7,-1.2);
                \draw (2.7,-0.4) node[anchor=east] {$\Pi_a^5$};
                \draw [-Triangle,line width=.5pt] (6,0.3) -- (6,-1.2);
                \draw (6,-0.4) node[anchor=east] {$\Pi_s^3$};
                \draw [-Triangle,line width=.5pt] (9.3,0.3) -- (9.3,-1.2);
                \draw (9.3,-0.4) node[anchor=east] {$\Pi_o^2$};

                \draw (3.3,0.9-2) node[anchor=north east] {$\mathcal{A}^5(\surf)$};
				\draw [-Triangle,line width=.5pt] (3.3,0.6-2) -- (4.8,0.6-2);
				\draw (4.,1.1-2) node[anchor=north] {$\airy$};
				\draw (5.3,0.9-2) node[anchor=north west] {$\HZ^3(\surf)$};
				\draw [-Triangle,line width=.5pt] (7.1,0.6-2) -- (8.6,0.6-2);
				\draw (7.8,1.1-2) node[anchor=north] {$\Di$};
				\draw (8.6,0.9-2) node[anchor=north west] {$[\mathit{D}^2(\surf)]^2$};
			\end{tikzpicture}

%% file: figs/com_rt.tex
\begin{tikzpicture}[line cap=round,line join=round,>=triangle 45,x=1.0cm,y=1.0cm]
				\clip(2,-1.7) rectangle (10,1.5);
				\draw (3.3,0.9) node[anchor=north east] {$\Hone(\surf)$};
				\draw [-Triangle,line width=.5pt] (3.3,0.6) -- (4.8,0.6);
				\draw (4.,1.1) node[anchor=north] {$\bm{R} \nabla$};
				\draw (5.1,0.9) node[anchor=north west] {$\Hd{,\surf}$};
				\draw [-Triangle,line width=.5pt] (7.1,0.6) -- (8.6,0.6);
				\draw (7.8,1.1) node[anchor=north] {$\di$};
				\draw (8.6,0.9) node[anchor=north west] {$\Le(\surf)$};

                \draw [-Triangle,line width=.5pt] (2.7,0.3) -- (2.7,-1.2);
                \draw (2.7,-0.4) node[anchor=east] {$\Pi^p_g$};
                \draw [-Triangle,line width=.5pt] (6,0.3) -- (6,-1.2);
                \draw (6,-0.4) node[anchor=east] {$\Pi^{p-1}_d$};
                \draw [-Triangle,line width=.5pt] (9.3,0.3) -- (9.3,-1.2);
                \draw (9.3,-0.4) node[anchor=east] {$\Pi^{p-1}_o$};

                \draw (3.3,0.9-2) node[anchor=north east] {$\U^p(\surf)$};
				\draw [-Triangle,line width=.5pt] (3.3,0.6-2) -- (4.8,0.6-2);
				\draw (4.,1.1-2) node[anchor=north] {$\bm{R} \nabla$};
				\draw (5.1,0.9-2) node[anchor=north west] {$\RT^{p-1}(\surf)$};
				\draw [-Triangle,line width=.5pt] (7.1,0.6-2) -- (8.6,0.6-2);
				\draw (7.8,1.1-2) node[anchor=north] {$\di$};
				\draw (8.5,0.9-2) node[anchor=north west] {$\mathit{D}^{p-1}(\surf)$};
			\end{tikzpicture}

%% file: figs/trimap.tex
\definecolor{asl}{rgb}{0.4980392156862745,0.,1.}
		\definecolor{asb}{rgb}{0.,0.4,0.6}
		\begin{tikzpicture}[line cap=round,line join=round,>=triangle 45,x=1.0cm,y=1.0cm]
			\clip(-1,-0.5) rectangle (12,3.5);
			\draw (0,0) node[circle,fill=asb,inner sep=1.5pt] {};
			\draw (0,2) node[circle,fill=asb,inner sep=1.5pt] {};
			\draw (2,0) node[circle,fill=asb,inner sep=1.5pt] {};
			\draw [color=asb,line width=.6pt] (0,0) -- (0,2) -- (2,0) -- (0,0);
			\fill[opacity=0.1, asb] (0,0) -- (0,2) -- (2,0) -- cycle;
			\draw (0,0) node[color=asb,anchor=north west] {$_{v_1}$};
			\draw (2,0) node[color=asb,anchor=south west] {$_{v_3}$};
			\draw (0,2) node[color=asb,anchor=north east] {$_{v_2}$};
			\draw (0.6,0.6) node[color=asb] {$\Gamma$};
			
			\draw [-to,color=asl,line width=1.pt] (0.75,1.25) -- (1.25,0.75);
			\draw [to-,color=asl,line width=1.pt] (1,1) -- (1.5,1.5);
			\draw (1.5,1.5) node[color=asl,anchor=south west] {$\bm{\nu}$};
			\draw (1.25,0.75) node[color=asl,anchor=west] {$\bm{\tau}$};
			
			\draw [-to,color=black,line width=1.pt] (0,0) -- (3,0);
			\draw [-to,color=black,line width=1.pt] (0,0) -- (0,3);
			\draw (3,0) node[color=black,anchor=west] {$\xi$};
			\draw (0,3) node[color=black,anchor=south] {$\eta$};
			
			\draw (9,1.5) node[circle,fill=asb,inner sep=1.5pt] {};
			\draw (10,3) node[circle,fill=asb,inner sep=1.5pt] {};
			\draw (11,1) node[circle,fill=asb,inner sep=1.5pt] {};
			\draw (9,1.5) node[color=asb,anchor=south east] {$_{\vb{x}_1}$};
			\draw (10,3) node[color=asb,anchor=south west] {$_{\vb{x}_3}$};
			\draw (11,1) node[color=asb,anchor=north east] {$_{\vb{x}_2}$};
			\draw [color=asb,line width=.6pt] (9,1.5) -- (10,3) -- (11,1) -- (9,1.5);
			\fill[opacity=0.1, asb] (9,1.5) -- (10,3) -- (11,1) -- cycle;
			\draw (10,1.85) node[color=asb] {$\surf_e$};
			\draw [-to,color=asl,line width=1.pt] (10.75,1.5) -- (10.25,2.5);
			\draw (10.3,2.5) node[color=asl,anchor=west] {$\vb{t}$};
			\draw [-to,color=asl,line width=1.pt] (10.5,2) -- (11.5,2.5);
			\draw (11,2.25) node[color=asl,anchor=north] {$\vb{n}$};
			
			\draw [-to,color=black,line width=1.pt] (8,0) -- (10,0);
			\draw [-to,color=black,line width=1.pt] (8,0) -- (8,2);
			\draw (10,0) node[color=black,anchor=west] {$x$};
			\draw (8,2) node[color=black,anchor=south] {$y$};
			
			\draw (5.3,1.5) node[color=black,anchor=south] {$\vb{x}:\Gamma \to \surf_e$};
			\draw [-Triangle,color=black,line width=1.pt] (4.5,1.5) -- (6.5,1.5);
		\end{tikzpicture}

%% file: figs/template.tex
\definecolor{asl}{rgb}{0.4980392156862745,0.,1.}
    		\definecolor{asb}{rgb}{0.,0.4,0.6}
    		\begin{tikzpicture}[line cap=round,line join=round,>=triangle 45,x=1.0cm,y=1.0cm]
    			\clip(-2,-1.5) rectangle (12.5,4.5);
    			\draw (-0.5,-0.5) node[circle,fill=asb,inner sep=1.5pt] {};
    			\draw (-0.5,4) node[circle,fill=asb,inner sep=1.5pt] {};
    			\draw (4,-0.5) node[circle,fill=asb,inner sep=1.5pt] {};
    			\draw [color=asb,line width=.6pt] (-0.5,0) -- (-0.5,3);
    			\draw [color=asb,line width=.6pt] (0,-0.5) -- (3,-0.5);
    			\draw [color=asb,line width=.6pt] (0.3,3.3) -- (3.3,0.3);
    			\draw [dotted,color=asb,line width=.6pt] (0,0) -- (0,3) -- (3,0) -- (0,0);
    			\fill[opacity=0.1, asb] (0,0) -- (0,3) -- (3,0) -- cycle;
    			\draw (-0.5,-0.5) node[color=asb,anchor=north east] {$_{v_1}$};
    			
                \draw [asl,domain=0:360,line width=1.pt] plot ({0.15*cos(\x-180)-0.5}, {0.15*sin(\x-180)-0.5});
    			
    			\draw (4,-0.5) node[color=asb,anchor=north west] {$_{v_3}$};
    			\draw (-0.5,4) node[color=asb,anchor=south east] {$_{v_2}$};

                \draw [asl,domain=0:360,line width=1.pt] plot ({0.15*cos(\x-180)-0.5}, {0.15*sin(\x-180)+4});
                \draw [asl,domain=0:360,line width=1.pt] plot ({0.15*cos(\x-180)+4}, {0.15*sin(\x-180)-0.5});
    			
    			\draw (-0.4,1.5) node[color=asb,anchor=east] {$_{e_{12}}$};
                \draw (-0.5,1.5) node[color=asl,anchor=west] {$\parallel$};
                \draw (-0.4,1.5) node[color=asl,anchor=north east] {\rotatebox{90}{$\bm{\parallel}$}};
                \draw (-0.4,1.5) node[color=asl,anchor=south east] {\rotatebox{-90}{$\bm{\perp}$}};
    			
    			\draw (1.55,-.5) node[color=asb,anchor=north] {$_{e_{13}}$};

                \draw (1.5,-.55) node[color=asl,anchor=south] {\rotatebox{90}{$\parallel$}};
                \draw (1.5,-.4) node[color=asl,anchor=north east] {$\bm{\parallel}$};
                \draw (1.5,-.4) node[color=asl,anchor=north west] {\rotatebox{180}{$\bm{\perp}$}};
    			
    			\draw (1.65,1.75) node[color=asb,anchor=south west] {$_{e_{23}}$};

                \draw (2,2) node[color=asl,anchor=north east] {\rotatebox{45}{$\parallel$}};
                \draw (1.8,1.4) node[color=asl,anchor=south west] {\rotatebox{135}{$\bm{\parallel}$}};
                \draw (1.4,1.8) node[color=asl,anchor=south west] {\rotatebox{135}{$\bm{\perp}$}};
    			
    			\draw (0.92,0.9)
    			node[color=asb] {$_{c_{123}}$};
                \draw [asl,domain=0:360,line width=1.pt, dashed] plot ({0.3*cos(\x-180)+0.92}, {0.3*sin(\x-180)+0.9});
    			
                \draw [asl,domain=0:360,line width=1.pt] plot ({0.15*cos(\x-180)+6.5}, {0.15*sin(\x-180)+2.25});
    			\draw (7,2.25)
    			node[color=asl,anchor=west] {Cartesian vertex template tensors};
    			\draw (6.5,1.5) node[color=asl] {$\bm{\parallel}$,$\bm{\perp}$};
    			\draw (7,1.5)
    			node[color=asl,anchor=west] {Edge template tensors};
                \draw (6.5,0.75) node[color=asl] {$\parallel$};
    			\draw (7,0.75)
    			node[color=asl,anchor=west] {Edge-cell template tensors};
                \draw [asl,domain=0:360,line width=1.pt, dashed] plot ({0.15*cos(\x-180)+6.5}, {0.15*sin(\x-180)});
    			\draw (7,0)
    			node[color=asl,anchor=west] {Cartesian cell template tensors};
    		\end{tikzpicture}

%% file: figs/curved.tex
\definecolor{asl}{rgb}{0.4980392156862745,0.,1.}
		\definecolor{asb}{rgb}{0.,0.4,0.6}
		\begin{tikzpicture}[line cap=round,line join=round,>=triangle 45,x=1.0cm,y=1.0cm]
			\clip(-1,-0.5) rectangle (12,3.5);
			\draw (0,0+0.5) node[circle,fill=asb,inner sep=1.5pt] {};
			\draw [color=asb,line width=.6pt] (0,0+0.5) -- (0,2+0.5) -- (2,0+0.5) -- (0,0+0.5);
			\fill[opacity=0.1, asb] (0,0+0.5) -- (0,2+0.5) -- (2,0+0.5) -- cycle;
			\draw (0,0+0.5) node[color=asb,anchor=north east] {$_{v}$};
            \draw [-to,color=asl,line width=1.pt] (0,0.5) -- (1,0.5);
            \draw (1,0.5) node[color=asl,anchor=north west] {$\vb{e}_1$};
            \draw [-to,color=asl,line width=1.pt] (0,0.5) -- (0,1.5);
            \draw (0,1.5) node[color=asl,anchor=south east] {$\vb{e}_2$};
            
			\draw (0.6,0.6+0.5) node[color=asb] {$\Gamma$};

           \draw [asb,domain=0:90,line width=.6pt] plot ({3*cos(\x+90)+10}, { 3*sin(\x+90)});
           \draw (10,3) node[circle,fill=black,inner sep=.6pt] {};
           \draw (7,0) node[circle,fill=black,inner sep=.6pt] {};
           \draw (7+0.8786796564403572,2.121320343559643) node[circle,fill=asb,inner sep=1.5pt] {};
           \draw (10,0) node[circle,fill=black,inner sep=.6pt] {};
           \fill [asb,domain=0:90, opacity=0.1] plot ({3*cos(\x+90)+10}, { 3*sin(\x+90)});
           \fill[opacity=0.1, asb] (7,0) -- (10,0) -- (10,3) -- cycle;

           \draw [color=asb,line width=.6pt] (7,0) -- (10,0) -- (10,3);
           \draw [color=asb,line width=.6pt] (7+0.8786796564403572,2.121320343559643) -- (10,0);

           \draw (7+0.87867965644035723,2.121320343559643)node[color=asb,anchor=south east] {$_{v}$};
           \draw [-to,color=asl,line width=1.pt] (7+0.87867965644035723,2.121320343559643) -- (7+0.8786796564403572+0.7071067811865475,2.121320343559643+0.7071067811865475);
           \draw (7+0.8786796564403572+0.7071067811865475,2.121320343559643+0.7071067811865475) node[color=asl,anchor=south west] {$\vb{d}_1$};
           \draw [-to,color=asl,line width=1.pt] (7+0.87867965644035723,2.121320343559643) -- (7+0.8786796564403572+0.7071067811865475,2.121320343559643-0.7071067811865475);
           \draw (7+0.8786796564403572+0.7071067811865475,2.121320343559643-0.7071067811865475) node[color=asl,anchor=north] {$\vb{d}_2$};

           \draw (8,0.75) node[color=asb] {$A_1$};
           \draw (9.25,2) node[color=asb] {$A_2$};
		
			\draw (5.3-1,1.5) node[color=black,anchor=south] {$\bm{Q}:\vb{e}_i \to \vb{d}_i$};
			\draw [-Triangle,color=black,line width=1.pt] (4.5-1,1.5) -- (6.5-1,1.5);
		\end{tikzpicture}

%% file: figs/ex1_e-1_w.tex
\begin{tikzpicture}
    			\definecolor{asl}{rgb}{0.4980392156862745,0.,1.}
    			\definecolor{asb}{rgb}{0.,0.4,0.6}
    			\begin{loglogaxis}[
    				/pgf/number format/1000 sep={},
    				axis lines = left,
    				xlabel={degrees of freedom},
    				ylabel={$\| \widetilde{w} - w^h \|_{\Le} / \| \widetilde{w} \|_{\Le} $ },
    				xmin=40, xmax=100000,
    				ymin=1e-6, ymax=50,
    				xtick={1e2,1e3,1e4,1e5},
    				ytick={1e-6,1e-4,1e-2, 1},
    				legend style={at={(0.95,1)},anchor= north east},
    				ymajorgrids=true,
    				grid style=dotted,
    				]
    				\addplot[color=asl, mark=triangle] coordinates {
    					( 90 , 1.1664246865854835 )
( 308 , 0.1577016942956067 )
( 1140 , 0.008366292981387493 )
( 4388 , 0.0004651026837020827 )
( 17220 , 2.729552659762359e-05 )
( 68228 , 1.6619553419387787e-06 )
    				};
    				\addlegendentry{TFSRM}
        \addplot[color=cyan, mark=diamond] coordinates {
    					( 113 , 1.3013465614821762 )
( 403 , 0.12320266074626204 )
( 1523 , 0.01893929814817134 )
( 5923 , 0.002607777639032599 )
( 23363 , 0.0003361292974670594 )
( 92803 , 4.235860606339746e-05 )
    				};
    				\addlegendentry{QFSRM}
    				
    				\addplot[color=violet, mark=o] coordinates {
    					( 70 , 1.4593871384856685 )
( 241 , 0.13358913496823077 )
( 889 , 0.007005683953781371 )
( 3409 , 0.00044831739716716526 )
( 13345 , 2.837151699037994e-05 )
( 52801 , 1.7772412837918691e-06 )
    				};
    				\addlegendentry{TDNNS}
    				
    				\addplot[color=asb, mark=square] coordinates {
                        ( 56 , 0.9327699074620524 )
( 179 , 0.10085127026314603 )
( 635 , 0.007578615209972425 )
( 2387 , 0.0004601929214491738 )
( 9251 , 2.8556647917600216e-05 )
( 36419 , 1.7801184282877486e-06 )
    				};
    				\addlegendentry{MITC}
    				
    				\addplot[color=blue, mark=pentagon] coordinates {
    					( 48 , 1.050470802373797 )
( 147 , 0.21930363849441245 )
( 507 , 0.012156844122413586 )
( 1875 , 0.0005250082488940946 )
( 7203 , 2.796440773278351e-05 )
( 28227 , 1.6708411596781342e-06 )
    				};
    				\addlegendentry{PRM}
    				
    				\addplot[dashed,color=black, mark=none]
    				coordinates {
    					(200, 0.0025118864315095794)
    					(2000, 2.5118864315095795e-05)
    				};    				

                   \addplot[dashed,color=black, mark=none]
    				coordinates {
    					(500, 1)
    					(2000, 0.125)
    				};    	
    			\end{loglogaxis}
    			\draw (2.2,1.1) 
    			node[anchor=south]{$\mathcal{O}(h^{4})$};
                \draw (3.05,4.1) 
    			node[anchor=south]{$\mathcal{O}(h^{3})$};
    		\end{tikzpicture}

%% file: figs/ex1_e-1_p.tex
\begin{tikzpicture}
    			\definecolor{asl}{rgb}{0.4980392156862745,0.,1.}
    			\definecolor{asb}{rgb}{0.,0.4,0.6}
    			\begin{loglogaxis}[
    				/pgf/number format/1000 sep={},
    				axis lines = left,
    				xlabel={degrees of freedom},
    				ylabel={$\| \widetilde{\bm{\phi}} - \bm{\phi}^h \|_{\Le} / \| \widetilde{\bm{\phi}} \|_{\Le} $ },
    				xmin=40, xmax=100000,
    				ymin=1e-6, ymax=50,
    				xtick={1e2,1e3,1e4,1e5},
    				ytick={1e-6,1e-4,1e-2, 1},
    				legend style={at={(0.95,1)},anchor= north east},
    				ymajorgrids=true,
    				grid style=dotted,
    				]
    				\addplot[color=asl, mark=triangle] coordinates {
    					( 90 , 1.0291875476077723 )
( 308 , 0.26638912327369313 )
( 1140 , 0.03913725659476879 )
( 4388 , 0.005188994381857399 )
( 17220 , 0.0006582122833061479 )
( 68228 , 8.289175116288201e-05 )
    				};
    				\addlegendentry{TFSRM}

        \addplot[color=cyan, mark=diamond] coordinates {
    					( 113 , 1.0688020675176826 )
( 403 , 0.2256637258267147 )
( 1523 , 0.036260967053064724 )
( 5923 , 0.005066737674784952 )
( 23363 , 0.0006566385441073786 )
( 92803 , 8.287782307160397e-05 )
    				};
    				\addlegendentry{QFSRM}
    				
    				\addplot[color=violet, mark=o] coordinates {
    					( 70 , 1.1393670279654358 )
( 241 , 0.4105476785753485 )
( 889 , 0.05300499248136813 )
( 3409 , 0.007236690538776733 )
( 13345 , 0.0009366210524952977 )
( 52801 , 0.00011826073259978438 )
    				};
    				\addlegendentry{TDNNS}
    				
    				\addplot[color=asb, mark=square] coordinates {
                        ( 56 , 0.968351606381999 )
( 179 , 0.13437541258655344 )
( 635 , 0.012917360471143908 )
( 2387 , 0.0008552245852441684 )
( 9251 , 5.2613704728259024e-05 )
( 36419 , 3.224365569190157e-06 )
    				};
    				\addlegendentry{MITC}
    				
    				\addplot[color=blue, mark=pentagon] coordinates {
    					( 48 , 1.0068021887896084 )
( 147 , 0.2915682463100399 )
( 507 , 0.020959291576392428 )
( 1875 , 0.0009999149671503945 )
( 7203 , 5.272120895393852e-05 )
( 28227 , 3.104806711869252e-06 )
    				};
    				\addlegendentry{PRM}

                    \addplot[dashed,color=black, mark=none]
    				coordinates {
    					(200, 1e-2)
    					(2000, 1e-04)
    				};   
    				
    				\addplot[dashed,color=black, mark=none]
    				coordinates {
    					(300, 1)
    					(3000, 0.03162277660168379)
    				};
    			\end{loglogaxis}
    			\draw (3,3.9) 
    			node[anchor=south]{$\mathcal{O}(h^{3})$};
    			\draw (2.2,1.5) 
    			node[anchor=south]{$\mathcal{O}(h^{4})$};
    		\end{tikzpicture}

%% file: figs/ex1_e-1_M.tex
\begin{tikzpicture}
    			\definecolor{asl}{rgb}{0.4980392156862745,0.,1.}
    			\definecolor{asb}{rgb}{0.,0.4,0.6}
    			\begin{loglogaxis}[
    				/pgf/number format/1000 sep={},
    				axis lines = left,
    				xlabel={degrees of freedom},
    				ylabel={$\| \widetilde{\bm{M}} - \bm{M}^h \|_{\Le} / \| \widetilde{\bm{M}} \|_{\Le} $ },
    				xmin=40, xmax=100000,
    				ymin=1e-6, ymax=50,
    				xtick={1e2,1e3,1e4,1e5},
    				ytick={1e-6,1e-4,1e-2, 1},
    				legend style={at={(0.95,1)},anchor= north east},
    				ymajorgrids=true,
    				grid style=dotted,
    				]
    				\addplot[color=asl, mark=triangle] coordinates {
    					( 90 , 1.1057528038595543 )
( 308 , 0.22024763180871387 )
( 1140 , 0.03664514809365704 )
( 4388 , 0.004227968755054931 )
( 17220 , 0.0003388265269079098 )
( 68228 , 2.272374186778709e-05 )
    				};
    				\addlegendentry{TFSRM}

        \addplot[color=cyan, mark=diamond] coordinates {
    					( 113 , 1.2946095248441862 )
( 403 , 0.17786461493221542 )
( 1523 , 0.020056245773143785 )
( 5923 , 0.0014436792638216002 )
( 23363 , 9.627056073909852e-05 )
( 92803 , 6.232696753308509e-06 )
    				};
    				\addlegendentry{QFSRM}
    				
    				\addplot[color=violet, mark=o] coordinates {
    					( 70 , 1.2263866541318211 )
( 241 , 0.40890830781106724 )
( 889 , 0.07735594447695991 )
( 3409 , 0.010627604855140334 )
( 13345 , 0.0013544610629046474 )
( 52801 , 0.00017017940231993072 )
    				};
    				\addlegendentry{TDNNS}
    				
    				\addplot[color=asb, mark=square] coordinates {
                        ( 56 , 0.8141650645767996 )
( 179 , 0.2177026534256839 )
( 635 , 0.05125354413999189 )
( 2387 , 0.007620759791066907 )
( 9251 , 0.0010157027464773288 )
( 36419 , 0.00012960731876161127 )
    				};
    				\addlegendentry{MITC}
    				
    				\addplot[color=blue, mark=pentagon] coordinates {
                             ( 48 , 1.0036382517491467 )
( 147 , 0.3937777390079299 )
( 507 , 0.07432872800363975 )
( 1875 , 0.008995172063471654 )
( 7203 , 0.0010845322750127821 )
( 28227 , 0.0001341713586402852 )
    				};
    				\addlegendentry{PRM}
    				
    				\addplot[dashed,color=black, mark=none]
    				coordinates {
    					(300, 1)
    					(3000, 0.03162277660168379)
    				};

                    \addplot[dashed,color=black, mark=none]
    				coordinates {
    					(1e+4, 2e-4)
    					(5e+4, 8e-6)
    				};
    			\end{loglogaxis}
    			\draw (3,3.9) 
    			node[anchor=south]{$\mathcal{O}(h^{3})$};
                \draw (5.2,0.6) 
    			node[anchor=south]{$\mathcal{O}(h^{4})$};
    		\end{tikzpicture}

%% file: figs/ex1_e-1_q.tex
\begin{tikzpicture}
    			\definecolor{asl}{rgb}{0.4980392156862745,0.,1.}
    			\definecolor{asb}{rgb}{0.,0.4,0.6}
    			\begin{loglogaxis}[
    				/pgf/number format/1000 sep={},
    				axis lines = left,
    				xlabel={degrees of freedom},
    				ylabel={$\| \widetilde{\vb{q}} - \vb{q}^h \|_{\Le} / \| \widetilde{\vb{q}} \|_{\Le} $ },
    				xmin=40, xmax=100000,
    				ymin=1e-6, ymax=50,
    				xtick={1e2,1e3,1e4,1e5},
    				ytick={1e-6,1e-4,1e-2, 1},
    				legend style={at={(0.95,1)},anchor= north east},
    				ymajorgrids=true,
    				grid style=dotted,
    				]
    				\addplot[color=asl, mark=triangle] coordinates {
    					( 90 , 0.9412019206546786 )
( 308 , 0.48549937042961977 )
( 1140 , 0.11647253335967538 )
( 4388 , 0.020587416195345144 )
( 17220 , 0.002977572204900863 )
( 68228 , 0.00038759815585101114 )
    				};
    				\addlegendentry{TFSRM}

        \addplot[color=cyan, mark=diamond] coordinates {
    					( 113 , 0.7642864259069242 )
( 403 , 0.21742308559850543 )
( 1523 , 0.0376687447534508 )
( 5923 , 0.00474667201389821 )
( 23363 , 0.0006000371681243421 )
( 92803 , 7.564423217102112e-05 )
    				};
    				\addlegendentry{QFSRM}
    				
    				\addplot[color=violet, mark=o] coordinates {
    					( 70 , 0.94384668806694 )
( 241 , 0.4452190137356832 )
( 889 , 0.07107126689788644 )
( 3409 , 0.009264497946508816 )
( 13345 , 0.0011545287649614884 )
( 52801 , 0.00014347460344119853 )
    				};
    				\addlegendentry{TDNNS}
    				
    				\addplot[color=asb, mark=square] coordinates {
                        ( 56 , 0.9354122999228609 )
( 179 , 0.4793284849167587 )
( 635 , 0.08054270065344643 )
( 2387 , 0.009926935167157029 )
( 9251 , 0.0011810527415542964 )
( 36419 , 0.00014436886259185507 )
    				};
    				\addlegendentry{MITC}
    				
    				\addplot[color=blue, mark=pentagon] coordinates {
                             ( 48 , 0.9117604039828263 )
( 147 , 0.7869717612349784 )
( 507 , 0.21287464630938374 )
( 1875 , 0.0347225904303337 )
( 7203 , 0.004590116513169589 )
( 28227 , 0.0005792668216848101 )
    				};
    				\addlegendentry{PRM}
    				
    				\addplot[dashed,color=black, mark=none]
    				coordinates {
    					(300, 3)
    					(3000, 0.09486832980505136)
    				};
    			\end{loglogaxis}
    			\draw (3,4.25) 
    			node[anchor=south]{$\mathcal{O}(h^{3})$};;
    		\end{tikzpicture}

%% file: figs/ex1_e-5_w.tex
\begin{tikzpicture}
    			\definecolor{asl}{rgb}{0.4980392156862745,0.,1.}
    			\definecolor{asb}{rgb}{0.,0.4,0.6}
    			\begin{loglogaxis}[
    				/pgf/number format/1000 sep={},
    				axis lines = left,
    				xlabel={degrees of freedom},
    				ylabel={$\| \widetilde{w} - w^h \|_{\Le} / \| \widetilde{w} \|_{\Le} $ },
    				xmin=40, xmax=100000,
    				ymin=1e-6, ymax=50,
    				xtick={1e2,1e3,1e4,1e5},
    				ytick={1e-6,1e-4,1e-2, 1},
    				legend style={at={(0.95,1)},anchor= north east},
    				ymajorgrids=true,
    				grid style=dotted,
    				]
    				\addplot[color=asl, mark=triangle] coordinates {
    					( 90 , 1.1295056010555566 )
( 308 , 0.16786876794688832 )
( 1140 , 0.010382379001524692 )
( 4388 , 0.000686238260858833 )
( 17220 , 4.3360073379663814e-05 )
( 68228 , 6.0204635896238975e-06 )
    				};
    				\addlegendentry{TFSRM}

                    \addplot[color=cyan, mark=diamond] coordinates {
    					( 113 , 1.2476373130846867 )
( 403 , 0.1191873496701265 )
( 1523 , 0.01566003784005293 )
( 5923 , 0.0022315551965655994 )
( 23363 , 0.0002914214466882587 )
( 92803 , 3.684872793794284e-05 )
    				};
    				\addlegendentry{QFSRM}
    				
    				\addplot[color=violet, mark=o] coordinates {
        ( 70 , 1.4380113470071056 )
( 241 , 0.14111270558986114 )
( 889 , 0.00685194957739399 )
( 3409 , 0.000433534217377639 )
( 13345 , 2.7718207297063424e-05 )
( 52801 , 1.546830723045226e-05 )
    				};
    				\addlegendentry{TDNNS}
    				
    				\addplot[color=asb, mark=square] coordinates {
                        ( 56 , 0.9509577375510833 )
( 179 , 0.08947400626540072 )
( 635 , 0.008990173409410825 )
( 2387 , 0.00047660482169916514 )
( 9251 , 2.8327878131682595e-05 )
( 36419 , 5.418101102197779e-06 )
    				};
    				\addlegendentry{MITC}
    				
    				\addplot[color=blue, mark=pentagon] coordinates {
    					( 48 , 1.0000000005691119 )
( 147 , 0.9999997262564514 )
( 507 , 0.17607691304604606 )
( 1875 , 0.039079805952224586 )
( 7203 , 0.009344006583981108 )
( 28227 , 0.0022982135857840934 )
    				};
    				\addlegendentry{PRM}
    				
    				\addplot[dashed,color=black, mark=none]
    				coordinates {
    					(200, 0.0025118864315095794)
    					(2000, 2.5118864315095795e-05)
    				};    		
                    \addplot[dashed,color=black, mark=none]
    				coordinates {
    					(500, 1)
    					(2000, 0.25)
    				};    		
    			\end{loglogaxis}
    			\draw (2.2,1.1) 
    			node[anchor=south]{$\mathcal{O}(h^{4})$};
                \draw (3,4.2) 
    			node[anchor=south]{$\mathcal{O}(h^{2})$};
    		\end{tikzpicture}

%% file: figs/ex1_e-5_p.tex
\begin{tikzpicture}
    			\definecolor{asl}{rgb}{0.4980392156862745,0.,1.}
    			\definecolor{asb}{rgb}{0.,0.4,0.6}
    			\begin{loglogaxis}[
    				/pgf/number format/1000 sep={},
    				axis lines = left,
    				xlabel={degrees of freedom},
    				ylabel={$\| \widetilde{\bm{\phi}} - \bm{\phi}^h \|_{\Le} / \| \widetilde{\bm{\phi}} \|_{\Le} $ },
    				xmin=40, xmax=100000,
    				ymin=1e-6, ymax=50,
    				xtick={1e2,1e3,1e4,1e5},
    				ytick={1e-6,1e-4,1e-2, 1},
    				legend style={at={(0.95,1)},anchor= north east},
    				ymajorgrids=true,
    				grid style=dotted,
    				]
    				\addplot[color=asl, mark=triangle] coordinates {
    					( 90 , 1.030152897310038 )
( 308 , 0.28299627001637107 )
( 1140 , 0.04804548400121678 )
( 4388 , 0.0068845352751627205 )
( 17220 , 0.0008879794546554563 )
( 68228 , 0.00011174628913154569 )
    				};
    				\addlegendentry{TFSRM}

                    \addplot[color=cyan, mark=diamond] coordinates {
    					( 113 , 1.0715393126532395 )
( 403 , 0.22350449603875674 )
( 1523 , 0.03586278563308439 )
( 5923 , 0.00506007795361093 )
( 23363 , 0.0006565724512657511 )
( 92803 , 8.287724114492166e-05 )
    				};
    				\addlegendentry{QFSRM}
    				
    				\addplot[color=violet, mark=o] coordinates {
        ( 70 , 1.1404588777159652 )
( 241 , 0.4121015626006295 )
( 889 , 0.05279369141362558 )
( 3409 , 0.007238282898663128 )
( 13345 , 0.0009367401429181414 )
( 52801 , 0.00011903951845284012 )
    				};
    				\addlegendentry{TDNNS}
    				
    				\addplot[color=asb, mark=square] coordinates {
                        ( 56 , 0.9653164467727271 )
( 179 , 0.1407213847129342 )
( 635 , 0.018539629311597538 )
( 2387 , 0.001174827854666955 )
( 9251 , 6.118396721703338e-05 )
( 36419 , 5.770236345621004e-06 )
    				};
    				\addlegendentry{MITC}
    				
    				\addplot[color=blue, mark=pentagon] coordinates {
    					( 48 , 1.0000000000881266 )
( 147 , 0.9999997474195524 )
( 507 , 0.2555439016839703 )
( 1875 , 0.0673767068447591 )
( 7203 , 0.016211173267663685 )
( 28227 , 0.003978607115545357 )
    				};
    				\addlegendentry{PRM}
    				
    				\addplot[dashed,color=black, mark=none]
    				coordinates {
    					(200, 1e-2)
    					(2000, 1e-04)
    				};   
    				
    				\addplot[dashed,color=black, mark=none]
    				coordinates {
    					(300, 1)
    					(3000, 0.1)
    				};
    			\end{loglogaxis}
    			\draw (3,4.05) 
    			node[anchor=south]{$\mathcal{O}(h^{2})$};
    			\draw (2.2,1.55) 
    			node[anchor=south]{$\mathcal{O}(h^{4})$};
    		\end{tikzpicture}

%% file: figs/ex1_e-5_M.tex
\begin{tikzpicture}
    			\definecolor{asl}{rgb}{0.4980392156862745,0.,1.}
    			\definecolor{asb}{rgb}{0.,0.4,0.6}
    			\begin{loglogaxis}[
    				/pgf/number format/1000 sep={},
    				axis lines = left,
    				xlabel={degrees of freedom},
    				ylabel={$\| \widetilde{\bm{M}} - \bm{M}^h \|_{\Le} / \| \widetilde{\bm{M}} \|_{\Le} $ },
    				xmin=40, xmax=100000,
    				ymin=1e-6, ymax=50,
    				xtick={1e2,1e3,1e4,1e5},
    				ytick={1e-6,1e-4,1e-2, 1},
    				legend style={at={(0.05,0.05)},anchor= south west},
    				ymajorgrids=true,
    				grid style=dotted,
    				]
    				\addplot[color=asl, mark=triangle] coordinates {
    					( 90 , 1.114543427364844 )
( 308 , 0.26484842577813617 )
( 1140 , 0.09227000231348946 )
( 4388 , 0.026592242031508367 )
( 17220 , 0.006776425066033525 )
( 68228 , 0.0016752305203324264 )
    				};
    				\addlegendentry{TFSRM}

                    \addplot[color=cyan, mark=diamond] coordinates {
    					( 113 , 1.2955842322158804 )
( 403 , 0.17249704797183613 )
( 1523 , 0.018358566848859947 )
( 5923 , 0.0012661602467930971 )
( 23363 , 8.3031216038524e-05 )
( 92803 , 5.348538685631e-06 )
    				};
    				\addlegendentry{QFSRM}
    				
    				\addplot[color=violet, mark=o] coordinates {
        ( 70 , 1.2275568328929685 )
( 241 , 0.40909666570274344 )
( 889 , 0.07719888469086597 )
( 3409 , 0.010611618388551911 )
( 13345 , 0.001353116537138617 )
( 52801 , 0.00017056625270538305 )
    				};
    				\addlegendentry{TDNNS}
    				
    				\addplot[color=asb, mark=square] coordinates {
                        ( 56 , 0.8182533207115036 )
( 179 , 0.2514875352902801 )
( 635 , 0.07360901431512705 )
( 2387 , 0.011034970595762198 )
( 9251 , 0.0013555433411832948 )
( 36419 , 0.00016544020435374282 )
    				};
    				\addlegendentry{MITC}
    				
    				\addplot[color=blue, mark=pentagon] coordinates {
        ( 48 , 1.0000000000602536 )
( 147 , 0.9999997575832174 )
( 507 , 0.4492909443706297 )
( 1875 , 0.23662481423123669 )
( 7203 , 0.11591707726736444 )
( 28227 , 0.05709292280367681 )
    				};
    				\addlegendentry{PRM}
    				
    				\addplot[dashed,color=black, mark=none]
    				coordinates {
    					(300, 1)
    					(3000, 0.31622776601683794)
    				};

                    \addplot[dashed,color=black, mark=none]
    				coordinates {
    					(1e+4, 2e-4)
    					(5e+4, 8e-6)
    				};

                    \addplot[dashed,color=black, mark=none]
    				coordinates {
    					(3e+3, 7e-2)
    					(4e+4, 0.00525)
    				};

                    \addplot[dashed,color=black, mark=none]
    				coordinates {
    					(1e+2, 2e-1)
    					(8e+2, 0.008838834764831844)
    				};
    			\end{loglogaxis}
    			\draw (3,4.15) 
    			node[anchor=south]{$\mathcal{O}(h^{1})$};
                \draw (5.15,3.05) 
    			node[anchor=south]{$\mathcal{O}(h^{2})$};
                \draw (1.35,2.85) 
    			node[anchor=south]{$\mathcal{O}(h^{3})$};
                \draw (5.2,0.6) 
    			node[anchor=south]{$\mathcal{O}(h^{4})$};
    		\end{tikzpicture}

%% file: figs/ex1_e-5_q.tex
\begin{tikzpicture}
    			\definecolor{asl}{rgb}{0.4980392156862745,0.,1.}
    			\definecolor{asb}{rgb}{0.,0.4,0.6}
    			\begin{loglogaxis}[
    				/pgf/number format/1000 sep={},
    				axis lines = left,
    				xlabel={degrees of freedom},
    				ylabel={$\| \widetilde{\vb{q}} - \vb{q}^h \|_{\Le} / \| \widetilde{\vb{q}} \|_{\Le} $ },
    				xmin=40, xmax=100000,
    				ymin=1e-6, ymax=50,
    				xtick={1e2,1e3,1e4,1e5},
    				ytick={1e-6,1e-4,1e-2, 1},
    				legend style={at={(0.05,0.05)},anchor= south west},
    				ymajorgrids=true,
    				grid style=dotted,
    				]
    				\addplot[color=asl, mark=triangle] coordinates {
    					( 90 , 0.9471222039056438 )
( 308 , 0.5632422761094182 )
( 1140 , 0.32411240701607574 )
( 4388 , 0.1854499999510951 )
( 17220 , 0.0953984741117806 )
( 68228 , 0.04680509898852171 )
    				};
    				\addlegendentry{TFSRM}

                    \addplot[color=cyan, mark=diamond] coordinates {
    					( 113 , 0.7653716016805013 )
( 403 , 0.245659698226654 )
( 1523 , 0.047219287859403 )
( 5923 , 0.006216289371960734 )
( 23363 , 0.0007776999882370374 )
( 92803 , 9.75035175573159e-05 )
    				};
    				\addlegendentry{QFSRM}
    				
    				\addplot[color=violet, mark=o] coordinates {
        ( 70 , 0.9447822364572867 )
( 241 , 0.4519609787386566 )
( 889 , 0.07316849359737641 )
( 3409 , 0.009951439288359457 )
( 13345 , 0.001374464185152153 )
( 52801 , 0.00021176702737748176 )
    				};
    				\addlegendentry{TDNNS}
    				
    				\addplot[color=asb, mark=square] coordinates {
                        ( 56 , 0.9386044606010054 )
( 179 , 0.6880950396431051 )
( 635 , 0.3286211929222826 )
( 2387 , 0.1138808340801877 )
( 9251 , 0.028509466996062816 )
( 36419 , 0.0069191032597119145 )
    				};
    				\addlegendentry{MITC}
    				
    				\addplot[color=blue, mark=pentagon] coordinates {
    					( 48 , 0.9058993533977119 )
( 147 , 9.065552288036736 )
( 507 , 9.580583660718263 )
( 1875 , 22.589542164450034 )
( 7203 , 46.051390228017006 )
( 28227 , 91.85373572217678 )
    				};
    				\addlegendentry{PRM}

        \addplot[dashed,color=black, mark=none]
    				coordinates {
    					(1e+4, 4e-4)
    					(5e+4, 3.5777087639996646e-05)
    				};

                    \addplot[dashed,color=black, mark=none]
    				coordinates {
    					(2e+3, 1)
    					(2e+4, 0.31622776601683794)
    				};

    			\end{loglogaxis}
    			\draw (5.25,0.95) 
    			node[anchor=south]{$\mathcal{O}(h^{3})$};
                \draw (4.6,4.2) 
    			node[anchor=south]{$\mathcal{O}(h^{1})$};
    		\end{tikzpicture}

%% file: figs/ex3_M.tex
\begin{tikzpicture}
    			\definecolor{asl}{rgb}{0.4980392156862745,0.,1.}
    			\definecolor{asb}{rgb}{0.,0.4,0.6}
    			\begin{loglogaxis}[
    				/pgf/number format/1000 sep={},
    				axis lines = left,
    				xlabel={degrees of freedom},
    				ylabel={$\norm{\Pi_g^p \bm{M}^h - \bm{M}^h}_\Le / \norm{\Pi_g^p \bm{M}^h}_\Le$},
    				xmin=700, xmax=150000,
    				ymin=1e-6, ymax=0.1,
    				xtick={1e3, 1e4,1e5},
    				ytick={1e-4,1e-2, 1},
    				legend style={at={(0.05,0.05)},anchor= south west},
    				ymajorgrids=true,
    				grid style=dotted,
    				]
    				\addplot[color=asl, mark=triangle] coordinates {
    					( 876 ,  0.05778989945401464 )
( 2784 ,  0.02591608203132607 )
( 3926 ,  0.0209251624880688 )
( 5124 ,  0.014265845034297345 )
( 6322 ,  0.010214023329476664 )
( 8332 ,  0.00594498515943203 )
( 9936 ,  0.004387528721117926 )
( 13738 ,  0.002874973135187766 )
( 17576 ,  0.0016356968024771451 )
( 21114 ,  0.0010737794030161348 )
( 26408 ,  0.0005351129490585203 )
( 29342 ,  0.00041658628250489674 )
( 32398 ,  0.0003568481092653749 )
( 44698 ,  0.0002485601233445473 )
( 63678 ,  0.00013641643031095773 )
( 75186 ,  8.904749876554364e-05 )
( 91588 ,  5.70304148184286e-05 )
( 108096 ,  3.897912099190602e-05 )
    				};
    				\addlegendentry{$p=3$}

    				\addplot[color=asb, mark=square] coordinates {
                        ( 2196 ,  0.023557524946529285 )
( 4808 ,  0.017382189197234083 )
( 7420 ,  0.01058371466726179 )
( 10032 ,  0.006180791456402256 )
( 12644 ,  0.0036791352369926877 )
( 15256 ,  0.0023302745537793703 )
( 17868 ,  0.0016540048954873335 )
( 21716 ,  0.0010991185216452477 )
( 27832 ,  0.000565142326975267 )
( 33088 ,  0.00032527247605145334 )
( 37076 ,  0.00023263959809698506 )
( 45100 ,  0.00012545959334450463 )
( 50152 ,  0.00012666469115214465 )
( 51560 ,  9.552103717618236e-05 )
( 62884 ,  5.43041795327268e-05 )
( 70656 ,  3.54730988384745e-05 )
( 81948 ,  2.0717392511489723e-05 )
( 90268 ,  1.5334604005133834e-05 )
( 108504 ,  7.677346787148695e-06 )
    				};
    				\addlegendentry{$p=5$}
    				
    				\addplot[color=blue, mark=pentagon] coordinates {
        ( 4092 ,  0.020919759234464874 )
( 9026 ,  0.013320993711578744 )
( 13960 ,  0.007778686315555429 )
( 18894 ,  0.004490540495678389 )
( 23828 ,  0.002616228518100854 )
( 28762 ,  0.0015581659090421038 )
( 33696 ,  0.0009707982131287379 )
( 38630 ,  0.0006602909801046882 )
( 45890 ,  0.0004305211027944902 )
( 55780 ,  0.00019036829776665267 )
( 60714 ,  0.00012244869253453145 )
( 65648 ,  8.673595958398042e-05 )
( 77212 ,  6.476504444087567e-05 )
( 85102 ,  3.778837543851916e-05 )
( 98252 ,  2.3220214176011403e-05 )
( 113076 ,  1.2415952797533033e-05 )
    				};
    				\addlegendentry{$p=7$}

            				\addplot[color=cyan, mark=diamond] coordinates {
                ( 584 ,  0.033402457180437356 )
( 2024 ,  0.01581617585793452 )
( 9548 ,  0.009613712410560162 )
( 37668 ,  0.005490375204312305 )
( 150496 ,  0.0032507896098540857 )
    				};
    				\addlegendentry{Uniform $p=5$}
        
    				\addplot[dashed,color=black, mark=none]
    				coordinates {
    					(6e+3, 5e-3)
    					(2e+4, 0.0002464751508773248)
    				};

                    \addplot[dashed,color=black, mark=none]
    				coordinates {
    					(4e+4, 5e-5)
    					(1e+5, 2.023857702507763e-06)
    				};

                    \addplot[dashed,color=black, mark=none]
    				coordinates {
    					(4e+4, 1e-3)
    					(1e+5, 0.00025298221281347036)
    				};

                    \addplot[dashed,color=black, mark=none]
    				coordinates {
    					(2e+4, 2e-2)
    					(1e+5, 0.010937454114084214)
    				};
    			\end{loglogaxis}
    			\draw (3.1,2.9) 
    			node[anchor=south]{$\mathcal{O}(h^{5})$};
                \draw (5.9,3.1) 
    			node[anchor=south]{$\mathcal{O}(h^{3})$};
                \draw (5.35, 0.6) 
    			node[anchor=south]{$\mathcal{O}(h^{7})$};
       \draw (5.35, 4.7) 
    			node[anchor=south]{$\mathcal{O}(h^{3/4})$};
    		\end{tikzpicture}